\documentclass[a4paper]{amsart}

\usepackage{amssymb,pstricks,amscd,epsfig}
\usepackage[totalwidth=17cm,totalheight=24cm]{geometry}
\usepackage{graphicx}
\usepackage{psfrag}

\makeatletter

\newenvironment{figurehere}
  {\def\@captype{figure}}
  {}
\makeatother

\begin{document}
\newtheorem{cor}{Corollary}[section]
\newtheorem{theorem}[cor]{Theorem}
\newtheorem{prop}[cor]{Proposition}
\newtheorem{lemma}[cor]{Lemma}
\newtheorem{sublemma}[cor]{Sublemma}
\theoremstyle{definition}
\newtheorem{defi}[cor]{Definition}
\theoremstyle{remark}
\newtheorem{remark}[cor]{Remark}
\newtheorem{example}[cor]{Example}
\newtheorem{question}[cor]{Question}

\newcommand{\cD}{{\mathcal D}}
\newcommand{\cF}{{\mathcal F}}
\newcommand{\cL}{{\mathcal L}}
\newcommand{\cM}{{\mathcal M}}
\newcommand{\cN}{{\mathcal N}}
\newcommand{\cP}{{\mathcal P}}
\newcommand{\cT}{{\mathcal T}}
\newcommand{\cCP}{{\mathcal C\mathcal P}}
\newcommand{\cML}{{\mathcal M\mathcal L}}
\newcommand{\cFML}{{\mathcal F\mathcal M\mathcal L}}
\newcommand{\cGH}{{\mathcal G\mathcal H}}
\newcommand{\cQF}{{\mathcal Q\mathcal F}}
\newcommand{\cQS}{{\mathcal Q\mathcal S}}
\newcommand{\C}{{\mathbb C}}
\newcommand{\HH}{{\mathbb H}}
\newcommand{\PP}{\mathbb{P}}
\newcommand{\N}{{\mathbb N}}
\newcommand{\R}{{\mathbb R}}
\newcommand{\Z}{{\mathbb Z}}
\newcommand{\Kt}{\tilde{K}}
\newcommand{\Mt}{\tilde{M}}
\newcommand{\dr}{{\partial}}
\newcommand{\betab}{\overline{\beta}}
\newcommand{\kappab}{\overline{\kappa}}
\newcommand{\pib}{\overline{\pi}}
\newcommand{\taub}{\overline{\tau}}
\newcommand{\ub}{\overline{u}}
\newcommand{\Sigmab}{\overline{\Sigma}}
\newcommand{\gd}{\dot{g}}
\newcommand{\diff}{\mbox{Diff}}
\newcommand{\dev}{\mbox{dev}}
\newcommand{\devb}{\overline{\mbox{dev}}}
\newcommand{\devt}{\tilde{\mbox{dev}}}
\newcommand{\vol}{\mbox{Vol}}
\newcommand{\hess}{\mbox{Hess}}
\newcommand{\cb}{\overline{c}}
\newcommand{\db}{\overline{\partial}}
\newcommand{\Sigmat}{\tilde{\Sigma}}
\newcommand{\Hyp}{\mathbb{H}}
\newcommand{\AdS}{\mathbb{A}\mathrm{d}\mathbb{S}}
\newcommand{\dS}{\mathrm{d}\mathbb{S}}
\newcommand{\cbull}{\bullet}
\newcommand{\psl}{\mathfrak{sl}}
\newcommand{\SL}{\mathrm{SL}}
\newcommand{\PSL}{\mathrm{PSL}}
\newcommand{\dual}{\star}

\newcommand{\cunc}{{\mathcal C}^\infty_c}
\newcommand{\cun}{{\mathcal C}^\infty}
\newcommand{\dd}{d_D}
\newcommand{\dmin}{d_{\mathrm{min}}}
\newcommand{\dmax}{d_{\mathrm{max}}}
\newcommand{\Dom}{\mathrm{Dom}}
\newcommand{\dn}{d_\nabla}
\newcommand{\ded}{\delta_D}
\newcommand{\delmin}{\delta_{\mathrm{min}}}
\newcommand{\delmax}{\delta_{\mathrm{max}}}
\newcommand{\hmin}{H_{\mathrm{min}}}
\newcommand{\maxi}{\mathrm{max}}
\newcommand{\oL}{\overline{L}}
\newcommand{\oP}{{\overline{P}}}
\newcommand{\xb}{{\overline{x}}}
\newcommand{\yb}{{\overline{y}}}
\newcommand{\Ran}{\mathrm{Ran}}
\newcommand{\tgamma}{\tilde{\gamma}}
\newcommand{\gammab}{\overline{\gamma}}
\newcommand{\cotan}{\mbox{cotan}}
\newcommand{\area}{\mbox{Area}}
\newcommand{\lambdat}{\tilde\lambda}
\newcommand{\xt}{\tilde x}
\newcommand{\Ct}{\tilde C}
\newcommand{\St}{\tilde S}
\newcommand{\tr}{\mbox{\rm tr}}
\newcommand{\tgh}{\mbox{th}}
\newcommand{\sh}{\mathrm{sinh}\,}
\newcommand{\ch}{\mathrm{cosh}\,}
\newcommand{\grad}{\mbox{grad}}
\newcommand{\hc}{H\C^3}
\newcommand{\cp}{\C P^1}
\newcommand{\rp}{\R P^1}

\newcommand{\II}{I\hspace{-0.1cm}I}
\newcommand{\III}{I\hspace{-0.1cm}I\hspace{-0.1cm}I}
\newcommand{\note}[1]{{\small {\color[rgb]{1,0,0} #1}}}

\title{A cyclic extension of the earthquake flow}
\author{Francesco Bonsante}
\address{Universit\`a degli Studi di Pavia\\
Via Ferrata, 1\\
27100 Pavia, Italy}
\email{bonsante@sns.it}
\thanks{F.B. is partially supported by the A.N.R. through project Geodycos }
\author{Gabriele Mondello}
\address{Universit\`a di Roma ``La Sapienza'' - Dipartimento di Matematica
``Guido Castelnuovo'' \\
piazzale Aldo Moro 5 \\
00185 Roma, Italy}
\email{mondello@mat.uniroma1.it}
\author{Jean-Marc Schlenker}
\thanks{J.-M. S. was partially supported by the A.N.R. through projects
ETTT, ANR-09-BLAN-0116-01, and GeomEinstein, ANR-09-BLAN-0116-01.}
\address{Institut de Math\'ematiques de Toulouse, UMR CNRS 5219 \\
Universit\'e Toulouse III \\
31062 Toulouse cedex 9, France}
\email{schlenker@math.univ-toulouse.fr}

\date{June 2011 (v1)}

\begin{abstract}
Let $\cT$ be Teichm\"uller space of a closed surface of genus at least
$2$. For any point $c\in \cT$, we describe an action of the circle on $\cT\times \cT$,
which limits to the earthquake flow when
one of the parameters goes to a measured lamination in the Thurston boundary of $\cT$. 
This circle action shares some of the main properties of the earthquake 
flow, for instance it satisfies an extension of Thurston's Earthquake Theorem
and it has a complex extension which is analogous
and limits to complex earthquakes. Moreover, a related circle action on 
$\cT\times \cT$ extends to the product of two copies of
the universal Teichm\"uller space.
\end{abstract}

\maketitle

\tableofcontents

\section{Introduction}

In all the paper we consider a closed, oriented surface $S$ of genus at least $2$. 
We denote by $\cT_S$, or simply by $\cT$, Teichm\"uller space of $S$, and by $\cML_S$, or simply by
$\cML$, the space of measured laminations on $S$.

\subsection{Earthquakes on hyperbolic surfaces}

Given a measured lamination $\lambda\in \cML_S$, we denote by
$E_\lambda$ the {\it left earthquake} along $\lambda$ on $S$.
$E_\lambda$ is a real-analytic map from $\cT_S$ to $\cT_S$, see
\cite{thurston-notes, kerckhoff,mcmullen:complex}. Recall that, in the simplest
case where $\lambda$ is supported on the simple closed
curve $\gamma$ with mass $a$, if $h\in \cT_S$ is a 
hyperbolic metric on $S$, $E_\lambda(h)$ is obtained by cutting
$(S,h)$ open along the minimizing geodesic homotopic to $\gamma$, 
rotating the left-hand side of $\gamma$ by $a$, and gluing back.

We consider here the {\it earthquake flow}, which can be defined as
a map:
$$
\begin{array}{cccc}
E: & \R\times \cT\times \cML & \rightarrow & \cT\times \cML \\
& (t,h,\lambda) & \mapsto & (E_{t\lambda}(h),\lambda) 
\end{array} $$
We call $E_t$ the corresponding map from $\cT\times \cML$
to $\cT\times \cML$, and will also use the notation 
$E_\lambda(h):=E_1(h,\lambda)$.

Earthquakes have a number of interesting properties, of which we
can single the following.
\begin{enumerate}
\item The earthquake flow defined above is indeed a flow: for all
$s,t\in \R$, $E_s\circ E_t=E_{s+t}$.
\item Thurston's Earthquake Theorem (see 
\cite{kerckhoff}): for any $h,h'\in \cT$, there is a unique $\lambda\in \cML$
such that $E_\lambda(h)=h'$.
\item For fixed $\lambda\in \cML_S$ and $h\in \cT_S$, the map 
$$ 
\begin{array}{rcl}
\R & \rightarrow & \cT_S \\
t & \mapsto & E_{-t\lambda}(h) 
\end{array}
$$ 
extends to a holomorphic map on a simply connected domain in 
$\C$ containing all complex 
numbers with non-negative imaginary part, see \cite{mcmullen:complex}. This
defines the notion of ``complex earthquake''.
\item When considered on imaginary numbers, complex earthquakes correspond to 
{\it grafting} maps.
The conformal grafting map $gr:\R_{\geq 0}\times \cT\times \cML\rightarrow \cT$
is actually obtained by composing the projective grafting map
$$ Gr:\R_{\geq 0}\times \cT\times \cML\rightarrow \cP~, $$
where $\cP$ is the space of complex projective
structures on $S$, with the forgetful
map $\cP\rightarrow \cT$ sending a $\mathbb{CP}^1$-structure to the 
underlying complex structure. 
\item Thurston (see \cite{kulkarni-pinkall}) proved that, for all $s>0$, $Gr_s$ 
it is a homeomorphism from $\cT\times \cML$ to $\cP$. 
\end{enumerate}

We introduce a flow on Teichm\"uller space
which in a way extends the earthquake flow, and which shares the properties described
above. The corresponding deformations are ``smoother'' than earthquakes, but 
earthquakes are limits in a natural sense. This motivates the term ``landslide''
that we use here. This deformation depends not on a measured
lamination but rather on a hyperbolic metric $h^\dual\in \cT$ and it determines
a flow:
$$
\begin{array}{rrcl}
  L:& \cT\times \cT\times S^1 & \rightarrow & \cT \times \cT \\
  & (h,h^\dual,e^{i\theta}) & \mapsto & L_{e^{i\theta}, h^\dual}(h)~. 
\end{array}
$$
We denote by $L_{e^{i\theta}}:\cT\times \cT\rightarrow \cT\times \cT$ 
the corresponding map seen as depending on the parameter $e^{i\theta}$.

We will also use the notation $L^1$ for the composition of $L$ with the projection
on the first factor, so that $L^1$ is a map from $S^1\times \cT\times \cT$ to $\cT$.
$L^1_{e^{i\theta}}$ will denote the same map, considered as depending on the 
parameter $e^{i\theta}\in S^1$, so it is a map from $\cT\times \cT$ to $\cT$.
Thus $L^1_{e^{i\theta}}$ is the analog of the earthquake map $E:\cT\times \cML\rightarrow
\cT$.

When the metrics $h$ and  $h^\dual$ are fixed, we can consider the image
of the map $S^1\ni e^{i\theta}\mapsto L^1_{e^{i\theta}}(h,h^\dual) \in \cT$
as a circle in Teichm\"uller space in which $h$ and $h^\dual$ are antipodal points.
When $h^\dual$ converges (projectively) to a measured lamination $\lambda$
at Thurston boundary of Teichm\"uller space, such a circle converges to 
the earthquake line $E_{t\lambda}(h)$. 
A more precise statement, Theorem \ref{tm:limit}, can be found below. 

This ``landslide flow'' shares the main properties of the earthquake flow 
recalled above: 
\begin{enumerate}
\item $L$ is a flow on $\cT\times \cT$ -- depending on the definition, checking this can be 
non trivial, see Theorem \ref{tm:cyclic}.
\item We prove an analog of Thurston's Earthquake Theorem, see Theorem \ref{tm:earthquake}.
\item For fixed $h,h^\dual\in \cT$, the map $L_\cbull(h,h^\dual):S^1\rightarrow\mathcal T$ 
extends to a holomorphic map from the
closed unit disk $\overline{\Delta}$ to $\cT$, see Theorem \ref{tm:complex}. This defines the ``complex landslide''
which are analogs of the ``complex earthquakes''.
\item For $r\in(0,1)$, the complex landslide $L_r$ corresponds to what we
call here ``smooth grafting'', which is analog to grafting in our context and we denote by
$sgr_{r}:\cT\times \cT\rightarrow \cT$ the map defined as $L_r:\cT\times \cT\rightarrow \cT\times \cT$ 
followed by projection on the first factor. 
It is obtained by composing a map $SGr:(0,1)\times \cT\times \cT 
\rightarrow \cP$ with the natural projection from $\cP$ to $\cT$.
\item For all $r\in (0,1)$, the map $SGr(r,\cbull,\cbull):(0,1)\times \cT\times \cT 
\rightarrow \cP$ is a homeomorphism.
\end{enumerate}

Our notations mean that we parameterize the complex landslides 
using the unit disk in $\C$ rather than the upper half-plane as is customary for complex earthquakes.
This notation is clearly equivalent but using the disk appears more 
natural in the context of the landslides considered here.

Considered as a circle action on $\cT\times\cT$,
the flow $L$ extends to a circle action
on the universal Teichm\"uller space, see Section \ref{sc:universal}.

\subsection{Harmonic maps and the landslide flow}
\label{ssc:harmonic}

Consider two hyperbolic metrics $c$ and $h$ on $S$. A map
$f:(S,c)\rightarrow (S,h)$ is said to be {\it harmonic} if it is a critical
point of the energy $E$.
The energy considered here is:
$$ E(f) = \frac{1}{2}\int_S \| df\|^2 \omega_c $$
where $\omega_c$ is the area element of $(S,c)$.
Although it is not immediately apparent in this definition, this notion
of harmonicity is conformally invariant on the domain, so that
we can regard $c$ as a conformal structure on $S$ rather than a metric.

\begin{theorem}[Sampson \cite{sampson}, Schoen and Yau \cite{schoen-Yau:78}] 
\label{tm:harmonic}
Let $c$ be a conformal class on $S$, and let $h\in \cT$ be a hyperbolic
metric. There is a unique harmonic map $f:(S,c)\rightarrow (S,h)$ 
isotopic to the identity. Moreover, $f$ is a diffeomorphism.
\end{theorem}

Consider a $C^1$ map $f:(S,h^\#)\rightarrow (S,h)$,
where $h^\#$ is a metric in the conformal class of $c$.
The Hopf differential $\Phi(f)$ of $f$ is a quadratic differential that measures
the traceless part of the pull-back of $h$ by $f$
and it is defined by the formula
$$ f^*h = e h^\# + \Phi + \overline{\Phi}~, $$
where $e=\frac{1}{2}tr_{h^\#}(f^*h)$.
If $f$ is harmonic, then $\Phi$ is holomorphic.
For $f$ $C^2$, also the converse holds.
It follows from its definition that $\Phi(f)$ is invariant under
conformal changes of the metric $h^\#$ on $S$.  

Conversely, given a holomorphic quadratic differential $\Phi$ on $(S,c)$,
there exists a unique hyperbolic metric $h$ on $S$ such that
the identity map $(S,c)\rightarrow (S,h)$ is harmonic with
Hopf differential $\Phi$, see \cite{sampson,wolf:teichmuller}.

This leads to the definition of a flow on $\cT$ depending on a ``center''
$c\in \cT$.

\begin{defi}
Let $c, h\in \cT$ and let $e^{i\theta}\in S^1$. We define $R_{c,e^{i\theta}}(h)$
as the (unique) hyperbolic metric $h'$ on $S$ such that, if
$f:(S,c)\rightarrow (S,h)$ and
$f':(S,c)\rightarrow(S,h')$ are the harmonic maps isotopic to the identity,
then $$ \Phi(f')=e^{i\theta} \Phi(f)~. $$
\end{defi}

This simple definition is strongly related to the flow $L$ mentioned above, 
but the relation is not obvious (see Corollary~\ref{cr:CD}), and using directly the definition of $R$ 
given here is not convenient.  
For this reason we give below a different definition of $L$, which
is more geometric, less directly accessible, but leads to straightforward 
arguments.

There is another, superficially similar flow on Teichm\"uller space,
the elliptic flow defined by one of us (Mondello), see \cite{mondello}. 
There are only limited similarities between the two flows, as 
should be clear from the sequel. 

\subsection{Minimal Lagrangian maps between hyperbolic surfaces}
\label{ssc:minilag}

The constructions considered here depend strongly on the notion of
minimal Lagrangian map between hyperbolic surfaces. Recall that, given
two hyperbolic metrics $h$ and $h^\dual$ on $S$, a diffeomorphism 
$m:(S,h)\rightarrow (S,h^\dual)$ is {\it minimal Lagrangian} if:
\begin{itemize}
\item it is area-preserving and orientation-preserving,
\item its graph is a minimal surface in $(S\times S, h\oplus h^\dual)$. 
\end{itemize}
 
\begin{theorem}[Schoen \cite{schoen:role}, Labourie \cite{L5}] \label{tm:LS}
Let $h, h^\dual$ be two hyperbolic metrics on $S$. There exists a unique 
minimal Lagrangian diffeomorphism $m:(S,h)\rightarrow (S,h^\dual)$ isotopic to the
identity. 
\end{theorem}

Minimal Lagrangian maps actually have a description in terms of hyperbolic surfaces
only, as follows (see e.g. \cite{L5}).

\begin{prop} \label{pr:minlag}
If $m:(S,h)\rightarrow (S,h^\dual)$ is minimal Lagrangian, then $m^*h^\dual=h(b\cbull,b\cbull)$,
where $b:TS\rightarrow TS$:
\begin{enumerate}
\item is self-adjoint for $h$,
\item has determinant $1$,
\item satisfies the Codazzi equation: $d^\nabla b=0$, where $\nabla$ is the Levi-Civita
connection of $h$. 
\end{enumerate}
Conversely, if $m:S\rightarrow S$ is a diffeomorphism satisfying those properties, 
then it is minimal Lagrangian. 
\end{prop}

\begin{cor} \label{cr:b}
Let $h,h^\dual$ be two hyperbolic metrics on $S$. There exists a unique
bundle morphism $b:TS\rightarrow TS$ which is self-adjoint for $h$, of determinant equal to
$1$ everywhere, satisfies Codazzi equation $d^\nabla b=0$, where $\nabla$ is the
Levi-Civita connection of $h$, and such that $h(b\cbull,b\cbull)$ is isotopic to $h^\dual$.    
\end{cor}

A consequence of this proposition is that for any $\tau,\tau^\dual\in \cT$,
we can realize $\tau$ and $\tau^\dual$
as a pair of hyperbolic metrics $h$ and $h^\dual$ ({\it not} considered up to isotopies)
so that $h^\dual=h(b\cbull, b\cbull)$, where $b$ is self-adjoint for $h$, of determinant $1$,
and satisfies Codazzi equation $d^\nabla b=0$. A pair of metrics with this property
will be a called a {\it normalized representative} of $(\tau,\tau^\dual)$. Clearly a normalized
representative of $(\tau,\tau^\dual)$ is uniquely determined up to isotopies acting
diagonally on both $h$ and on $h^\dual$.

By abuse of notation, we will sometimes denote by $(h,h^\dual)$ both a couple of normalized
hyperbolic metrics and its correspondent point in $\cT\times\cT$.

\subsection{The landslide action on $\cT\times \cT$}

We now introduce the action $L$ of $S^1$ on $\cT\times \cT$. We will see below that it is
strongly related to the map $R$ introduced above.

\begin{defi}\label{ldsl:defi}
Let $h,h^\dual$ be two hyperbolic metrics on $S$, and let $\theta\in \R$. We consider the
bundle morphism $b:TS\rightarrow TS$ given by Corollary \ref{cr:b}, and set
\begin{equation}\label{operator:eq}
 \beta_\theta := \cos(\theta/2) E + \sin(\theta/2) Jb~, 
 \end{equation}
where $E:TS\rightarrow TS$ is the identity map and  $J$ is the complex structure 
of $h$ on $S$. We then call 
$$ L_{e^{i\theta}}(h,h^\dual) := (h(\beta_\theta\cbull, \beta_\theta\cbull),
h(\beta_{\theta+\pi}\cbull, \beta_{\theta+\pi}\cbull))~. $$
\end{defi}

Notice that by construction, $L_1(h,h^\dual)=(h,h^\dual)$, while $L_{-1}(h,h^\dual)=(h^\dual,h)$. 

\begin{prop} \label{pr:basics}
For all $\theta\in \R$, $h(\beta_\theta\cbull, \beta_\theta\cbull)$ is a hyperbolic
metric on $S$.
\end{prop}

\begin{theorem} \label{tm:cyclic}
Let $h,h^\dual$ be two hyperbolic metrics on $S$, let $\theta,\theta'\in \R$. 
Then 
$$ L_{e^{i\theta'}}(L_{e^{i\theta}}(h,h^\dual))=L_{e^{i(\theta'+\theta)}}(h,h^\dual)~. $$
In other terms, $L$ defines an action of $S^1$ on $\cT\times \cT$.   
We call $L$ the {\bf landslide flow}, or {\bf landslide action} on $\cT\times \cT$.
\end{theorem}

The proofs of Proposition \ref{pr:basics} and of Theorem \ref{tm:cyclic} are in Section \ref{ssc:cyclic}.

\subsection{Relations to AdS geometry}

We briefly recall some properties
of globally hyperbolic anti-de Sitter manifolds. More details can be found
e.g. in \cite{mess,mess-notes}.

The anti-de
Sitter space is a Lorentz analog of hyperbolic 3-space, it can be defined
as the quadric:
$$ \AdS^3 = \{ x\in \R^{2,2}~|~ \langle x,x\rangle =-1\}~, $$
where $\R^{2,2}=(\R^4,-dx_0^2-dx_1^2+dx_2^2+dx_3^2)$.
It is a complete Lorentz manifold of constant 
curvature $-1$ with fundamental group isomorphic to $\Z$. 

A manifold $N$ with an AdS metric -- a Lorentz metric locally modeled
on $\AdS^3$ -- is {\it maximal globally hyperbolic} (MGH) if:
\begin{itemize}
\item $N$ contains a closed space-like surface $F$,
\item any inextendible time-like curve in $N$ intersects $F$ exactly once,
\item $N$ is maximal for inclusion, under these properties.
\end{itemize}
Mess \cite{mess,mess-notes} proved that,
if $N$ is (GH) and $\bar\phi:S\rightarrow N$
is an embedding onto a closed space-like surface $F$, then $N$ is
the quotient of a convex domain $\Omega$ in $\AdS^3$ by an action
of the fundamental group of $S$.

A key feature of $\AdS^3$ is that the identity component of its isometry
group is isomorphic to $\SL_2(\R)\times \SL_2(\R)/\Z_2$, which is the double
cover of $\PSL_2(\R)\times \PSL_2(\R)$. As a consequence,  the action of 
$\pi_1(S)$ on $\Omega$ decomposes as $(\rho_l,\rho_r)$, where $\rho_l$
and $\rho_r$ are morphisms from $\pi_1(S)$ to $\PSL_2(\R)$. It was proved
in \cite{mess} that these morphisms have maximal Euler number, so that
they correspond to points in the Teichm\"uller space of $S$. Maximal
globally hyperbolic AdS spaces are uniquely determined by these left and
right representations, see \cite{mess,mess-notes}.

\begin{lemma} \label{lm:ads}
Let $h,h^\dual$ be  a pair of normalized metrics, 
let $\theta\in (0,\pi)$. There exists a unique
equivariant embedding $(\phi,\rho)$ of $\tilde{S}$ in
$\AdS^3$ with induced metric 
$\cos^2(\theta/2)h$ and third fundamental form $\sin^2(\theta/2)h^\dual$.
Moreover, $\rho$ is the holonomy representation of a globally hyperbolic
AdS manifold $N$: the first factor in $L_{e^{i\theta}}(h,h^\dual)$ is the
left representation of $N$ and the first factor in  $L_{e^{-i\theta}}(h,h^\dual)$
is the right representation of $N$.
\end{lemma}

The proof is in Section \ref{ssc:ads}.

\subsection{The center of circles in Teichm\"uller space}
\label{ssc:center}

Let $h,h^\dual\in \cT$. For each $\theta\in \R$, let $(h_\theta,h^\dual_\theta)=
L_{e^{i\theta}}(h,h^\dual)$. According to Theorem \ref{tm:LS}, there is a unique
minimal Lagrangian diffeomorphism $m_\theta:(S,h_\theta)\rightarrow (S,h^\dual_\theta)$
isotopic to the identity. We can then consider on $S$ the conformal structure
$c_\theta$ of the metric $h_\theta + m_\theta^* h^\dual_\theta$. We call $c_\theta$ the
{\it center} of $(h_\theta,h^\dual_\theta)$. 
This conformal class of metrics has some interesting properties, proved in Section \ref{ssc:c}.

\begin{theorem} \label{thm:c}
\begin{enumerate}
\item $m_\theta$ is the identity -- that is, the identity is minimal Lagrangian between
$(S,h_\theta)$ and $(S, h^\dual_\theta)$. 
\item $c_\theta$ does not depend on $\theta$, it is equal to a fixed conformal class $c$. 
\item Let $f_\theta:(S,c)\rightarrow (S,h_\theta)$ and $f_\theta^\dual:(S,c)\rightarrow (S,h_\theta^\dual)$ be the
unique harmonic maps isotopic to the identity. Then $f_\theta$ and $f_\theta^\dual$
have opposite Hopf differentials. 
\item For any $\theta\in \R$, 
$$ \Phi(f_\theta) = e^{i\theta} \Phi(f_0)~. $$
\end{enumerate}
\end{theorem}

\subsection{Obtaining $R$ from $L$}

As a consequence of Theorem \ref{thm:c}, we find a simple relation between the
map $R$ defined earlier in terms of Hopf differential, and the map $L$.

\begin{cor} \label{cr:CD}
Let $(h,h^\dual)$ be a couple of normalized metrics
and let $c$ be the conformal class of $h+h^\dual$. For any
$e^{i\theta}\in S^1$, we have
$$ L_{e^{i\theta}}(h,h^\dual) = (R_{c,\theta}(h),R_{c,\theta+\pi}(h^\dual))~. $$
\end{cor}

\subsection{The earthquake flow as a limit}

 
\begin{theorem} \label{tm:limit}
Let $h\in \cT$, let $(h^\dual_n)_{n\in \N}$ be a sequence of hyperbolic metrics
and let $\lambda\neq 0$ be a measured lamination.
Consider a sequence $(\theta_n)_{n\in \N}$
of positive real numbers such
that $\lim_{n\rightarrow \infty} \theta_n \ell_{h^\dual_n}=\iota(\lambda,\cbull)$
in the sense of convergence of the length spectra of simple closed curves.
Then
$$ \lim_{n\rightarrow +\infty} h^1_n = E_{\lambda/2}(h)\qquad
\qquad
\lim_{n\rightarrow+\infty}\theta_n\ell_{h^2_n}=\iota(\lambda,\cbull) $$  
where $(h^1_n,h^2_n):=L_{e^{i\theta_n}}(h,h_n^\dual)$.
\end{theorem}


At first sight it would appear more natural to take the limit 
where the sequence of centers $(c_n)$ converges projectively
to $\lambda$. However the statement obtained by replacing $L$ by $R$ and
$h_n^\dual$ by $c_n$ in Theorem \ref{tm:limit} turns out to be false, as
proved -- in one example -- in Section \ref{ssc:limit_cn}, see Corollary
\ref{cor:center}.

The heuristic argument motivating Theorem \ref{tm:limit} involves 
the convergence of constant Gauss curvature surfaces to a pleated surface
in $\AdS^3$. However, writing a proof based on these ideas turns out to be
more difficult than it appears. A key technical statement is that 
minimal Lagrangian maps have a close proximity to Thurston compactification
of Teichm\"uller space: minimal Lagrangian maps isotopic to the identity
provide ``the'' correct normalization to understand the convergence of
a sequence of hyperbolic metrics to a projective measured lamination in 
Thurston boundary of $\cT$.

\begin{theorem} \label{tm:boundary}
Let $h$ be a hyperbolic metric on $S$, and let $(h^\dual_n)_{n\in\mathbb{N}}$
be a sequence of
hyperbolic metrics such that
$\theta_n \ell_{h_n^\dual}\rightarrow \iota(\lambda,\cbull)$, where $\lambda$
is a measured geodesic lamination, the $\theta_n$ are positive numbers, 
and the convergence is in the sense of the length spectrum. 
For each $n$, let $m_n:(S,h)\rightarrow (S,h^\dual_n)$ be the
minimal Lagrangian diffeomorphism isotopic to the identity. Then, for
every smooth arc
$\alpha$ in $S$ that meets the $h$-geodesic representative of $\lambda$
transversely and with endpoints not in the support of
(the $h$-geodesic representative of) $\lambda$, the
length for $\theta_n^2 m_n^*(h^\dual_n)$ of the geodesic segment homotopic to $\alpha$
(with fixed endpoints)
converges to the intersection between $\alpha$ and $\lambda$.
\end{theorem}

The proof of this theorem involves the convergence of smooth surfaces
to a pleated limit, but in the hyperbolic, rather than the anti-de Sitter, context.

\subsection{An extension of the Earthquake Theorem}

We can now state an extension to the landslide flow $L$ of Thurston's
Earthquake Theorem (see \cite{kerckhoff}). Recall that this theorem
states that, given two hyperbolic metrics $h$ and $h'$ on a surface,
there is a unique measured lamination $\lambda$ such that the left
earthquake along $\lambda$ sends $h$ to $h'$.

\begin{theorem} \label{tm:earthquake}
Let $h,h'\in \cT$ and let $e^{i\theta}\in S^1\setminus \{ 0\}$. 
There is a unique $h^\dual\in \cT$ such that $L^1_{e^{i\theta}}(h,h^\dual)=h'$.  
\end{theorem}

We give in section \ref{ssc:bbz} a simple proof based on a recent result
of Barbot, B\'eguin and Zeghib \cite{BBZ2} on the existence and uniqueness of
constant Gauss curvature foliations in globally hyperbolic AdS manifolds. 

As an easy consequence, a similar statement holds also for the flow $R$.

\subsection{A complex extension}

The earthquake flow has an extension as a map
$E: \overline{\HH}\times\cT\times \cML\rightarrow 
\cT$, where $\overline{\HH}$ is the set of complex numbers with nonnegative imaginary part. This 
map has the property that, for any $h\in \cT$ and any $\lambda\in \cML$, the map
$z\mapsto E(z,h,\lambda)$ is holomorphic, see \cite{mcmullen:complex}. It can be defined in terms
of grafting, or (for small $\lambda$) in terms of pleated surfaces in hyperbolic 3-space. 

In Section \ref{sc:complex}
we prove that the landslide map $L$ defined above has
a similar holomorphic extension
where the parameter $e^{i\theta}$ is replaced by a 
complex number $\zeta$ in the closed unit disk.
This defines many holomorphic disks in Teichm\"uller space,
see Theorem \ref{tm:complex}. 
Similarly to what happens for complex earthquakes,
this construction factors through the space of
complex projective structures on $S$ for $\zeta\neq 0$,
and the complex cyclic flow provides punctured
holomorphic disks in this space. This factorization however does
not extend for $\zeta=0$.

The complex landslide map limits to complex earthquakes just as the ``real''
landslide flow limits to the earthquake flow, see Theorem \ref{conv:thm}.

We hope at some point in the future to give another proof of the holomorphicity of this complex
landslide map, based on a geometric argument taking place in the complexification
of $\Hyp^3$. This line of argument should also provide a straightforward and geometric
way to understand why complex earthquakes are holomorphic disks. 

\subsection{Landslide on the universal Teichm\"uller space}

Recall that a homeomorphism of the circle is {\it quasi-symmetric} if and only 
if it is the boundary value of a quasi-conformal diffeomorphism of the disk. 

\begin{defi}
The universal Teichm\"uller space $\cT_U$ is the quotient of the group $\cQS$ of quasi-symmetric 
homeomorphisms of the circle by left composition by projective transformations. 
\end{defi}

The universal Teichm\"uller space contains embedded copies of the Teichm\"uller
space of all closed surfaces. Indeed, consider a closed surface $S$ of genus at 
least $2$, a fixed hyperbolic metric $h^\#$ on $S$, and its holonomy representation 
$\rho^\#:\pi_1(S)\rightarrow \PSL_2(\R)$. 
Given another hyperbolic metric $h$ on $S$ and its holonomy representation
$\rho:\pi_1(S)\rightarrow \PSL_2(\R)$, there is a quasiconformal map 
$\tilde{f}:\Hyp^2\rightarrow \Hyp^2$ conjugating $\rho^\#$ and $\rho$. Moreover, the 
boundary value $\partial \tilde{f}:\partial_\infty \Hyp^2\rightarrow \partial_\infty \Hyp^2$
is uniquely determined by $\rho^\#$ and $\rho$, and the map sending $h$ to 
$\partial \tilde{f}$ is an embedding of $\cT_S$ in $\cT_U$, see e.g. \cite{gardiner-harvey}.

Let $\psi:S^1\rightarrow S^1$ be a quasi-symmetric homeomorphism. 
There is (see \cite{maximal}) a unique
minimal Lagrangian quasiconformal diffeomorphism
$m:\Hyp^2\rightarrow \Hyp^2$
with $\partial m=\psi$.
As for closed surfaces, there is then a unique bundle morphism 
$b:T\Hyp^2\rightarrow T\Hyp^2$ such that 
\begin{itemize}
\item $b$ is self-adjoint,
\item it satisfies the Codazzi equation $d^\nabla b=0$,
\item $m^*g=g(b\cbull, b\cbull)$, where $g$
is the hyperbolic metric on $\Hyp^2$.
\end{itemize}
For every $\theta\in \R$ we then consider
$\beta_\theta:=\cos(\theta/2)E+\sin(\theta/2)b$,
where $E$ is the identity and $g_\theta:=
g(\beta_\theta \cbull,\beta_\theta \cbull)$.

\begin{lemma} \label{lm:univ}
$g_\theta$ is a complete hyperbolic metric on $\Hyp^2$.
The identity map between $(\Hyp^2,g)$ and $(\Hyp^2,g_\theta)$ 
is quasiconformal (and minimal Lagrangian), and its extension
$\psi_\theta:S^1\rightarrow S^1$
to the boundary of $\Hyp^2$
is quasi-symmetric, 
so that it defines a point in $\cQS$.
\end{lemma}

In Section 8 we show how to use this fact to construct
an extension of $L$ to a non-trivial circle action 
$\cL$ on $\cT_U\times \cT_U$ (see Theorem \ref{tm:univ}).

\subsection{Content of the paper}

In Section 2 we present the background material, concerning in particular minimal Lagrangian
diffeomorphisms between hyperbolic surfaces and globally hyperbolic AdS manifolds. In Section
3 we define the landslide flow and prove that it is indeed a flow (Theorem \ref{tm:cyclic})
as well as Theorem \ref{thm:c}. In Section 4 we give the proof of the extension to the landslide
flow of Thurston's Earthquake Theorem (Theorem \ref{tm:earthquake}). Then in Section 5 we construct
the complex landslide map, actually as a map from
$\cT\times \cT\times \dot{\overline{\Delta}}$ to $\cP$, where 
$\dot{\overline{\Delta}}$ is the pointed
closed unit disk in $\C$, and we prove that it is holomorphic
and that it extends over $\overline{\Delta}$ as a map to $\cT$. 
Section 6 considers the limit when the parameter $h^\dual$
converges projectively to a measured lamination at Thurston boundary of Teichm\"uller
space, and contains the proof of Theorem \ref{tm:limit} as well as its complex extension,
Theorem \ref{conv:thm}. In Section 7, on the other hand, we show that the situation is
not as simple for the ``center'' $c$ determined by a fixed metric $h$ and a sequence
$h_n^\dual$ going to a point at infinity in Thurston compactification of $\cT$: the limit
of the corresponding sequence of centers does not depend only on the limit of $(h_n^\dual)$.
Section 8 deals with
the circle action on the universal Teichm\"uller space, while Section 9 contains some
remarks and open questions.

\subsection*{Acknowledgment}

We are grateful to the Institute of Mathematical Sciences of the National University of
Singapore, where most of the results presented here were obtained.
The second named author would like to thank Mike Wolf for helpful
clarifications on his paper \cite{wolf:infinite}.

\section{Minimal lagrangian maps and AdS geometry}

We present in this section some background material used in the paper. 

\subsection{Notations}

In all the paper we consider a closed, oriented surface $S$ of genus
at least $2$. 

We consider $\AdS^3$, as well as all AdS manifolds, as oriented and time-oriented.
All the embeddings of $S$ that we consider will implicitly be considered as
time-oriented, that is, the oriented normal to the image is future-oriented.
Moreover, the convex embeddings will always be considered to be positively convex,
that is, the oriented normal is future-directed and
pointing towards the convex side.
We recall that it is possible to identify the isometries of $\AdS^3$ with double cover of
$\PSL_2(\R)\times\PSL_2(\R)$ in such a way that,
if $S$ is a positively convex pleated surface in $\AdS^3$, bent along $\lambda$
and with first fundamental form $h$,
then the first (resp. second) factor corresponds to the holonomy of the
hyperbolic surface obtained from $h$ performing a left (resp. right) earthquake
along $\lambda$.

\subsection{Hyperbolic ends}

The 3-dimensional hyperbolic space can be defined as a quadric in the 4-dimensional
Minkowski space $\R^{1,3}=(\R^4,-dx_0^2+dx_1^2+dx_2^2+dx_3^2)$,
with the induced metric.
$$ \Hyp^3 = \{ x\in \R^{1,3}~|~\langle x,x\rangle=-1~ \mbox{and}~ x_0>0\}~. $$
It is a simply connected, complete manifold with constant curvature $-1$.

A quasifuchsian hyperbolic manifold is a 3-dimensional manifold locally isometric
to $\Hyp^3$, homeomorphic to $S\times \R$, which contains a non-empty compact convex subset. 

Such a quasifuchsian manifold $M$ contains a smallest non-empty convex
subset $C(M)$ called its convex core. $M$ is {\it Fuchsian} if $C(M)$ is a 
totally geodesic surface, otherwise the boundary of $C(M)$ is the disjoint union
of two pleated surfaces. 

Each connected component of the complement of $C(M)$ in $M$ is an instance of
a {\it hyperbolic end}: a hyperbolic manifold homeomorphic to $S\times \R_{>0}$,
complete on one side and bounded by a locally concave pleated surface on the 
other. There is a one-to-one correspondence between hyperbolic ends homeomorphic
to $S\times \R_{>0}$ and complex projective structures on $S$, which associates
to a hyperbolic end the natural complex projective structure on its boundary
at infinity, see e.g. \cite{kulkarni-pinkall}.

Labourie \cite{L5} proved that any hyperbolic end has a unique foliation by
convex, constant curvature surfaces. The curvature varies monotonically from 
$-1$ close to the pleated surface boundary, to $0$ close to the boundary at
infinity. 

Given an oriented surface $\Sigma$ in a hyperbolic end $M$
(or in $\Hyp^3$) we will usually denote
by $I$ its induced metric, and by $B$ its shape operator, considered as a
bundle morphism from $T\Sigma$ to $T\Sigma$.
It is defined by $BX=\nabla_X\nu$,
where $\nu$ is the oriented unit normal to
$\Sigma$ and $\nabla$ is the Levi-Civita
connection of $M$. We will also denote by $E:TS\rightarrow TS$ the identity. 

\begin{defi} \label{df:grafted}
Let $\Sigma$ be a convex surface embedded in a hyperbolic end $M$ with
embedding data $(I_\Sigma,B_\Sigma)$. The {\it{grafted metric}}
on $\Sigma$ is $I^\#_\Sigma=I_\Sigma((E+B_\Sigma)\cbull,(E+B_\Sigma)\cbull)$.
\end{defi}

A basic and well-known property of this metric $I^\#_\Sigma$ is that the 
hyperbolic Gauss map -- sending a point $x\in \Sigma$ to the endpoint
at infinity of the geodesic ray starting at $x$ orthogonal to $\Sigma$ --
is a conformal map between $(\Sigma, I^\#_\Sigma)$ and $\partial_\infty M$
with its conformal structure.
More details will be found in Section \ref{sc:limit}.

\subsection{The duality between hyperbolic and de Sitter ends}

The 3-dimensional de Sitter space can be defined, as the hyperbolic space, as
a quadric in the 4-dimensional Minkowski space, with the induced metric. 
$$ \dS^3 = \{ x\in \R^{1,3}~|~\langle x,x\rangle=1\}~. $$
There is a one-to-one correspondence between points in $\dS^3$ and
oriented totally geodesic planes in $\Hyp^3$, see e.g. \cite{RH,shu}. 
Given an oriented surface $S\subset \Hyp^3$, its dual is the set $S^\dual$ of 
points of $\dS^3$ corresponding to oriented planes tangent to $S$
in $\Hyp^3$.  If $S$ is smooth and locally strictly convex, then $S^\dual$
is also smooth and locally strictly convex.

Consider a quasifuchsian hyperbolic manifold $M$, and let $E$ be one 
of the ends of $M$, that is, one of the connected components of the 
complement of $C(M)$ in $M$. The universal cover of $M$ is identified
with $\Hyp^3$, and the universal cover $\tilde{E}$ of $E$ 
is then identified with a connected component of the
complement of the convex hull of the limit set $\Lambda$ of $\pi_1(M)$
in $\partial_\infty \Hyp^3$. 
Let $\tilde{E}^\dual$ be the set of points of $\dS^3$ corresponding to
oriented planes in $\Hyp^3$
contained in $\tilde{E}$.
Then
(see \cite{mess}) $\pi_1(M)$ acts properly discontinuously
on $\tilde{E}^\dual$,
and the quotient is a de Sitter {\it domain of dependence}, that is, 
a globally hyperbolic maximal de Sitter manifold (see below for the
definition in the AdS case).

This construction actually extends (see \cite{mess}) to a hyperbolic end
$E$ which is not necessarily one of the ends of a quasifuchsian manifold.
In this manner, any hyperbolic end $E$ has a ``dual'' de Sitter domain of
dependence $E^\dual$, and conversely. 

One feature of this duality which will be used below is that if $S$
is a surface in $E$ with constant curvature $K$, then there is a dual
surface $S^\dual$ in $E^\dual$. (It is the quotient by $\pi_1(M)$ of the surface
in $\dS^3$ dual to the universal cover of $S$ in $\Hyp^3$.) The curvature 
of $S^\dual$ is then also constant, and equal to $K/(K+1)$. In this manner
a foliation of $E$ by constant curvature surfaces gives rise to a dual
foliation of $E^\dual$ by constant curvature surfaces (see \cite{BBZ2} for
more details).

\subsection{Globally hyperbolic AdS manifolds}
\label{ssc:ads}

The definition of $\AdS^3$ and of globally hyperbolic AdS manifolds are recalled in
the introduction. There are many similarities between quasifuchsian hyperbolic
manifolds and globally hyperbolic AdS manifolds, some of which -- being used
in the arguments below -- are recalled here.

Let $N$ be a globally hyperbolic AdS 3-manifold, and let $F$ be a closed,
space-like surface in $N$ for which the induced metric has negative sectional
curvature (or, equivalently, the determinant of the second fundamental form
of $F$ is everywhere larger than $-1$). 
Let $I$ and $\nu$ be the induced metric
on this surface and a unit normal vector.
Let $J$ be the complex structure induced by $\nu$ on $F$: namely
$J(v)=\nu\times v$ where $\times $ is the vector product on $T\AdS^3$.
Finally, let $B=\nabla\nu$ be the the shape operator of $F$,
where $\nabla$ is the the Levi-Civita connection of $\AdS^3$. We consider
the Riemannian metrics on $F$ 
$$ h_l = I((E+JB)\cbull,(E+JB)\cbull)~, h_r = I((E-JB)\cbull,(E-JB)\cbull)~. $$
Then $h_l$ and $h_r$ are two smooth hyperbolic metrics on $F$
(see \cite{minsurf}).
This can be used when $F$
is a maximal or constant mean curvature surface in $N$, 
but also when $F$ is a constant Gauss curvature surface.
\begin{remark}
Notice that, even if $J$ and $B$ depend on the choice
of a normal vector, $JB$ and the metrics $h_l$ and $h_r$ are independent of it.

According to our orientation and time-orientation
of $\AdS^3$, the holonomy of the metric $h_l$ is equal to the first factor
of the holonomy of $N$ and the holonomy of $h_r$
is equal to the second factor of the holonomy of $N$, see \cite{minsurf}. 
This last observation can be used to prove Lemma \ref{lm:ads}.
\end{remark}

\begin{proof}[Proof of Lemma~\ref{lm:ads}]
Let $b:TS\rightarrow TS$ be the positive $h$-self-adjoint operator
given by Corollary \ref{cr:b}
such that $h^\dual=h(b\cbull, b\cbull)$.
Since $h$ and $h^\dual$ are normalized metrics we know that $b$ is
a solution of Codazzi equation and $\det b=1$.

We consider now the pair
\[
I_{\theta}=\cos^2(\theta/2)h\qquad B_\theta=\tan(\theta/2)b\,.
\]
Clearly, $B_\theta$ is a solution of Codazzi equation. Moreover, 
since $I_\theta$ is a metric of constant curvature $K=-\frac{1}{\cos^2(\theta/2)}$,
we easily get that
\[
  K_\theta=-1-\det B_\theta
\]
so $(I_\theta, B_\theta)$ is also a solution of Gauss equation for 
spacelike AdS surfaces (see \cite{minsurf}).

This implies that there is an equivariant map
\[
\phi:\tilde S\rightarrow \tilde{F}\subset \AdS^3
\]
whose embedding data are $\tilde I_\theta$ and $\tilde B_\theta$.
The map $\phi$ is unique up to isometries of $\AdS^3$.
We will also require that
\begin{itemize}
\item the normal field $\tilde \nu$
that induces the right orientation on $\tilde S$ points toward the convex side
\item $\tilde \nu$ is a future-directed vector field.
\end{itemize}

\end{proof}

Lemma \ref{lm:ads} and \cite{minsurf} imply that $L_{e^{i\theta}}(h,h^\dual)$
is a couple of hyperbolic metrics. For convenience of the reader,
we will give a simple proof of this fact in Section \ref{ssc:civita}.\\

Globally hyperbolic AdS manifolds have a unique foliation by constant mean
curvature surfaces, see \cite{BBZ}.

A globally hyperbolic AdS manifold $N$ contains a smallest non-empty
convex subset $C(N)$, called its convex core: $N$ is called {\it Fuchsian} if
$C(N)$ is a totally geodesic surface;
otherwise, the boundary of $C(N)$ is the
disjoint union of two pleated locally convex surfaces in $N$, so that 
its induced metric is hyperbolic and its pleating is described by a measured
lamination (see \cite{mess}). The complement of $C(N)$ in $N$ has two connected
components, one future convex, the other past convex. Barbot, B\'eguin and
Zeghib \cite{BBZ2} proved that $N\setminus C(N)$ has a unique foliation by
convex, constant Gauss curvature surfaces. The curvature is monotonic along the
foliation, and varies from $-1$ in the neighborhood of the convex core to 
$-\infty$ near the initial/final singularity.

\subsection{Minimal Lagrangian maps between hyperbolic surfaces}
\label{ssc:24}

The definition of minimal Lagrangian diffeomorphisms has been recalled in the
introduction. Remark that the definition directly shows that the inverse of
a minimal Lagrangian diffeomorphism is also minimal Lagrangian.

Let us mention here that they occur in several distinct geometric
contexts, and the interplay between the different occurences is used below, 
in particular in Section 4. 
\begin{itemize}
\item If $S$ is a surface of constant curvature $K$ in a constant curvature, Riemannian or Lorentzian
$3$-manifold $M$, then the third fundamental form $\III$ of $S$ also has constant curvature $K^\dual$,
where $K^\dual$ depends on $K$, on the ambient curvature, and on whether the ambient space is
Riemannian or Lorentzian (for instance, if $M=\Hyp^3$, then $K^\dual=K/(K+1)$). If both $K$ and $K^\dual$
are negative, then $|K| I$ and $|K^\dual| \III$ are hyperbolic metrics, and the identity, considered
as a map from $(S,|K| I)$ to $(S, |K^\dual| \III)$, is minimal Lagrangian.
\item If $M$ is an ``almost Fuchsian'' manifold, that is, $M$ is a quasifuchsian hyperbolic 3-manifolds
containing a closed, embedded minimal surface $S$ with principal curvatures everywhere in $(-1,1)$, then
$S$ is the unique closed minimal surface in $M$. The ``hyperbolic Gauss maps'' send $S$ to each
connected component of the boundary at infinity of $M$, and both maps are diffeomorphisms. Composing
these maps one finds a diffeomorphism between one component of $\dr_\infty M$ and the other. This 
diffeomorphism is minimal Lagrangian if each boundary component is endowed with the (unique) 
hyperbolic metric in its conformal class. (See e.g. \cite{minsurf} for details and proofs.)
\item Similarly, if $N$ is a globally hyperbolic AdS manifold, then it contains a unique
closed, space-like maximal surface $F$.
Consider the metrics $h_l$ and $h_r$ defined above on $F$.
Then $h_l$ and $h_r$ are the left and right hyperbolic metrics of $N$, respectively, and moreover
the identity between $(F,h_l)$ and $(F, h_r)$ is minimal Lagrangian (see \cite{minsurf} for details).
\end{itemize}
It is the first of these occurences which will play the largest role here.

Minimal Lagrangian maps between hyperbolic surfaces are intimately related to harmonic maps: let
$m:(S,h)\rightarrow (S',h')$ be a minimal Lagrangian diffeomorphism between two hyperbolic surfaces
and consider the conformal structure $c$ on $S$ of the metric $h+m^*h'$. Then
\begin{itemize}
\item the identity is harmonic between $(S,c)$ and $(S,h)$,
\item $m$ is harmonic between $(S,c)$ and $(S', h')$,
\item those two harmonic maps have opposite Hopf differentials. 
\end{itemize}
The converse is also true. Details can be found e.g. in \cite{schoen:role}.

\section{Definition of the cyclic flow}\label{section:definition}

In this section we consider two fixed normalized hyperbolic metrics $h,h^\dual$ on 
$S$, and call $b$ the bundle morphism given by Corollary \ref{cr:b}.
Let $\beta_\theta$ be the family of operators defined in (\ref{operator:eq}).

\begin{defi}
Given $\theta\in \R$, we call
$$ h_\theta = h(\beta_\theta\cbull,\beta_\theta\cbull)~. $$
\end{defi}
Comparing with definition \ref{ldsl:defi}, 
we have $L_{e^{i\theta}}(h,h^\dual)=(h_\theta, h_{\theta+\pi})$.

\subsection{The Levi-Civita connection of $h_\theta$}\label{ssc:civita}

\begin{lemma} \label{lm:beta}
For all $\theta\in \R$, $d^\nabla \beta_\theta=0$, where $\nabla$ is the Levi-Civita
connection of $h$.
\end{lemma}

\begin{proof}
Let $u,v$ be two vector fields on $S$. 
Note first that $\beta_\theta$ satisfies Codazzi equation: 
\begin{eqnarray*}
(d^\nabla \beta_\theta)(u,v) & = & \nabla_u(\beta_\theta v) - \nabla_v(\beta_\theta u)-
\beta_\theta([u,v]) \\
& = & \cos(\theta/2)(\nabla_u v - \nabla_v u - [u,v]) + \sin(\theta/2)(\nabla_u (Jbv) - \nabla_v (Jbu) - 
Jb[u,v]) \\
& = & \sin(\theta/2) J(\nabla_u (bv) - \nabla_v (bu) - b[u,v]) \\
& = & 0~,
\end{eqnarray*}
where the last equality follows from the fact that $b$ satisfies Codazzi equation.
\end{proof}

\begin{lemma} \label{lm:nabla}
The Levi-Civita connection $\nabla^\theta$ of $h_\theta$ is given by
$$ \nabla^\theta_u v =\beta_{-\theta}\nabla_u(\beta_\theta v)~. $$
\end{lemma}

\begin{proof}
Consider the connection $\nabla^\theta$ defined in the statement
of the lemma. It is sufficient to prove that it is compatible
with $h_\theta$ and torsion-free. 

Let $u,v,w$ be three vector fields on $S$. Then
\begin{eqnarray*}
  u\cdot h_\theta(v,w) & = & u\cdot h(\beta_\theta v,\beta_\theta w) \\
& = & h(\nabla_u(\beta_\theta v),\beta_\theta w) + h(\beta_\theta v,\nabla_u(\beta_\theta w))\\
& = & h_\theta(\nabla^\theta_u v,w)+h_\theta(v,\nabla^\theta_u w)~,
\end{eqnarray*}
and therefore $\nabla^\theta$ is compatible with $h_\theta$.

We can now compute the torsion of $\nabla^\theta$ on $u,v$.
\begin{eqnarray*}
\nabla^\theta_u v - \nabla^\theta_v u -[u,v] & = & 
\beta_{-\theta}(\nabla_u(\beta_\theta v)-\nabla_v(\beta_\theta u) - \beta_\theta[u,v]) \\
& = & \beta_{-\theta}(d^\nabla \beta_\theta)(u,v) \\
& = & 0~. 
\end{eqnarray*}
Thus $\nabla^\theta$ is torsion-free and compatible with $h_\theta$, and so it is the
Levi-Civita connection of $h_\theta$.
\end{proof}

The computation of the connection of $h_\theta$ allows to prove easily that
$h_\theta$ is a hyperbolic metric.

\begin{proof}[Proof of Proposition \ref{pr:basics}]
Let $(e_1,e_2)$ be an orthonormal frame for $h$ on $S$, with connection 1-form $\omega$ 
with respect to the Levi-Civita connection $\nabla$ of $h$.
Let $(e'_1,e'_2)=(\beta_\theta^{-1}e_1, \beta_\theta^{-1}e_2)$, it is an orthonormal moving
frame for $h_\theta$. Since the Levi-Civita connection $\nabla^\theta$ of $h_\theta$ is
given by Lemma \ref{lm:nabla}, a direct computation shows that the connection 1-form of
$(e'_1,e'_2)$ for $\nabla^\theta$ is also equal to $\omega$. Since the identity map
between $(S,h)$ and $(S,h_\theta)$ is area-preserving, it follows that $h$ and $h_\theta$
have the same curvature, so that $h_\theta$ is also hyperbolic.
\end{proof}

\subsection{A cyclic property}

Our next goal is to prove Theorem \ref{tm:cyclic} and Theorem \ref{thm:c}. The proofs are in the next sections, after
some preliminary lemmas. 

\begin{lemma}
Let $\theta,\theta'\in \R$. Then $\beta_\theta\circ \beta_{\theta'}=\beta_{\theta+\theta'}$.
\end{lemma}

\begin{proof}
By definition we have:
\begin{eqnarray*}
\beta_\theta\circ\beta_{\theta'} & = & (\cos(\theta/2)E+\sin(\theta/2)Jb)\circ 
(\cos(\theta'/2)E+\sin(\theta'/2)Jb) \\
& = & \cos(\theta/2)\cos(\theta'/2) E + \sin(\theta/2)\sin(\theta'/2) JbJb +
(\cos(\theta/2)\sin(\theta'/2)+\sin(\theta/2)\cos(\theta'/2))Jb~.
\end{eqnarray*}
But $JbJb=-E$ because $b$ is self-adjoint and of determinant $1$. It follows
that 
$$ \beta_\theta\circ\beta_{\theta'} = \cos((\theta+\theta')/2) E +
\sin((\theta+\theta')/2) Jb = \beta_{\theta+\theta'}~. $$
\end{proof}

For instance, since $Jb=\beta_\pi$, it follows that $Jb\circ \beta_\theta = 
\beta_\theta \circ  Jb =\beta_{\theta+\pi}$. 

\begin{lemma} \label{lm:b}
Let $b_\theta:TS\rightarrow TS$ be the bundle morphism associated to 
$h_\theta$ and $h^\dual_\theta$ by Corollary \ref{cr:b}. Then 
$$ b_\theta = \beta_{-\theta}\circ b\circ \beta_\theta~. $$ 
\end{lemma}

\begin{proof}
Let $b_\theta = \beta_{-\theta}\circ b\circ \beta_\theta$. 
This lemma will follow if we prove the following facts:
\begin{enumerate}
\item $h_\theta(b_\theta\cbull, b_\theta\cbull)$ is isotopic to 
$h^\dual_\theta$ (we will actually prove that it is equal to 
$h^\dual_\theta$).
\item $b_\theta$ has determinant $1$.
\item $b_\theta$ is self-adjoint for $h_\theta$.
\item $d^{\nabla^\theta}\beta_\theta=0$. 
\end{enumerate}
For the first point, note that $Jb=\beta_\pi$, and therefore:
\begin{eqnarray*}
h_\theta(b_\theta\cbull, b_\theta\cbull) & = & 
h(b\beta_{\theta}\cbull, b\beta_\theta\cbull)
=  h(Jb\beta_{\theta}\cbull, Jb\beta_\theta\cbull)\\
& = & h(\beta_\pi\circ\beta_\theta\cbull,\beta_\pi\circ\beta_\theta\cbull)
= h(\beta_{\pi+\theta}\cbull,\beta_{\pi+\theta}\cbull) \\
& = & h^\dual_\theta~. 
\end{eqnarray*}

The second point is clear since $\det(b_\theta)=
\det(\beta_{-\theta})\det(b)\det(\beta_\theta)=1$.

For point (3) let $u,v$ be two vector fields on $S$, then
\begin{eqnarray*}
h_\theta(b_\theta u,v) & = & h(\beta_\theta b_\theta u, \beta_\theta v)
= h(b \beta_\theta u, \beta_\theta v) \\
& = & h(\beta_\theta u, b\beta_\theta v)
= h_\theta(u,b_\theta v)~. 
\end{eqnarray*}

For point (4), let again be $u,v$ be two vector fields on $S$. Then,
using the expression of $\nabla^\theta$ in Lemma \ref{lm:nabla}, we have:
\begin{eqnarray*}
(d^{\nabla^\theta}b_\theta)(u,v) & = & \nabla^\theta_u(b_\theta v) - 
\nabla^\theta_v(b_\theta u) - b_\theta[u,v] \\
& = & \beta_{-\theta}\nabla_u(\beta_\theta b_\theta v) - 
\beta_{-\theta}\nabla_v(\beta_\theta b_\theta u) - b_\theta[u,v] \\
& = & \beta_{-\theta} (d^\nabla \beta_{\theta+\pi})(u,v) 
\end{eqnarray*}
so that $(d^{\nabla^\theta}b_\theta)(u,v)=0$ by Lemma \ref{lm:beta}.
This completes the proof.
\end{proof}
\begin{remark}
Lemma \ref{lm:b} shows that $h_\theta$ and $h_{\theta+\pi}$ are normalized metrics.
\end{remark}
\subsection{Proof of Theorem \ref{tm:cyclic}}
\label{ssc:cyclic}

Let $\theta,\theta'\in \R$. By definition of $L_{\theta'}$, 
$$ L_{\theta'}(L_\theta(h,h^\dual)) = (h_\theta(\betab_{\theta'}\cbull,
\betab_{\theta'}\cbull),(h_{\theta+\pi}(\betab_{\theta'}\cbull,
\betab_{\theta'}\cbull))~, $$
where 
$$ \betab_{\theta'} = \cos(\theta'/2) E + \sin(\theta'/2) J_\theta b_\theta~, $$
where $J_\theta$ is the complex structure of $h_\theta$. Clearly, 
$J_\theta=\beta_{-\theta}\circ J\circ \beta_\theta$, so that 
\begin{eqnarray*}
\betab_{\theta'} & = & \beta_{-\theta}\circ(\cos(\theta'/2)E+\sin(\theta'/2)Jb)\circ\beta_\theta \\
& = &  \beta_{-\theta}\circ \circ\beta_{\theta'}\circ\beta_\theta \\
& = & \beta_{\theta'}~. 
\end{eqnarray*}
We now see that 
\begin{eqnarray*}
L_{\theta'}(L_\theta(h,h^\dual)) & = & (h_\theta(\beta_{\theta'}\cbull,
\beta_{\theta'}\cbull),(h_{\theta+\pi}(\beta_{\theta'}\cbull, \beta_{\theta'}\cbull)) \\
& = & (h(\beta_\theta\circ\beta_{\theta'}\cbull,\beta_\theta\circ\beta_{\theta'}\cbull),
h(\beta_{\theta+\pi}\circ\beta_{\theta'}\cbull,\beta_{\theta+\pi}\circ\beta_{\theta'}\cbull)) \\
& = & (h(\beta_{\theta'+\theta}\cbull, \beta_{\theta'+\theta}\cbull), 
h^\dual(\beta_{\theta'+\theta}\cbull, \beta_{\theta'+\theta}\cbull)) \\
& = & L_{\theta'+\theta}(h,h^\dual)~.
\end{eqnarray*}
This proves Theorem \ref{tm:cyclic}.

\subsection{Proof of Theorem \ref{thm:c}}
\label{ssc:c}

%

The fact that $h_\theta$ and $h_{\theta+\pi}$ are normalized metrics follows from 
Lemma \ref{lm:b}.

We 
compute the expression of $c_\theta$.
\begin{eqnarray*}
c_\theta & = & h_\theta + h_{\theta+\pi} \\
& = & h((\cos(\theta/2)E+\sin(\theta/2)b)\cbull,(\cos(\theta/2)E+\sin(\theta/2)b)\cbull) + \\
& & + h((-\sin(\theta/2)E+\cos(\theta/2)b)\cbull,(-\sin(\theta/2)E+\cos(\theta/2)b)\cbull) \\
& = & (\cos^2(\theta/2)+\sin^2(\theta/2))h + (\cos^2(\theta/2)+\sin^2(\theta/2))h(b\cbull,b\cbull) \\
& = & c~, 
\end{eqnarray*}
so $c_\theta$ is indeed independent of $\theta$.

The fact that the identity $(S,c)\rightarrow(S,h_\theta)$ is harmonic
follows from the last paragraph of Subsection \ref{ssc:24}.

For point (4) recall that $b_\theta$ has determinant $1$ and that 
$h_\theta^\dual=h_\theta(b_\theta\cbull,b_\theta\cbull)$. 
A simple computation then shows
that $h_\theta-h_\theta^\dual$ is traceless with respect to $h_\theta+h_\theta^\dual$. The definition
of Hopf differential (see Section \ref{ssc:harmonic}) therefore shows that, if
$\Phi_\theta$ is the Hopf differential of the identity from $(S,c)$ to $(S,h_\theta)$, 
then 
$$ \frac{h_\theta-h^\dual_\theta}{2} =  2 Re(\Phi_\theta)~. $$
Therefore
\begin{eqnarray*}
4 Re(\Phi_\theta) & = & h((\cos(\theta/2)E+\sin(\theta/2)b)\cbull,(\cos(\theta/2)E+\sin(\theta/2)b)\cbull) - \\
&  & - h((-\sin(\theta/2)E+\cos(\theta/2)b)\cbull,(-\sin(\theta/2)E+\cos(\theta/2)b)\cbull) \\
& = & (\cos^2(\theta/2)-\sin^2(\theta/2))(h-h^\dual) + 4\sin(\theta/2)\cos(\theta/2)h(b\cbull,\cbull) \\
& = & \cos(\theta)(h-h^\dual) +2\sin(\theta)h(b\cbull,\cbull)~.
\end{eqnarray*}
So $Re(\Phi_\theta)=Re(e^{i\theta}\Phi)$, where 
\begin{equation}\label{eq:hopf}
\Phi = \frac{h-h^\dual + 2ih(b\cbull,\cbull)}{4}=
\frac{h((E-b^2)\cbull,\cbull) + 2ih(b\cbull,\cbull) }{4}
\end{equation}
is the Hopf differential of the identity from $(S,c)$ to $(S,h)$.

\subsection{Centers}

We conclude this section by some remarks on the respective behavior of $h,h^\dual$
and $c$. 

\begin{remark}
$c$ is uniquely determined by $h,h^\dual\in \cT$. Conversely, $h^\dual$ is uniquely
determined by $h$ and $c$.
\end{remark}

\begin{proof}
Given $h$ and $h^\dual$, we have seen that there is a unique minimal Lagrangian diffeomorphism
$m:(S,h)\rightarrow (S,h^\dual)$ isotopic to the identity. By definition, $c$ is the conformal
class of $h+m^*h^\dual$. 

For the second point, let $f:(S,c)\rightarrow (S,h)$ be the unique harmonic map, given by
Theorem \ref{tm:harmonic}, and let $\Phi$ be its Hopf differential. There is then a unique
harmonic map $f^\dual$ from $(S,c)$ to a hyperbolic surface $(S,h^\dual)$ with Hopf differential equal
to $-\Phi$. The content of Section \ref{ssc:24} then indicates that $f^\dual\circ f^{-1}$ is
minimal Lagrangian, so that $c$ is the center of $(h,h^\dual)$.
\end{proof}

\section{An extension of Thurston's Earthquake Theorem}

In this section we recall a recent result of \cite{BBZ2} on constant
curvature folations of GH AdS manifolds, and use it to prove Theorem
\ref{tm:earthquake}.

\subsection{Constant curvature foliations in AdS geometry}

We recall here one of the main result of \cite{BBZ2}. Let $N$ be a
MGH AdS 3-dimensional manifold, let $C(N)$ be its convex core.

\begin{theorem}[Barbot, B\'eguin, Zeghib \cite{BBZ2}] \label{tm:bbz}
The complement of $C(N)$ in $N$ is foliated by surfaces of constant 
(Gauss) curvature. Moreover, for any $k\in (-\infty, -1)$, there exists
a unique future-convex (resp. past-convex) surface of constant 
curvature $k$ in $N$, and it is a leaf of the foliation.
\end{theorem}

\subsection{Proof of Theorem \ref{tm:earthquake}}
\label{ssc:bbz}

We first translate Theorem \ref{tm:bbz} in terms of the lanslide flow,
using Lemma \ref{lm:ads}.

\begin{cor} \label{cr:exp}
Choose $(\rho_l,\rho_r)\in \cT$ and $\alpha\in (0,\pi)$. There exists a 
unique $(h,h^\dual)\in \cT\times \cT$ such that
$$ L^1_{e^{i\alpha}}(h,h^\dual)=\rho_l ~ \mbox{and}~  L^1_{e^{-i\alpha}}(h,h^\dual)=\rho_r~. $$ 
\end{cor}

\begin{proof}
It is a direct consequence of Theorem \ref{tm:bbz}. Given 
$\rho_l$ and $\rho_r$ there is a unique GHMC AdS manifold $N\cong S\times\R$ of which
they are the left and right representations, respectively. $N$ contains
a unique past-convex surface $F$ with constant curvature $-1/\cos^2(\alpha/2)$,
which comes with a diffeomorphism $\bar{\phi}:S\rightarrow F$ (well-defined up to
isotopy).
We call $I$ and $\III$ the induced metric and third fundamental form on $S$,
respectively. Then $\III$ has constant curvature $-1/\sin^2(\alpha/2)$.
We then set $h=(1/\cos^2(\alpha/2))I, h^\dual=(1/\sin^2(\alpha/2)) \III$, 
so that $h$ and $h^\dual$ are normalized hyperbolic metrics on $S$
(see \cite{minsurf}).
Lemma \ref{lm:ads} then shows that $\rho_l=L^1_{e^{i\alpha}}(h,h^\dual)$, while
$\rho_r=L^1_{e^{-i\alpha}}(h,h^\dual)$.

Conversely, given $h,h^\dual\in \cT$ such that $\rho_l=L^1_{e^{i\alpha}}(h,h^\dual)$, while
$\rho_r=L^1_{e^{-i\alpha}}(h,h^\dual)$, we can consider the unique equivariant
embedding $\phi:\tilde{S}\rightarrow\AdS^3$
onto a past-convex surface $\tilde{F}$,
with induced metric $\cos^2(\alpha/2)h$ and
third fundamental form $\sin^2(\alpha/2)h^\dual$. Then $\tilde{F}$ is the lift to 
$\AdS^3$ of a past-convex surface $F$ in a GHMC AdS manifold $N$, and the
left and right representations of $N$ are $\rho_l$ and $\rho_r$ by 
Lemma \ref{lm:ads}. This shows the uniqueness.
\end{proof}

\begin{proof}[Proof of Theorem \ref{tm:earthquake}]
Apply Corollary \ref{cr:exp} with $\rho_r=h$, $\rho_l=h'$, and with
$\alpha=\theta/2$. It shows there exists a unique $h_0\in\cT$ and a unique
$h_0^\dual\in\cT$ such that
\begin{equation}\label{dd:eq}
    h=L^1_{e^{-i\alpha}}(h_0, h_0^\dual)\,,\quad h'=L^1_{e^{i\alpha}}(h_0, h_0^\dual)~.
\end{equation}
As a consequence, putting $h^\dual=L^2_{e^{-i\alpha}}(h_0, h_0^\dual)$ we have
$L^1_{e^{i\theta}}(h, h^\dual)=h'$.

Conversely given $h^\dual\in\cT$ such that $L^1_{e^{i\theta}}(h,h^\dual)=h'$, letting
$(h_0,h_0^\dual)=L_{e^{i\theta/2}}(h,h^\dual)$ we easily see that (\ref{dd:eq})
is verified. The uniqueness in Corollary \ref{cr:exp} therefore implies the uniqueness
here.
%
\end{proof}

\section{The complex cyclic map} \label{sc:complex}

This section describes a natural extension of the cyclic flow from a real
to a complex parameter. This is analogous to the complex earthquake introduced
by McMullen \cite{mcmullen:complex}. We will actually show in the next section 
that the ``complex cyclic flow''
introduced here limits, in a suitable sense, to the complex earthquake. 

\subsection{Main statements}
\label{ssc:main}

Let $\mathcal P$ be the space of complex projective surface on $S$. 
The space $\mathcal P$ is naturally a complex manifold of real dimension $12g-12$. 
Moreover the natural map $\mathcal P\rightarrow\mathcal T$ that associates 
to a complex projective surface the underlying complex structure is holomorphic.
A projective structure is {\it Fuchsian} 
if its universal covering is projectively equivalent to a round disk
in $\mathbb{CP}^1$.

Let $\mathbb{H}$ be the upper half-plane in $\mathbb C$.
We define a map 
\[
    P':\overline{\mathbb H}\times\cT\times\mathcal T\rightarrow\mathcal P
\]
where
\begin{itemize}
\item For every fixed $h,h^\dual$ the map $z\rightarrow P_{z}(h,h^\dual)$ is holomorphic.
\item For $t$ real, $P'_t(h,h^\dual)$ is the Fuchsian projective surface corresponding to
$L^1_{e^{-it}}(h,h^\dual)$.
\end{itemize}

Notice that here $z$ lives in the upper half-plane, while, in the introduction, the
flow usually depended on a complex parameter $\zeta$ in the unit disk.
Both parameterizations are quite
useful here. Taking $\zeta$ in the unit disk is natural when dealing with the landslide
flow, while taking $z$ in the upper half-space is natural when thinking of complex
earthquakes as a limit (since complex earthquakes are usually parameterized by the
upper half-space). Until Section \ref{ssc:cocycles} we consider the parameterization
by $z$ in the upper half-space, while in Section \ref{ssc:another} we will make
the connection to the parameterization by $\zeta$ in the unit disk. We use a prime to
denote the various maps when $z$ is in the upper half-plane, this explains the
notation $P'$ above.

The construction of the map $P'$ is the analog of 
the construction of the complex earthquake due to McMullen \cite{mcmullen:complex}. 
The first point is to define the analog of the grafting.

Given two normalized hyperbolic metrics $h$ and $h^\dual$ 
on $\mathcal T$, let $b$ be the operator defined in Subsection \ref{ssc:minilag}.
Given a positive number $s>0$, we consider the metric $I_s=\cosh^2(s/2)h$ and
the operator
$B_s=-\tanh(s/2)b$.
It is easy to see that $(I_s, B_s)$
satisfies the Gauss-Codazzi equation for immersed surfaces
in $\Hyp^3$, that is
\[
   d^{\nabla}B_s=0~,\qquad\qquad K_s=-1+\det B_s~,
\]
where $\nabla$ is the Levi-Civita connection for $I_s$ 
(which is equal to the Levi-Civita connection for $h$) 
and $K_s$ is the curvature of $I_s$ (which is constant and equal to $-1/\cosh^2(s/2)$).

As a consequence there exists a convex equivariant immersion
\begin{equation}\label{graft:eq}
  \sigma_s:\tilde S\rightarrow \Hyp^3 
\end{equation}
whose first fundamental form is the pull-back $\tilde{I}_s$ of $I_s$ 
and whose shape operator is the pull-back $\tilde B_s$ of $B_s$. 
This map $\sigma_s$ is uniquely determined up to elements of $\PSL_2(\mathbb C)$,
once we state that the orientation on $\tilde S$ at $\sigma_s(\tilde{p})$ coincides with 
the orientation induced by the normal vector $\tilde\nu_s(\tilde{p})$ pointing towards the concave part 
(this is the reason why the sign of $B_s$ is negative).

Given $\tilde{p}\in\tilde S$, let $dev_s(\tilde{p})\in S^2_\infty=\mathbb{CP}^1$
be the endpoint of the geodesic ray 
starting from $\sigma_s(\tilde{p})$ with velocity $\tilde\nu_s(\tilde{p})$. The map
\[
   dev_s:\tilde S\rightarrow \mathbb{CP}^1 
\]
is a developing map for a complex projective structure $SGr'_{s}(h,h^\dual)$ on $S$. (The notation
$SGr$ is used to keep in mind the analogy to the grafting map $Gr$).

Notice that, if $h=h^\dual$, then $b=E$
and 
$\sigma_s=d_s\circ\sigma_0$, where $d_s:\sigma_0(\tilde S)=\mathbb \Hyp^2\rightarrow \Hyp^3$
is the map associating to $\sigma_{s_0}(\tilde{p})$ the end-point of the geodesic segment of length $s$ starting from
$\sigma_{s_0}(\tilde{p})$ with velocity $\tilde\nu_{s_0}(\tilde{p})$.
So, in this case,  $SGr'_{s}(h,h)$ is the
Fuchsian projective structure associated to $h$.

Finally, given a complex number $z=t+is$ with $s\geq 0$, we define
\[
     P'_z(h,h^\dual)=SGr'_{s}(L'_{-t})(h,h^\dual))
\]
where $L'_{-t}(h,h^\dual):=L_{e^{-it}}(h,h^{\dual})$.

Most of the remaining part of this section is devoted to proving the following theorem.

\begin{theorem} \label{tm:complex}
The map
\[
   z\mapsto P'_{z}(h,h^\dual)\in\mathcal P
\]
is holomorphic.
\end{theorem}

Composing $P'$ with the forgetful map from $\cP$ to $\cT$, we obtain for each 
 $z$ in the upper half-plane a map 
$$ C'_z:\cT\times\cT\rightarrow \cT $$
sending $(h,h^\dual)\in \cT\times\cT$ to the complex structure underlying the complex projective structure 
$P'_{z}(h,h^\dual)$. 

\begin{cor}
The map $z\rightarrow C'_{z}(h,h^\dual)$ is holomorphic in the upper half-plane.
\end{cor}

This clearly follows from Theorem \ref{tm:complex} since the forgetful map from 
$\cP$ to $\cT$ is holomorphic.

To prove Theorem \ref{tm:complex},
we  will show that the holonomy $\rho_z$ of $P'_{z}$, holomorphically
depends on $z$.
In fact the derivatives 
\[
    \xi_z=\frac{\partial \rho_z}{\partial t}~,~ \qquad
\eta_z=\frac{\partial \rho_z}{\partial s}
\]
are $\psl_2(\mathbb{C})$-valued cocycles in $H^1(\pi_1(S), Ad\circ \rho_z)$ and we will show that
\begin{equation}\label{hol:eq}
     \eta_z=i\xi_z\,.
\end{equation}

Let us remark that, since $P'_{z+t}(h,h^\dual)=P'_z(L_{e^{-it}}(h,h^\dual))$ for any 
$z$ in the upper plane and $t$ real, it is sufficient to check
(\ref{hol:eq}) at imaginary points $z_0=is_0$.

To compute the cocycles we consider the family of convex immersions
\[
    \sigma_s,\tau_t:\tilde S\rightarrow \Hyp^3
\]
such that $\sigma_s$ corresponds to $SGr'_{s}(h,h^\dual)$ and $\tau_t$ corresponds to
$SGr'_{s_0}(L_{e^{-it}}(h,h^\dual))$.

The first-order variations of $\sigma_s$ and  $\tau_t$ are the fields
\begin{equation} \label{eq:XY}
      X=\frac{\partial\sigma_s}{\partial s}|_{s=s_0}~,  \qquad Y=\frac{\partial\tau_t}{\partial t}|_{t=0}
\end{equation}
regarded as sections of the fiber-bundle $\Theta=\sigma_{s_0}^*(T\Hyp^3)$ on $\tilde S$.
Imposing the equivariance of $\sigma_s$ under $\rho_{is}$, and of $\tau_t$ under $\rho_{t+is_0}$,
we deduce that
\begin{eqnarray}\label{field:eq}
   X(\tilde{p})= (\rho(\gamma))_* X(\gamma^{-1}\tilde{p})+\eta(\gamma)(\tilde{p})\\
   Y(\tilde{p})=(\rho(\gamma))_*Y(\gamma^{-1}\tilde{p})+\xi(\gamma)(\tilde{p})\label{field2}
\end{eqnarray}
where we have put $\rho=\rho_{z_0}$, $\eta=\eta_{z_0}$, $\xi=\xi_{z_0}$ and we are identifying
the elements of $\psl_2(\mathbb C)$ with the Killing vector fields on $\Hyp^3$.

We will find some explicit relation between $X$ and $Y$ 
that, used in (\ref{field:eq}) and (\ref{field2}), shall show equation (\ref{hol:eq}).

\subsection{General formulas}

We consider any smooth family of immersions
\[
    \sigma_s:\tilde S\rightarrow \Hyp^3
\]
and let $\tilde{I}_s$ be the first fundamental form on $\tilde S$ and $\tilde{B}_s$ the
shape operator associated with $\sigma_s$.

Fix $s_0>0$: we study $\sigma_{s}$ around $s=s_0$.
Let us denote by $\Theta$ the vector bundle $\sigma_{s_0}^*(T\Hyp^3)$
and let $X=\frac{\partial\sigma}{\partial s}|_{s=s_0}$, seen as a section of $\Theta$.
In this section, we will express the derivatives of $\tilde{I}_s$ and $\tilde{B}_s$ at $s=s_0$ in terms of the field $X$,
and we will show that these quantities determine $X$ up to global Killing vector fields.

We notice that there is a natural inclusion of $T\tilde S$ into $\Theta$ given by $d\sigma_{s_0}$.
For the sake of simplicity, we will identify $T\tilde S$ with its image in $\Theta$.
Given a point $\tilde{p}\in\tilde S$ and $s>0$ we denote by $\tilde{\nu}_s(\tilde{p})$ the unit vector at
$\sigma_s(\tilde{p})$ orthogonal to $d\sigma_s(\tilde{p})(T_{\tilde{p}}\tilde S)$
such that, for every positive basis $\{e_1, e_2\}$ of $T_{\tilde{p}}\tilde S$, the vectors $\tilde{\nu}_s(\tilde{p}), d\sigma_s(\tilde{p})(e_1),
d\sigma_s(\tilde{p})(e_2)$ form a positively oriented basis of $T_{\sigma_s(\tilde{p})}\Hyp^3$.
In this way, if $\tilde{J}$ denotes the complex structure on $\tilde S$ (and by abuse of
notation on $\sigma_s(\tilde S)$), we have
\begin{equation}\label{complex:eq}
  \tilde{J}_s\tilde{v}=\tilde{\nu}_s\times \tilde{v}
\end{equation}
for all $\tilde{v}\in d\sigma_{s}(T\tilde S)$.

A linear connection $D$ is defined on $\Theta$ 
by pulling back the Levi-Civita connection on $T\Hyp^3$.

The covariant derivative of a section $Y$ of $\Theta$ is a linear operator
\begin{equation}\label{decomp:eq}
D Y: T \tilde S\rightarrow T\Hyp^3\,.
\end{equation}
Such an operator can be decomposed as the sum of a self-adjoint operator
of $T\tilde S$ (identified to a subspace of $T\Hyp^3$) and the restriction of a 
skew-symmetric operator  of $T\Hyp^3$.


\begin{lemma}
Given a section $V$ of $\Theta$, there exist
\begin{itemize}
\item a self-adjoint operator $A^V$ of $T\tilde S$;
\item a section $S^V$ of $\Theta$
\end{itemize}
such that
\begin{equation}\label{dec:eq}
   D_{\tilde{v}}V=A^V(\tilde{v})+S^V\times \tilde{v}
\end{equation}
for every $\tilde{v}\in T\tilde S$.
Moreover both $A^V$ and $S^V$ are uniquely determined.
\end{lemma}

\begin{proof}
At every point, $D V$ can be decomposed in a tangential part and normal part:
\[
   D_{\tilde{v}} V=\alpha(\tilde{v})+\langle D_{\tilde{v}} V,\tilde{\nu}\rangle\tilde{\nu}\,.
\]
Clearly $\alpha$ is an operator of $T\tilde S$, so it can be decomposed into a self-adjoint part
$A(\alpha)$ and a skew-symmetric part $S(\alpha)$.
Notice that the skew-symmetric part is a multiple of $\tilde{J}$, in particular there is $a\in\mathbb{R}$ such that
$S(\alpha)(\tilde{v})=a \tilde{J}(\tilde{v})=a\tilde{\nu}\times \tilde{v}$.
On the other hand, there exists a tangent vector $\tilde{w}$ such that
\[
    \langle D_{\tilde{v}}V,\tilde{\nu}\rangle=\langle \tilde{v},\tilde{w}\rangle
\]
for every $\tilde{v}\in T\tilde S$.
In particular the normal part of $D_{\tilde{v}} V$ can be regarded as the restriction
on $T\tilde S$ of the skew-symmetric operator
\[
    \tilde{v} \mapsto \langle \tilde{v},\tilde{w}\rangle\tilde{\nu}-\langle \tilde{v},\tilde{\nu}\rangle \tilde{w}=
    (\tilde{w}\times\tilde{\nu} )\times \tilde{v}\,.
\]
So if we put $A^V=A(\alpha)$ and $S^V=a\tilde{\nu}+(\tilde{w}\times\tilde{\nu})$,
Equation (\ref{dec:eq}) is verified.

We show now that this decomposition is unique.
Suppose $A$ is a self-adjoint operator of $T\tilde S$ and let $S$ be a
vector tangent to $\Hyp^3$ such that $ A^V(\tilde{v})+ S^V\times \tilde{v}=A(\tilde{v})+S\times \tilde{v}$
for all $\tilde{v}\in T\St$.
Clearly $(A^V-A)(\tilde{v})=(S-S^V)\times \tilde{v}$.
This implies that $S-S^V$ is a normal vector. Since $A-A^V$
is self-adjoint, it follows that $A^V-A=0$ and $S^V-S=0$.
\end{proof}

An important property of the covariant derivative of $X$ is that $D_{\tilde{v}} X$ is the variation 
of the image of $\tilde{v}$ in $\Hyp^3$ along the family $(\sigma_s)$. More precisely the following statement holds. 
Here we denote by $D/ds$ the covariant derivative along $(\sigma_s)$ associated to $D$.

\begin{lemma}\label{varvec:lem}
Given a tangent vector $\tilde{v}\in T_{\tilde{p}}\tilde S$, we consider the field $\tilde{v}_s= d\sigma_s(\tilde{v})$ 
along the curve $s\mapsto \sigma_s(\tilde{p})$. We have
\begin{equation}\label{var-vec:eq}
    \frac{D \tilde{v}_s}{ds}|_{s=s_0}=D_{\tilde{v}} X\,.
\end{equation}
\end{lemma}

\begin{proof}
Take a path $\upsilon:(-\delta,\delta)\rightarrow\tilde S$ such that
$\upsilon'(0)=\tilde{v}$ and consider the map
$\tilde{x}(\varepsilon,s)=\sigma_s(\upsilon(\varepsilon))$.
We have 
$\tilde{v}_s=\frac{\partial\tilde{x}}{\partial \varepsilon}(0,s)$
whereas $\frac{\partial\tilde{x}}{\partial s}(\varepsilon,s_0)=
X(\upsilon(\varepsilon))$.
A direct computations shows that
\[
\frac{D}{ds}\tilde{v}_s|_{s=s_0}=
\frac{D}{ds}\frac{\partial\tilde{x}}{\partial \varepsilon}(0,s_0)=
\frac{D}{d\varepsilon}\frac{\partial\tilde{x}}{\partial s}(0,s_0)=
\frac{D}{d\varepsilon}X(0,s_0)=
D_{\tilde{v}} X\,.
\]
\end{proof}

Now we apply the decomposition (\ref{dec:eq}) to the field $X$, so we call
$A^X$ the self-adjoint part of $D X$ and $X'$ the field $S^X$.
It turns out that the first order variation of $\tilde{I}_s$ is determined by $A^X$.
On the other hand, the field $X'$ determines the variation of the normal field along
the family of $\sigma_s$.

\begin{lemma} \label{lm:74}
Given $\tilde{u},\tilde{v}\in T_{\tilde{p}}\tilde S$, we have
\begin{equation}\label{var1:eq}
  \frac{d\,}{ds}\tilde{I}_s(\tilde{u},\tilde{v})|_{s=s_0}=2\tilde{I}_s(A^X(\tilde{u}),\tilde{v})\,,\qquad\qquad
  \frac{D\tilde{\nu}}{ds}|_{s=s_0}=X'\times\tilde{\nu}\,.
\end{equation}
\end{lemma}
\begin{proof}
We have that
\[
  \tilde{I}_s(\tilde{u},\tilde{v})=\langle d\sigma_s(\tilde{u}), d\sigma_s(\tilde{v})\rangle\,.
\]
Applying Lemma \ref{varvec:lem} we get
\[
 \frac{d}{ds}\tilde{I}_s(\tilde{u},\tilde{v})=\langle D_{\tilde{u}}X,\tilde{v}\rangle+\langle \tilde{u},D_{\tilde{v}}X\rangle\,.
\]
Since $D_{\tilde{u}}X=A^X(\tilde{u})+X'\times \tilde{u}$, where $A^X$ is a self-adjoint operator of $T\tilde S$,
we get the first formula.

To prove the second formula, first we notice that since
$\langle\tilde{\nu},\tilde{\nu}\rangle=1$, then $\frac{D\tilde{\nu}}{ds}|_{s=s_0}$ is a tangent field.
On the other hand, given a tangent vector $\tilde{v}$ we have
$\langle\tilde{\nu}, d\sigma_s(\tilde{v})\rangle=0$. Differentiating this identity we get
 \[
 \left\langle \frac{D\tilde{\nu}}{ds},\tilde{v}\right\rangle=-\langle\tilde{\nu},D_{\tilde{v}}X\rangle=
 \langle \tilde{v}, X'\times\tilde{\nu}\rangle~.
 \]
 Since $X'\times\tilde{\nu}$ is tangent, this proves that $\frac{D\tilde{\nu}}{ds}=X'\times\tilde{\nu}$.
\end{proof}

Using (\ref{var1:eq}) we get the first order variation of $\tilde{B}_s$ at $s=s_0$.
In the computation of such a variation, the covariant derivative of $X'$ appears.
It is useful to apply the decomposition (\ref{dec:eq}) to 
$D X'$. In particular, we put $X''=S^{X'}$, so that
\[
  D_{\tilde{v}} X'=A^{X'}(\tilde{v})+X''\times \tilde{v}\,.
\]

\begin{lemma} \label{lm:75}
Given $\tilde{v}\in T_{\tilde{p}}\tilde S$, we have
\begin{equation}\label{var2:eq}
  \frac{d\,}{ds}(\tilde{B}_s(\tilde{v}))|_{s=s_0}=
\tilde{J}A^{X'}(\tilde{v})-\langle X+X'',\tilde{\nu}\rangle \tilde{v}-A^X\circ \tilde{B}_{s_0}(\tilde{v})\,.
\end{equation}
\end{lemma}

\begin{proof}
Differentiating with respect to $s$ the identity
\[
  d\sigma_s(\tilde{B}_s(\tilde{v}))=-D_{\tilde{v}}\tilde{\nu}
\]
and evaluating at $s=s_0$, we obtain
\begin{equation}\label{p:eq}
D_{\tilde{B}_{s_0}(\tilde{v})}X+\dot{\tilde{B}} (\tilde{v})=-\frac{D}{ds}(D_{\tilde{v}}\tilde{\nu})\,.
\end{equation}
On the other hand, we have that
\begin{equation}\label{pp:eq}
 \frac{D}{ds}(D_{\tilde{v}}\tilde{\nu})=D_{\tilde{v}}\left(\frac{D\tilde{\nu}}{ds}\right)+\bar R(X, \tilde{v})\tilde{\nu}
\end{equation}
where $\bar R$ is the Riemann tensor of $\Hyp^3$.

By (\ref{var1:eq}) we have that
\begin{equation}\label{ppp:eq}
D_{\tilde{v}}\left(\frac{D\tilde{\nu}}{ds}\right)=D_{\tilde{v}}(X'\times\tilde{\nu})=
A^{X'}(\tilde{v})\times\tilde{\nu}+(X''\times \tilde{v})\times\tilde{\nu}-X'\times \tilde{B}_{s_0}(\tilde{v})\,.
\end{equation}
On the other hand, since $\Hyp^3$ has constant curvature $-1$, its Riemann tensor is
simply given by 
\begin{equation}\label{pppp:eq}
\bar R(X,\tilde{v})\tilde{\nu}=\langle X,\tilde{\nu}\rangle \tilde{v}-\langle \tilde{v},\tilde{\nu}\rangle X=\langle X,\tilde{\nu}\rangle \tilde{v}\,.
\end{equation}

Using (\ref{pp:eq}), (\ref{ppp:eq}), and (\ref{pppp:eq}) in (\ref{p:eq}) we get
\[
A^X(\tilde{B}_{s_0}(\tilde{v}))+X'\times \tilde{B}_{s_0}(\tilde{v})
+\dot{\tilde{B}} (\tilde{v})=-A^{X'}(\tilde{v})\times\tilde{\nu}-(X''\times \tilde{v})\times\tilde{\nu}+X'\times \tilde{B}_{s_0}(\tilde{v})-
\langle X,\tilde{\nu}\rangle \tilde{v}
\]
Since $\tilde{\nu}\times A^{X'}(\tilde{v})=\tilde{J}A^{X'}(\tilde{v})$ and
$(X''\times \tilde{v})\times\tilde{\nu}=\langle X'',\tilde{\nu}\rangle \tilde{v}$, Equation (\ref{var2:eq}) follows.
\end{proof}

Finally we show that Equations (\ref{var1:eq}) and (\ref{var2:eq}) determine
$X$ up to some global vector field. This is an easy consequence of the following lemma.

\begin{lemma}\label{uniq:lem}
Let $V$ be a section of $\Theta$ and let us put $V'=S^V$ and 
$V''=S^{V'}$.
Suppose that
\[
A^V=0\qquad\qquad \tilde{J}A^{V'}-\langle V''+V,\tilde{\nu}\rangle E-A^V\circ \tilde{B}_{s_0}=0\,.
\]
Then $V$ is the restriction of a global Killing field of $\Hyp^3$ on $\tilde S$.
\end{lemma}

\begin{proof} 
Under the hypothesis of the lemma, neither the induced metric nor the shape operator
of the surface vary under the first-order deformation defined by $V$. The conclusion
therefore follows from the Fundamental Theorem of surface theory, see e.g. \cite{spivak}.
\end{proof}

\subsection{The variation field of $SGr'$}

In this section we apply formulas obtained in the previous subsection
to the family of convex immersions $\sigma_s:\tilde S\rightarrow \Hyp^3$
defined in (\ref{graft:eq}).

\begin{lemma}\label{graftvar:lem}
For $X=\frac{\partial\sigma}{\partial s}|_{s=s_0}$, denote by $X'$ the section $S^X$ and by $X''$
the section $S^{X'}$. The following formulas hold:
\begin{eqnarray}
2A^X=\tanh(s_0/2)E \label{gr1:eq}\\
A^{X'}=[\tilde{J},\tilde{b}]/4  \label{gr2:eq}\\
\langle X+X'', \tilde{\nu}\rangle=\frac{\tr (\tilde{b})}{4}\label{gr3:eq}.
\end{eqnarray}
\end{lemma}

\begin{proof}
The embedding data corresponding to $\sigma_s$ are
\[
\tilde{I}_s=\cosh^2(s/2) \tilde{h} \qquad\qquad \tilde{B}_s=-\tanh(s/2) \tilde{b}
\]
so we easily get that
\[
  \frac{d}{ds}\tilde{I}_s(\tilde{u},\tilde{v})|_{s=s_0}=\tanh(s_0/2)\tilde{I}_{s_0}(\tilde{u},\tilde{v}).
\]
Comparing this formula with (\ref{var1:eq}), we get that $2A^X=\tanh(s_0/2)E$.

On the other hand, applying (\ref{var2:eq}) we get
\[
-\frac{1}{2\cosh^2(s_0/2)}\tilde{b}=\tilde{J}A^{X'}- \langle X''+X,\tilde{\nu}\rangle E+\frac{\tanh^2(s_0/2)}{2} \tilde{b}~,
\]
which can be also written 
\[
   \tilde{b}=-2\tilde{J}A^{X'} + 2\langle X''+X,\tilde{\nu}\rangle E~.
\]
Multiplying by $\tilde{J}$ we deduce that
\[
  \tilde{J}\tilde{b}=2A^{X'}+2\langle X''+X,\tilde{\nu}\rangle \tilde{J}\,.
\]
Notice that this must coincide with the  decomposition of $\tilde{J}\tilde{b}$ in symmetric and skew-symmetric part.
Since the adjoint of $\tilde{J}\tilde{b}$ is $-\tilde{b}\tilde{J}$ it follows that
$$ 2A^{X'}=\frac{\tilde{J}\tilde{b}-\tilde{b}\tilde{J}}{2}=\frac{[\tilde{J},\tilde{b}]}{2} $$
$$ 2\langle X''+X,\tilde{\nu}\rangle \tilde{J}=\frac{\tilde{J}\tilde{b}+\tilde{b}\tilde{J}}{2}=\frac{-\tilde{J}\tilde{b}\tilde{J}+\tilde{b}}{2}\tilde{J}=\frac{\tilde{b}^{-1}+\tilde{b}}{2}\tilde{J}=\frac{\tr(\tilde{b})}{2}\tilde{J} $$
where we have used that $\tilde{b}+\tilde{b}^{-1}=\tr(\tilde{b})E$.
\end{proof}

A consequence of Lemma \ref{graftvar:lem} is that $X''+X$ can be explicitly computed.

\begin{prop} \label{Zfield:prop}
With the notation of Lemma \ref{graftvar:lem} the following identity holds:
\[
   X''=-X+\frac{\tr(\tilde{b})}{4}\tilde{\nu}\,.
\]
\end{prop}

\begin{proof}
By (\ref{gr3:eq}), it is sufficient to prove that 
\[
   \langle X''+X, \tilde{v}\rangle=0
\]
for every tangent vector $\tilde{v}$.

Let $\tilde{u},\tilde{v}$ be two tangent fields on $\tilde S$. By using the identity
\[
   D_{\tilde{u}} X=\frac{1}{2}\tanh(s_0/2) \tilde{u}+X'\times \tilde{u} 
\]
an explicit computation shows that
\[
\bar R(\tilde{u},\tilde{v})X=
 D_{\tilde{u}}D_{\tilde{v}}(X)-D_{\tilde{v}}D_{\tilde{u}}(X)-D_{[\tilde{u},\tilde{v}]}X=
 D_{\tilde{u}}X'\times \tilde{v}-D_{\tilde{v}} X'\times \tilde{u}.
\]
Moreover,
\begin{equation}\label{q:eq}
D_{\tilde{u}}X'\times \tilde{v}-D_{\tilde{v}}X'\times \tilde{u}=
A^{X'}(\tilde{u})\times \tilde{v}-A^{X'}(\tilde{v})\times \tilde{u}\ 
+\ (X''\times \tilde{u})\times \tilde{v}- (X''\times \tilde{v})\times \tilde{u}\,. 
\end{equation}
By (\ref{gr2:eq}), $A^{X'}$ is a self-adjoint traceless operator, and it follows that 
the sum of the first two terms of (\ref{q:eq}) vanishes.
Eventually, we get
\[
\bar R(\tilde{u},\tilde{v})X=
(X''\times \tilde{u})\times \tilde{v}-(X''\times \tilde{v})\times \tilde{u}=
\langle X'',\tilde{v}\rangle \tilde{u}-\langle X'',\tilde{u}\rangle \tilde{v}\,.
\]
Since $\bar{R}(\tilde{u},\tilde{v})X=\langle X,\tilde{u}\rangle \tilde{v}-
\langle X,\tilde{v}\rangle \tilde{u}$, we easily deduce that
$\langle X''+X,\tilde{v}\rangle=0$ for all tangent vectors $\tilde{v}$.
\end{proof}

\begin{prop}\label{fieldsmain:prop}
Let $\tau_t:\tilde S\rightarrow \Hyp^3$ be the family of convex immersions corresponding to
the projective surface $SGr'_{s_0}\circ L_{e^{-it}}(h,h^\dual)$ 
and denote by $Y$ its first order variation at $t=0$.
Then $Y=-X'$ up to adding a global Killing vector field, where $X'$ is the vector field
defined in Lemma \ref{graftvar:lem}. 
\end{prop}

\begin{proof}
Let $\tilde{I}_t$ be the first fundamental form corresponding to $\tau_t$ and let 
$\tilde{B}_t$ be the corresponding shape operator.
According to Lemmas \ref{lm:74}, \ref{lm:75} and \ref{uniq:lem}, it is sufficient to show that
\begin{align}
\frac{d\,}{dt}\tilde{I}_t(\tilde{u},\tilde{v})|_{t=0}& =-2\tilde{I}_0( A^{X'}(\tilde{u}), \tilde{v}) \label{r:eq}\\
\frac{d\,}{dt}\tilde{B}_t(\tilde{u})|_{t=0} & =  -\left(\tilde{J}A^{X''}(\tilde{u})-\langle X'+X''',\tilde{\nu}\rangle \tilde{u}-A^{X'}\circ \tilde{B}_0(\tilde{u})\right)\,.\label{rr:eq}
\end{align}
where $X'''=S^{X''}$ is the vector field corresponding to the skew-symmetric part of  $D X''$.
Call $\beta_t=\cos (t/2)E-\sin (t/2)Jb$, so that we have
\begin{align}
I_t & =\cosh^2(s_0/2)h(\beta_t\cbull,\beta_t\cbull) \\
B_t & =-\tanh(s_0/2)\beta_{-t}b\beta_t\,.\label{eq:Bt}
\end{align}
It follows that
\begin{eqnarray*}
\frac{d}{dt} \tilde{I}_t(\tilde{u},\tilde{v})|_{t=0} & = & \frac{1}{2}
\left(-\tilde{I}_0(\tilde{J}\tilde{b}(\tilde{u}),\tilde{v})-\tilde{I}_0(\tilde{u},\tilde{J}\tilde{b}(\tilde{v}))\right)\\
& = & -\frac{1}{2}\tilde{I}_0([\tilde{J},\tilde{b}]\tilde{u},\tilde{v})=-2\tilde{I}_0(A^{X'}(\tilde{u}),(\tilde{v}))\,.
\end{eqnarray*}
To show equation (\ref{rr:eq}), we first compute $D X''$.
By Proposition \ref{Zfield:prop} we have
\begin{eqnarray*}
D_u X'' & = & -A^Xu-X'\times u-\frac{\tr(\tilde{b})}{4}\tilde{B}_0u+\left\langle\grad\left(\frac{tr(\tilde{b})}{4}\right),u\right\rangle\tilde{\nu} \\
& = & -\left(A^X+\frac{\tr(\tilde{b})}{4} \tilde{B}_0\right)u-(X'+\tilde{J}\grad(\tr(\tilde{b})/4))\times u\,.
\end{eqnarray*} 
where $\grad$ is the gradient on $\tilde S$ with respect to $\tilde{I}_0$.
In particular, 
$$ X'''=-X'-\tilde{J}\grad(\tr(\tilde{b})/4)~,\qquad
2A^{X''}=-\tanh(s_0/2)\left(E-\frac{\tr(\tilde{b})}{2}\tilde{b}\right)~. $$
Replacing these identities in the right hand side of (\ref{rr:eq}),
we deduce
\begin{eqnarray*}
2(\tilde{J}A^{X''}(\tilde{u})-\langle X'+X''',\tilde{\nu}\rangle \tilde{u}-A^{X'}\circ \tilde{B}_0(\tilde{u})) & = & 
-\tanh(s_0/2)(\tilde{J}-\frac{\tr(\tilde{b})}{2}\tilde{J}\tilde{b})(\tilde{u})+\tanh(s_0/2)\frac{(\tilde{J}\tilde{b}-\tilde{b}\tilde{J})}{2}\tilde{b}(\tilde{u}) \\
& = & -\tanh(s_0/2)(\tilde{J}-\frac{\tr \tilde{b}}{2}\tilde{J}\tilde{b}-\frac{1}{2}\tilde{J}\tilde{b}^2+\frac{\tilde{J}}{2})(\tilde{u})
\end{eqnarray*} 
Using the identity $\tilde{b}^2=\tr(\tilde{b})\tilde{b}-E$,
we obtain that the right hand side in (\ref{rr:eq}) is equal to
\begin{equation}\label{rrr:eq}
\tanh(s_0/2)(\tilde{J}-(\tr(\tilde{b})/2)\tilde{J}\tilde{b})
\end{equation}
On the other hand, (\ref{eq:Bt}) shows that the left hand side of (\ref{rr:eq}) is equal to
\begin{equation}\label{rrrr:eq}
\frac{\tanh(s_0/2)}{2}(\tilde{b}\tilde{J}\tilde{b}-\tilde{J}\tilde{b}^2)=\frac{\tanh(s_0/2)}{2}(\tilde{J}-\tilde{J}(\tr(\tilde{b})\tilde{b}-E))
= \frac{\tanh(s_0/2)}{2}(2\tilde{J}-\tr(\tilde{b})\tilde{J}\tilde{b})\,.
\end{equation}
Equation (\ref{rr:eq}) follows by comparing (\ref{rrr:eq}) with (\ref{rrrr:eq}).
\end{proof}

\subsection{The comparison of the cocycles}
\label{ssc:cocycles}

Any element of $K\in \psl_2(\mathbb C)$ can be regarded as a Killing vector field on $\Hyp^3$.
 Notice that by definition of Killing vector field, there is another field $K'$ associated to $K$ such that
 \[
     D_{\tilde{v}} K=K'\times \tilde{v}\,.
 \]
 Using the same argument as in Proposition \ref{Zfield:prop}, one can check that $K'$ is a Killing
vector field and in fact $K''=-K$. More precisely, we have the following lemma:
\begin{lemma}\label{complex:lm}
As elements of $\psl_2(\mathbb C)$ we have $K'=iK$.
\begin{flushright}$\Box$\end{flushright}
\end{lemma}

Given a point $\tilde{x}\in \Hyp^3$ we have a natural map
\[
   \mathrm{ev}_{\tilde{x}}: \psl_2(\mathbb C)\ni K
\mapsto (K(\tilde{x}), K'(\tilde{x}))\in 
T_{\tilde{x}}\Hyp^3\oplus T_{\tilde{x}}\Hyp^3\,.
\]
It is a well-known fact that such a map is an isomorphism for every
$\tilde{x}\in \Hyp^3$.

Because of Lemma \ref{complex:lm},
if $\mathrm{ev}_{\tilde{x}}(K)=(\tilde{w}_1,\tilde{w}_2)$, then
$\mathrm{ev}_{\tilde{x}}(iK)=(\tilde{w}_2,-\tilde{w}_1)$.

Given any section $V$ on $\Theta$,  we define
\[
   K^V:\tilde S\rightarrow \psl_2(\mathbb C)
\]
such that $\mathrm{ev}_{\tilde{x}}(K^V(\tilde{p}))=V(\tilde{p})$ and 
$\mathrm{ev}_{\tilde{x}}((K^V)'(\tilde{p}))=V'(\tilde{p})$ for every
$\tilde{x}=\sigma_{s_0}(\tilde{p})$.

In particular,
we have maps $K^X$ and $K^Y$ associated
to the fields $X,Y$ defined in (\ref{eq:XY}),
so that $\mathrm{ev}_{\tilde{x}}(K^X)=(X,X')$ and
$\mathrm{ev}_{\tilde{x}}(K^Y)=(-X',-X'')$ by Proposition \ref{fieldsmain:prop}.
We conclude by Proposition \ref{Zfield:prop} that
\[
    K^Y=-iK^X-K_0
\]
where $K_0(\tilde{p})=\mathrm{ev}_{\tilde{x}}^{-1}(0,X''+X)=
\mathrm{ev}_{\tilde{x}}^{-1}(0,\frac{\tr(\tilde{b})}{4}\tilde{\nu})$
for all $\tilde{x}=\sigma_{s_0}(\tilde{p})$.
Since the field $\frac{\tr(\tilde{b})}{4}\tilde{\nu}$
is invariant under the action of $\pi_1(S)$, 
we find that $K_0$ is equivariant, that is
\[
   K_0(\gamma \tilde{p})=Ad(\rho(\gamma))K_0(\tilde{p})\,.
\]
However, it follows from Equation (\ref{field:eq}) that
\begin{eqnarray}
   K^X(\gamma \tilde{p})-Ad(\rho(\gamma))K^X(\tilde{p})=\eta(\gamma)\\
   K^Y(\gamma \tilde{p})-Ad(\rho(\gamma))K^Y(\tilde{p})=\xi(\gamma)
\end{eqnarray}
and by these equations and identities (\ref{field:eq}) and (\ref{field2}) one deduces that
\[
     \xi(\gamma)=-i\eta(\gamma)-(Ad(\rho(\gamma))-1)K_0(\tilde{p})=-i\eta(\gamma).
\]
Thus (\ref{hol:eq}) is proved.

\subsection{Parameterization by the disk}
\label{ssc:another}

The parameterization of the complex landslide used above is well-suited
for a comparison with the complex earthquake.
However, another parameterization -- already used
in the introduction -- is perhaps more
convenient when considering the holomorphic disks in Teichm\"uller space 
obtained as the image of the complex flow. This new parameter 
$\zeta$ takes values in the unit disk. We develop here the relationship between
these two parameterizations and we investigate the regularity at $\zeta=0$.

Consider $z=t+is$ in the upper half-plane,
so that $t\in \R$ and $s\geq 0$, and we set $\zeta=\exp(iz)=\exp(-s+it)$, which
belongs to the punctured closed unit disk $\dot{\overline{\Delta}}=
\{\zeta\in\mathbb{C}\,|\,0<|\zeta|\leq 1\}$.
For $h,h^\dual\in \cT$ we then define
$$ P_{\zeta}(h,h^\dual) := P'_{t+is}(h,h^\dual)~\qquad C_{\zeta}(h, h^\dual) := C'_{t+is}(h,h^\dual)~. $$
These maps are well-defined since $P'$ and $C'$ are invariant under
$t\mapsto t+2\pi$. Clearly, the maps $\zeta\mapsto P_{\zeta}(h,h^\dual)$ and
$\zeta\mapsto C_{\zeta}(h, h^\dual)$ 
are holomorphic in the unit disk minus its center, for any fixed $h$ and $h^\dual$. 

We first give an explicit formula for $C'_{is}(h,h^\dual)$. 

\begin{lemma} \label{lm:Cprime}
Let $h,h^\dual\in \cT$ be two hyperbolic metrics on $S$, and let $b$ be the bundle morphism
appearing in Corollary \ref{cr:b}. For every $s\in \R_{\geq 0}$, 
$$ C'_{is}(h,h^\dual) = h(\gamma_s \bullet,\gamma_s\bullet) $$
as conformal structures, where 
$$ \gamma_s = \cosh(s/2)E+\sinh(s/2)b~. $$
\end{lemma}

\begin{proof}
By definition, $C'_{is}(h,h^\dual)$ is the conformal structure at infinity of the (unique) hyperbolic
end containing a convex surface with induced metric $I=\cosh^2(s/2)h$ and third fundamental form 
$\III=\sinh^2(s/2) h^\dual$. Its shape operator is then $B=\tanh(s/2)b$, and conformal structure at
infinity of the end is given (see e.g. \cite{horo}) by
$$ C'_{is}(h,h^\dual) = I((E+B)\bullet,(E+B)\bullet) = h(\gamma_s\bullet,\gamma_s\bullet)~. $$
\end{proof}

We can now give a general formula for $C'_{t+is}(h,h^\dual)$ for
$t+is\in \overline\HH$.

\begin{lemma} \label{lm:other}
Let $s,t\in \R$ with $s\geq 0$, and let $\zeta=\exp(-s+it)$. 
Then, for all $h,h^\dual\in \cT$,
$$ C_{\zeta}(h,h^\dual) = h\left(B^\#_\zeta \cbull, B^\#_\zeta \cbull\right)~, $$
where 
$$ B^\#_\zeta = \frac{\zeta+1}{2\sqrt{\zeta}} E - \frac{\zeta-1}{2\sqrt{\zeta}} b $$
and $\sqrt{\zeta}$ is a notation for $\exp((-s+it)/2)$.
\end{lemma}

Here we use the convention that the complex number $i$ acts as the 
complex structure $J$ on tangent vectors.

\begin{proof}
It follows from the definition of $C'$ and from the previous lemma that 
$$ C'_{t+is}(h,h^\dual) = h(\beta_{-t}\circ \gammab_s\bullet,\beta_{-t}\circ \gammab_s\bullet)~, $$
where (as in Section 3):
$$ \beta_t = \cos(t/2) E + \sin(t/2) Jb~, ~~
 \gammab_s = \cosh(s/2)E+\sinh(s/2)b_{-t}~,~~b_{-t} =  \beta_{t}\circ b\circ \beta_{-t}~. $$
It is then clear that 
$$ \beta_{-t}\circ \gammab_s = \gamma_s\circ \beta_{-t}~, $$
so that
$$ C'_{t+is}(h,h^\dual) = h(B^\#_\zeta\cbull,B^\#_\zeta\cbull)~, $$
where $B^\#_\zeta=\gamma_s\circ \beta_{-t}$. 

Using the fact that $bJbJ=-E$:
\begin{eqnarray*}
B^\#_\zeta & = & (\cosh(s/2)E + \sinh(s/2) b)\circ (\cos(t/2) E - \sin(t/2) Jb)\\
& = & (\cosh(s/2)\cos(t/2) E - \sinh(s/2)\sin(t/2)bJb) + 
(\cos(t/2)\sinh(s/2)E-\cosh(s/2)\sin(t/2)J)b \\
& = & (\cosh(s/2)\cos(t/2) E - \sinh(s/2)\sin(t/2)J)  + 
(\cos(t/2)\sinh(s/2)E-\cosh(s/2)\sin(t/2)J)b \\
& = & \cosh((-s+it)/2)E - \sinh((-s+it)/2)b~.
\end{eqnarray*}
Setting $\sqrt{\zeta}=\exp((-s+it)/2)$, we can write this relation as 
$$ B^\#_\zeta = \frac{\sqrt{\zeta}+1/\sqrt{\zeta}}{2} E - \frac{\sqrt{\zeta}-1/\sqrt{\zeta}}{2}b~. $$
\end{proof}

It follows from the definitions that $C$ is essentially the same as $C'$ with a different 
parameterization. 
The main properties of this map are as follows.

\begin{prop}\label{prop:disc-filled}
Let $h,h^\dual\in \cT$ and let $c$ be the ``center'' of $(h,h^\dual)$ as defined
in Section \ref{ssc:center}. Then:
\begin{enumerate}
\item $C_{\zeta}(h,h^\dual)$ is defined for all $\zeta\in \dot{\overline{\Delta}}$,
\item it is holomorphic in $\zeta$,
\item it extends continuously, and therefore holomorphically, at $\zeta=0$, with 
$C_{0}(h,h^\dual)=c$. 
\item it also extends holomorphically to the open disk of center $0$ and radius 
$(\kappa_0+1)/(\kappa_0-1)$, where $\kappa_0=\max_{x\in S}\kappa(x)$ and
$\kappa:S\rightarrow[1,\infty)$ is the bigger eigenvalue 
of the operator $b$ associated to
the minimal Lagrangian map isotopic to the identity between
$(S,h)$ and $(S,h^\dual)$ (see Section \ref{ssc:minilag}).
\end{enumerate}
\end{prop}

In particular, $c$ appears as a smooth point of the holomorphic disk 
defined by $C$, while it was obtained only in the limit $s\rightarrow \infty$
in the parameterization used by $C'$.

\begin{remark}
Unlike the map $C_\cbull(h,h^\dual)$, the map $P_\cbull(h,h^\dual)$
does not extend at $\zeta=0$.
Indeed, take any sequence of positive real numbers $\zeta_n\rightarrow 0$.
By definition of the map $P$, there is an embedding of $S$
into the hyperbolic end $M_n$, with first fundamental form equal to 
$I_n=\cosh^2(-\frac{1}{2}\log\zeta_n)h$ and shape operator 
$B_n=-\tanh(-\frac{1}{2}\log\zeta_n)b$,
which corresponds to the projective structure $P_n$.
In particular, $B_n$ converges to $-b$.
On the other hand, by Proposition 4.2 of \cite{LL},
if $P_n$ converges to a projective surface,
$B_n$ should converge to the identity.
\end{remark}

\begin{proof}[Proof of Proposition \ref{prop:disc-filled}]
The first two points are direct consequences of the definition of $C$ from
$C'$, and of Theorem \ref{tm:complex}. The third point follows from the expression
of $B^\#_\zeta$ in Lemma \ref{lm:other}, because $C_{\zeta}(h,h^\dual)$ is really considered
as a conformal structure, so it is not changed if we multiply $B^\#_\zeta$ by a
complex-valued function defined on $S$. In particular, we can multiply 
$B^\#_\zeta$ by $2\sqrt{\zeta}$, obtaining 
$$ 2\sqrt{\zeta}B^\#_\zeta = (1+\zeta)E+(1-\zeta)b~, $$
which is clearly continuous at $u=0$. 

For the last point note that the expression defining $C_{\zeta}(h,h^\dual)$ in
Lemma \ref{lm:other} can be analytically continued if $B^\#_\zeta$ is
non-singular at all points of $S$.
This happens if
$$ \frac{\zeta+1}{2\sqrt{\zeta}} + \frac{1-\zeta}{2\sqrt{\zeta}}\kappa \neq 0 $$
everywhere on $S$, which is certainly satisfied if 
$$ |\zeta| < \frac{\kappa_0+1}{\kappa_0-1}~. $$
\end{proof}

\section{The earthquake flow as a limit} \label{sc:limit}

The main goal of this section is to prove Theorem \ref{tm:limit}. The arguments
are based on comparing surfaces in $\AdS_3$ with constant Gauss curvature close
to $-1$ to pleated surfaces. The key step in the proof of Theorem \ref{tm:limit} will be
Theorem \ref{tm:boundary}.

We fix a hyperbolic metric $h$ on $S$ and a
divergent sequence of metrics $h_n^\dual\in\cT$.  We will study the
asymptotic behavior of the holomorphic map $z\mapsto P_{z}(h,h_n^\dual)$
assuming that $(h_n^\dual)_{n\in \N}$ converges
to a point in the Thurston boundary of $\mathcal T(S)$ which is the
projective class of some measured geodesic lamination $\lambda$ on
$S$.

Take any sequence $\theta_n>0$ with $\lim_{n\rightarrow\infty}\theta_n=0$ such that
$\theta_n\ell_{h_n^\dual}\rightarrow\iota(\lambda,\cbull)$,
i.e. for every free homotopy class
$\gamma$ of  closed curves of $S$, the $h_n^\dual$-length of the
$h_n^\dual$-geodesic representative of $\gamma$ rescaled by the factor $\theta_n$
converges to the intersection between $\gamma$ and $\lambda$.  Define
$P'_{n}:\overline{\mathbb H}\rightarrow\mathcal P(S)$ as
$P'_{n}(z)=P'_{\theta_n z}(h,h_n^\dual)$, which are holomorphic, and let
$P'_{\infty}:\overline{\mathbb H}\rightarrow\mathcal P(S)$ be
$P'_{\infty}(t+is)=Gr_{s\lambda/2}(E_{-t\lambda/2}(h))$.

\begin{theorem}\label{conv:thm}
For every $z\in\overline{\mathbb{H}}$, we have that $P'_{n}(z)
\rightarrow P'_{\infty}(z)$.
\end{theorem}

Notice that since the $P'_{n}$ are holomorphic, the convergence
$P'_{n}\rightarrow P'_{\infty}$ is, in fact, in $C^\infty$. 
Note also that in this section we use the parameterization by the upper
half-plane -- which is more practical when considering the limit to 
complex earthquakes -- so that we use the notations with primes for
$L, P, SGr$, etc.

\subsection{Convergence on the imaginary axis}
\label{ssc:imaginary}

In this subsection we prove that if
$\theta_n\ell_{h_n^\dual}\rightarrow\iota(\lambda,\cbull)$,
then
\begin{equation}\label{convim:eq}
P'_{i\theta_n}(h,h_n^\dual)\rightarrow Gr_{\lambda/2}(h)\,.
\end{equation}
In fact, in order to prove Theorem \ref{conv:thm}, we will need
to prove that the convergence is uniform with respect to $h$.

\begin{prop}\label{im:prop}
Let $(h_n)_{n\in \N}$ be a sequence of hyperbolic metrics converging to a
hyperbolic metric $h$ on $S$, and let $(h_n^\dual)_{n\in \N}$ be a sequence
of hyperbolic metrics converging projectively to $[\lambda]$ in
the Thurston boundary of $\mathcal T$.
If $\theta_n$ is a sequence of positive numbers such that
$\theta_n\ell_{h_n^\dual}\rightarrow\iota(\lambda,\cbull)$,
then $SGr'_{\theta_n}(h_n,h_n^\dual)$ converges to
$Gr_{\lambda/2}(h)$.
\end{prop}

Notice that (\ref{convim:eq}) corresponds to the particular
case of Proposition
\ref{im:prop} in which $h_n$ is constant.

By definition, $SGr'_{\theta_n}(h_n,h_n^\dual)=P'_{i\theta_n}(h_n,h_n^\dual)$ is the
projective structure on $S$ determined by
prescribing that the associated hyperbolic end $(M_n,g_{M_n})$ contains a
constant curvature surface $S_n$,
parametrized by $\bar\sigma_n:S\rightarrow S_n\subset M_n$,
with first fundamental form
$I_{S_n}=\cosh^2(\theta_n/2)h_n$ and third fundamental form
$\III_{S_n}=\sinh^2(\theta_n/2)h_n^\dual$.

In general, given an embedding $S\rightarrow\Sigma\subset X$
of $S$ inside a (hyperbolic, de Sitter or anti de Sitter) 3-manifold,
we will denote by the same symbol the first fundamental form
$I_\Sigma$ (resp. the third fundamental form $\III_\Sigma$) and
its pull-back on $S$.

Let $\partial M_n$ be the hyperbolic boundary of $M_n$, that carries a
hyperbolic induced metric $g_{\partial M_n}$ and is locally bent along a measured
geodesic lamination $\lambda_n$. 
By definition, $P'_{i\theta_n}(h_n,h_n^\dual)=Gr_{\lambda_n}(g_{\partial M_n})$. 
So, in order to prove Proposition
\ref{im:prop}, it is sufficient to check that 
$(S, g_{\partial{M}_n})$
converges to $(S,h)$, and $\lambda_n$ converges to $\lambda/2$
in $\mathcal {ML}(S)$.

\begin{lemma}\label{cpthyp:lm}
The hyperbolic metrics $g_{\partial M_n}$ are contained in a compact subset
of $\cT$.
\end{lemma}

\begin{proof}
The closest point projection $r_n:M_n\rightarrow \partial{M}_n$ is
$1$-Lipschitz. In particular, we have a $1$-Lipschitz map
$r_n|_{S_n}:S_n\rightarrow \partial{M}_n$.  This implies that the marked length
spectrum of $\partial{M}_n$ is bounded from above by the marked length spectrum
of $S_n$, that in turn is locally uniformly bounded.
 \end{proof}

\begin{lemma}\label{cptds:lm}
The bending laminations $\lambda_n$ are contained in a compact subset
of $\cML$.
\end{lemma}

In order to prove Lemma \ref{cptds:lm}, we will consider the de Sitter
spacetime $M_n^\dual$ dual to $M_n$: it is the set of complete geodesic
planes contained in $M_n$.
The de Sitter structure is induced by the
natural map $dev^\dual:\tilde{M}_n^\dual\rightarrow \dS^3$, where the model
of de Sitter geometry is the set of geodesic planes of $ \Hyp^3$
(see \cite{scannell, BeBo}).
 
 Scannell \cite{scannell} showed that $M_n^\dual$ is a globally hyperbolic
 spacetime diffeomorphic to $S\times\mathbb R$.  Following
 \cite{BeBo}, $\tilde M_n^\dual$ has a natural boundary corresponding to
 the set of support planes of $\partial{\tilde{M}}_n$. This boundary is called
 the initial singularity: the de Sitter metric extends to the boundary
 and makes it an achronal (but not spacelike) surface. So $\partial
 \tilde M_n^\dual$ carries a pseudometric $d_0$ induced by $\tilde M_n^\dual$. By
 \cite{BeBo, beguad} it turns out that the action of $\pi_1(S)$ on
 $\tilde M_n^\dual$ extends to the boundary (even if the action on the
 boundary is neither proper nor free) and the marked length spectrum
 of this action coincides with the intersection with the bending
 lamination:
 \begin{equation}\label{initsing:eq}
     \iota(\lambda_n, \gamma)=\inf\{d_0(x^\dual, \gamma x^\dual)|
x^\dual\in\partial \tilde M_n^\dual\}\,.
 \end{equation}
 
 Let $S_n^\dual$ be the surface in $M_n^\dual$ dual to $S_n$, corresponding to the
 set of support planes of $S_n$ (that by the convexity of $S_n$ are
 complete planes in $M_n$). There is a natural map $S_n\rightarrow S_n^\dual$
 sending $x$ to the dual of the plane tangent to $S_n$ at $x$.

A simple local computation shows that the
 first fundamental form of $S_n^\dual$ coincides with the third
 fundamental form of $S_n$ (throught the natural map $S_n\rightarrow
 S_n^\dual$).  In particular $S_n^\dual$ is a surface of constant curvature
 $-1/\sh^2(\theta_n/2)$.  Barbot, Beguin and
Zeghib \cite{BBZ2} have shown that there is a
 time function (K-time) $\mathfrak{t}_n:M_n^\dual\rightarrow (-\infty,0)$ such that
 $\mathfrak{t}^{-1}_n(k)$ is the unique surface in $M_n^\dual$ with constant
 curvature $k$.
 
Lemma \ref{cptds:lm} is a simple consequence of (\ref{initsing:eq})
and the following general Lemma of de Sitter geometry.

 \begin{lemma}[\cite{mehdi}] \label{cptlam:lm}
 Let $M^\dual$ be a de Sitter MGHC spacetime and $\partial \tilde M^\dual$ be the
 boundary of its universal covering. 
 If $S^\dual$ is a constant curvature surface in $M^\dual$, there is a natural
 $1$-Lipschitz equivariant map
 \[
  \tilde{r}^\dual: \tilde S^\dual\rightarrow\partial \tilde{M}^\dual
 \]
 such that $\tilde{r}^\dual(\tilde{x})\in I^-(\tilde{x})\cap\partial
\tilde{M}^\dual$.
 
 In particular, we have
 \[
   \ell_{\partial M^\dual}(\gamma)\leq\ell_{S^\dual}(\gamma)
 \]
 where $\ell_{\partial M^\dual}$ and $\ell_{S^\dual}$ are the marked length
spectra of $\partial M^\dual$ and $S^\dual$ respectively.
 \end{lemma}

\begin{prop}\label{limithyp:lm}
The hyperbolic metrics $g_{\partial M_n}$ converge to $h$.
\end{prop}

\begin{proof}
By Lemmas \ref{cpthyp:lm} and \ref{cptds:lm}, we have that up to
passing to a subsequence $M_n\rightarrow M_\infty$.  In particular, we
can concretely realize $M_n\cong (S\times[0,+\infty),g_{M_n})$ in such a way that
 $g_{M_n}$ converges to a hyperbolic metric $g_{M_\infty}$ such that
  $(S\times[0,+\infty),g_{M_\infty})\cong M_\infty$, where we have denoted
by the same symbol the metric $g_{M_n}$ on $M_n$ and the corresponding
metric on the model $S\times[0,+\infty)$.
  
By abuse of notation,
we denote again by $r_n:S\times[0,+\infty)\rightarrow S\times\{0\}$ the
  retraction corresponding to the retraction of $M_n$ onto
$\partial M_n$.  Let
  $\bar{\sigma}_n:(S,h)\rightarrow (S\times[0,+\infty),g_{M_n})$ be 
the embedding
    with first fundamental form $I_{S_n}=\ch^2(\theta_n/2)h_n$
and third fundamental form $\III_{S_n}=\sh^2(\theta_n/2)h^\dual_n$.

Notice that the composition
$i_n=r_n\circ\bar{\sigma}_n:(S,\ch^2(\theta_n/2)h_n)\rightarrow
(S\times\{0\}, g_{\partial M_n})$
is a $1$-Lipschitz homotopy equivalence.  So $i_n$ converges (up to
passing to a subsequence) to a $1$-Lipschitz homotopy equivalence
$i_\infty:(S,h)\rightarrow(S\times\{0\},g_{\partial M_\infty})$.
Since both $h$ and
$g_{\partial M_\infty}$ are hyperbolic metrics, we conclude that $i_\infty$ is an
isometry.
\end{proof}

Let $\lambda_\infty$ the bending lamination of $M_\infty$.  In order to
conclude the proof of (\ref{convim:eq}) we need to show that
$\lambda_\infty=\lambda/2$. 
In fact,
the following general result in Lorentzian geometry
and (\ref{initsing:eq}) show that
\begin{equation}\label{convlam:eq}
  \iota(\lambda_\infty,\gamma)=\lim_n \ell_{S_n^\dual}(\gamma)=
\lim_n\ell_{\III_{S_n}}(\gamma)=\lim_n\frac{\theta_n}{2}\ell_{h_n^\dual}(\gamma)=
\iota(\lambda/2,\gamma)\,. 
\end{equation}
for every 
closed curve $\gamma$.

\begin{prop}\label{mehdi:prop}\cite{mehdi}
Let $(X_n^\dual)_{n\in\N}$ be a sequence of MGH de Sitter (or anti de Sitter)
spacetimes homeomorphic to $S\times\mathbb R$. Suppose
that $X_n^\dual$ converges to a MGH spacetime $X_\infty^\dual$.
Take any sequence
of numbers $k_n\rightarrow-\infty$ and let $\Sigma_n^\dual$ be the
future-convex surface of constant curvature $k_n$ contained in $X_n^\dual$.

Denote by $\ell_0$ the length spectrum of the initial singularity of
$X_\infty^\dual$. Then, for
every $\gamma\in\pi_1(S)$ we have
\[
   \ell_{\Sigma_n^\dual}(\gamma)\rightarrow\ell_0(\gamma)
\]
as $n\rightarrow+\infty$.
\end{prop}

In the next section we will give a short description of the initial
singularity for anti de Sitter spacetimes and
we will apply Proposition \ref{mehdi:prop} to this case.

\subsection{Convergence on the real axis}\label{re:ssec}

In this section we will prove that
\begin{equation}\label{convre:eq}
   L'^1_{\theta_n}(h, h_n^\dual)\rightarrow E_{\lambda/2}(h)\,
\end{equation}

Recall that $L'^1_{\theta_n}$ is the composition of the map $L'_{\theta_n}:\cT\times\cT
\rightarrow \cT\times \cT$ with the projection on the first factor. 

As in the previous section, we will need a slightly stronger statement.

\begin{prop}\label{re:prop}
Let $(h_n)_{n\in \N}$ and $(h_n^\dual)_{n\in \N}$ be two sequences of hyperbolic metrics 
such that $(h_n)_{n\in\N}$ converges to a hyperbolic metric $h$ on $S$ and that $(h_n^\dual)$ 
converges to a point $[\lambda]$ in Thurston boundary of $\mathcal T$.
If $\theta_n$ is a sequence of positive numbers such that
$\theta_n\ell_{h_n^\dual}\rightarrow\iota(\lambda,\cbull)$,
then $L'^1_{\theta_n}(h_n,h_n^\dual)$ converges to
$E_{\lambda/2}(h)$.
\end{prop}

Recall that the holonomy of $L'^1_{\theta_n}(h_n,h_n^\dual)$ corresponds to the left
holonomy of the MGH AdS $N_n\cong S\times\mathbb R$
containing a future-convex $K$-surface $F_n$ with $I_{F_n}=\cos^2(\theta_n/2)h_n$
and $\III_{F_n}=\sin^2(\theta_n/2)h_n^\dual$.  In order to prove Proposition
\ref{re:prop}, we will show that $N_n$ converges to a MGH AdS
structure $N_\infty$ and $F_n$ converges to the lower boundary
$\partial_- C(N_\infty)$
of the convex core of $N_\infty$. Then we will prove that
$\partial_- C(N_\infty)$ is
isometric to $(S,h)$ and is bent along a lamination corresponding to
$\lambda/2$. By a result of Mess \cite{mess}, the left holonomy of
$N_\infty$ (that is, by definition, the limit of the left holonomies of
$N_n$) is equal to the holonomy of $E_{\lambda/2}(h)$, and
Proposition \ref{re:prop} follows.

In order to prove that $N_n$ converges to some structure, we
will consider the lifting $\phi_n:\tilde S\rightarrow
\tilde{N}_n\subset \AdS^3$ corresponding to the embedding
$\bar\phi_n:S\rightarrow F_n\subset N_n$.  The map $\phi_n$
is determined up to isometry of $\AdS^3$ and we will normalize it
by requiring that, for some fixed $\tilde{p}_0\in\tilde S$,
$\phi_n(\tilde{p}_0)=\tilde{x}_0$ and
the normal vector to $\tilde F_n=\phi_n(\tilde{S})$ at $\tilde{x}_0$ is equal to
$\tilde{\nu}_0$ for some fixed $\tilde{x}_0,\tilde{\nu}_0$ in $\AdS^3$.

The first step to prove the convergence of $N_n$ is to show
that $\phi_n$ converges to a spacelike embedding into $\AdS^3$.

\begin{prop}\label{splklimit:prop}
Up to passing to a subsequence,
$\tilde{F}_n$ converges to a spacelike surface $\tilde{F}_\infty$  in
$\AdS^3$
and the map $\phi_n$ converges to an embedding
\[
\phi_\infty:\tilde S\rightarrow \AdS^3
\]
whose image is $\tilde{F}_\infty$.
\end{prop}

The easy part of the proof is to show that $\tilde{F}_n$ converges to an
embedded surface $\tilde{F}_\infty$ in $\AdS^3$ that is achronal. The main issue is to
show that the surface $\tilde{F}_\infty$ is spacelike.  The proof relies on the fact
that, for some fixed $\tilde{p}\in\tilde S$, the tangent planes of $\tilde{F}_n$ at
$\phi_n(\tilde{p})$ are uniformly spacelike, in the sense that they cannot
approximate lightlike planes.

The proof of this fact is based on the technical Lemma \ref{boundgr:lm}.

\begin{lemma}\label{boundgr:lm}
Let $b_n$ be the $h_n$-self adjoint operator such that $h^\dual_n=h_n(b_n\cbull,
b_n\cbull)$ and let $\bar\sigma_n:S\rightarrow S_n\subset M_n$ be the
embedding introduced in Section \ref{ssc:imaginary}.
Let $I_{\tilde{S}_n}^\#$ be the lifting to $\tilde S$ of the grafted metric $I_{S_n}^\#$ 
introduced in Definition \ref{df:grafted}.
Then, for every compact set $K\subset\tilde S$, there
is a constant $C_K$ such that the diameter of $K$ with respect to
$I_{S_n}^\#$ is bounded by $C_K$ for every $n$.
\end{lemma}

\begin{proof}
For any $k\in [-1,0)$, $M_n$ contains exactly one
  $K$-surface of constant curvature $k$, denoted here by $M_n(k)$
  (where by $M_n(-1)$ we mean the boundary of $M_n$).  For each $n$,
  let $G_n=SGr'_{\theta_n}(h_n,h_n^\dual)$ be the projective surface at
  infinity of $M_n$.  
  Let us consider the natural retraction
$\Pi_{M_n(k)}: G_n\rightarrow M_n(k)$, which is the limit
of the closest point projections $M_n(K)\rightarrow M_n(k)$
onto the convex surface $M_n(k)$
as $K>k$ converges to $0$ (see Figure \ref{hmap:fig}).

\begin{center}
\begin{figurehere}
\psfrag{x}{$x$}
\psfrag{Gn}{$G_n$}
\psfrag{Mn}{$M_n$}
\psfrag{Mnk}{$M_n(k)$}
\psfrag{Pik}{$\Pi_{M_n(k)}(x)$}
\psfrag{pMn}{$\partial M_n$}
\psfrag{Pi}{$\Pi_{M_n(-1)}(x)$}
\psfrag{ln}{$\lambda_n$}
\includegraphics[width=0.6\textwidth]{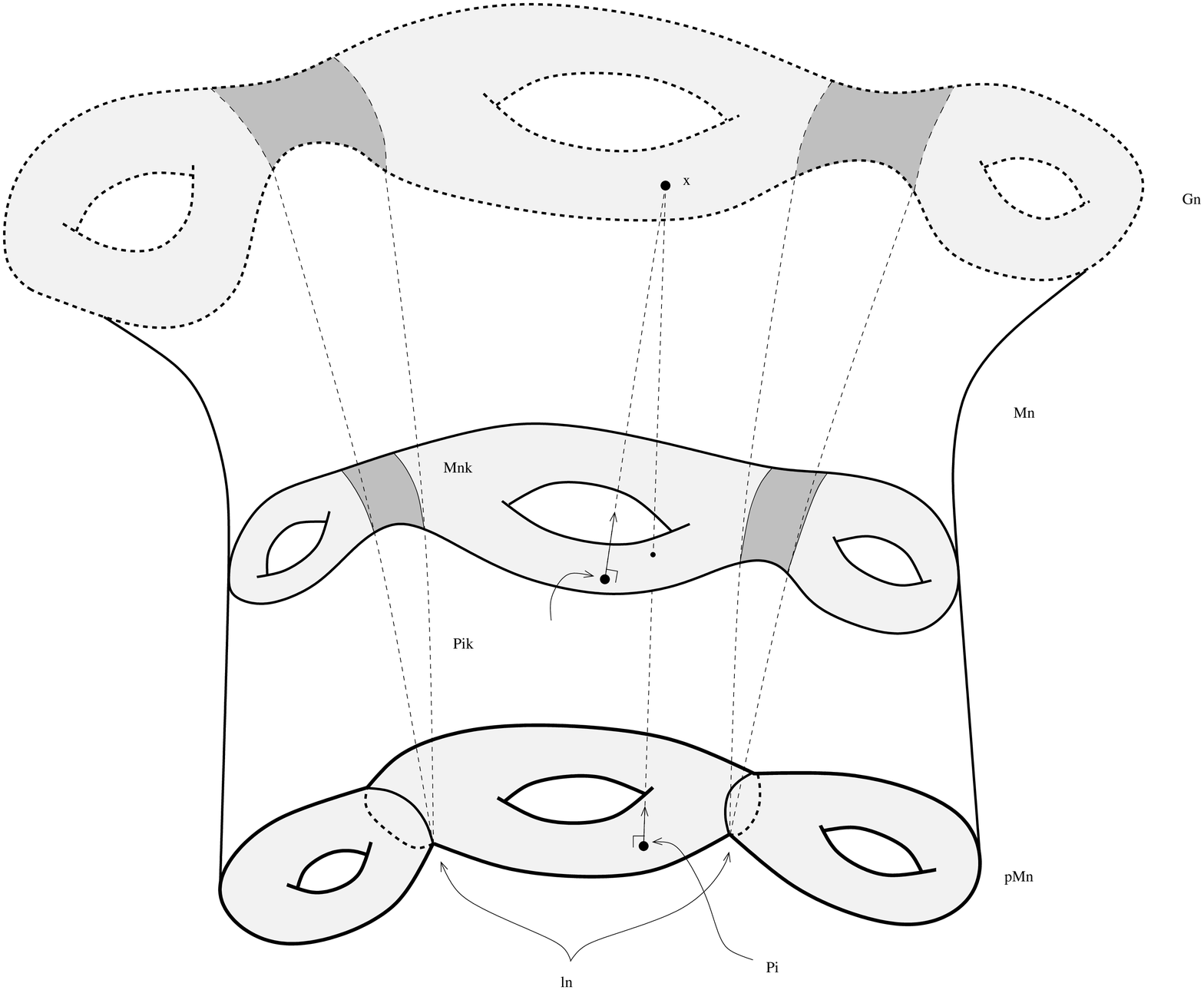}
\caption{{\small The retraction $\Pi_{M_n(k)}$.}}\label{hmap:fig}
\end{figurehere}
\end{center}

On the universal covering, $\Pi_{\tilde{M}_n(k)}$
sends a point $\tilde{x}\in\tilde{G}_n$ to the tangency point of the unique
horocycle centered at $\tilde{x}$ and tangent to $\tilde M_n(k)$.
For $k>-1$, $\Pi_{M_n(k)}$ is a diffeomorphism and the inverse is the map
obtained by sending each point of $y\in M_n(k)$ to the final point of
the geodesic ray starting from $y$ and orthogonal to $M_n(k)$.  If
$k=-1$, $\Pi_{M_n(-1)}$ is not injective in general,
since points on $M_n(-1)$ can admit several normal directions.
Nevertheless, $\Pi_{M_n(-1)}:G_n\rightarrow M_n(-1)$  is a
homotopy equivalence.  

In \cite{horo}, it has been showed that this diffeomorphism is conformal
with respect to the grafted metric $I^\#_{M_n(k)}$ of $M_n(k)$.
The conformal factor is an increasing function of $k$:
this precisely means that the conformal map
\[
   \Pi_{M_n(k')}\circ\Pi_{M_n(k)}^{-1}:(M_n(k), I^\#_{M_n(k)})\rightarrow
   (M_n(k'), I^\#_{M_n(k')})
\]
decreases the lengths when $k>k'$.

Now notice that $S_n$ is equal to $M_n(K_n)$ for $K_n=-1/\cosh^2(\theta_n)$.
As definitively $K_n<-1/2$, the map
\[
    j_n=\bar\sigma_n^{-1}\circ\Pi_{S_n}\circ \Pi_{M_n(-1/2)}^{-1}:(M_n(-1/2),
    I^\#_{M_n(-1/2)})\rightarrow (S, I^\#_{S_n})
\]
decreases the lengths.

Since $M_n$ converges to an hyperbolic end $M_\infty$,
the surface $M_n(-1/2)$ converges to $M_\infty(-1/2)$ in $C^\infty$-sense.
This means that $M_n$ can be  concretely realized  as a hyperbolic metric
$g_{M_n}$
on $S\times[0,+\infty)$ such that $M_n(-1/2)=S\times\{1\}$ and such that
$g_{M_n}$ converges to a hyperbolic metric $g_{M_\infty}$ and
$M_\infty(-1/2)=S\times\{1\}$.
Then the family of $1$-Lipschitz maps $j_n$ converges
to the map $j_\infty=\bar\sigma_\infty^{-1}\circ\Pi_{S_\infty}
\circ\Pi_{M_\infty(-1/2)}^{-1}:M_\infty(-1/2)\rightarrow S$,
which is a homotopy equivalence.

Let $\tilde j_\infty$ and $\tilde j_n$ be the lifting of those maps to the
universal covering.  Notice that $\tilde j_\infty$ is a proper map.  If $K$
is a compact subset of $\tilde S$, then
$K'=\tilde j_\infty^{-1}(K)$ is a compact subset of
$\tilde S\times\{1\}$, and $K'_n=\tilde j_n^{-1}(K)$ is contained in
some compact neighborhood of $K'$ for every $n$. In particular, there
exists a constant $C'_K$ such that the diameter of $K'_n$ with respect
to $\tilde I^\#_{M_n(-1/2)}$ is bounded by $C'_K$ for all $n$.
Taking $C_K$ bigger than $C'_K$, it
follows that the diameter of every $K'_n$ with respect to
$\tilde I^\#_{M_n(-1/2)}$ is bounded by $C_K$. Since
$\tilde{j}_n$ decreases the lengths, we have that the diameter of $K$ with
respect to $I_{\tilde{S}_n}^\#$ is bounded by $C_K$ for $n$ large enough.
\end{proof}

\begin{lemma}\label{unifsplk:lm}
For every $d>0$ there is a compact set $K$ in $\AdS^3$ such that for
$\tilde{p}\in\tilde S$ with $d_{h}(\tilde{p},\tilde{p}_0)<d$,
the normal vector $\tilde{\nu}_n(\tilde{p})$ of
$\tilde{F}_n$ at $\phi_n(\tilde{p})$ lies in $K$.
\end{lemma}

\begin{proof}
Lemma \ref{boundgr:lm} implies that for any $d>0$, there is $D>0$
such that for any $n$ and any
$\tilde{p}\in B_{\tilde{h}}(\tilde{p}_0, d)$   there exists a path
$\tilde\varsigma:[0,1]\rightarrow \tilde S$ connecting
$\tilde{p}_0$ to $\tilde{p}$ such that $\ell_{I_{\tilde{S}_n}^\#}(\tilde\varsigma)$
is bounded by $D$.

We claim (and will prove below) that this implies that 
\begin{eqnarray}
    |\langle \tilde{x}_0,\tilde{\nu}_n(\tilde{p}) \rangle| & \leq & 2e^{2D} \label{fstcptest:eq}\\
    |\langle \tilde{\nu}_0,\tilde{\nu}_n(\tilde{p})\rangle| & \leq & 2e^{2D}\,.\label{scdcptest:eq}
\end{eqnarray}
It follows from this claim that $\tilde{\nu}_n(\tilde{p})$ is contained in 
$$ K=\{w\in \AdS^3|\langle \tilde{x}_0, w\rangle\leq 2e^{2D}
\qquad\textrm{ and }\langle \tilde{\nu}_0, w\rangle\leq 2e^{2D}\} $$
which is a compact subset of $\AdS^3$, and the lemma follows. We now
turn to the proof of the claim.

We fix $n$ and consider the following functions:
\begin{eqnarray*}
a(t)=-\langle \tilde{x}_0, \phi_n(\tilde\varsigma(t))\rangle, &  &
a_\perp(t)=-\langle \tilde{x}_0, \tilde{\nu}_n(\tilde\varsigma(t))\rangle~, \\ 
a^\dual(t)=-\langle \tilde{\nu}_0, \phi_n(\tilde\varsigma(t))\rangle, &  &
a_\perp^\dual(t)=-\langle \tilde{\nu}_0, \tilde{\nu}_n(\tilde\varsigma(t))\rangle~.
\end{eqnarray*}
Notice that $a$ is a positive function since $\tilde{x}_0$ and
$\phi_n(\tilde\varsigma(t))$ are
contained in a spacelike surface. Moreover, since the surface $\tilde
F_n$ is convex, the plane orthogonal to
$\tilde{\nu}_n(\tilde\varsigma(t))$ is a support
plane for $\tilde F_n$, so it is not difficult to check that also $a_\perp$ is
positive (see Figure \ref{supp:fig}).

\begin{center}
\begin{figurehere}
\psfrag{x0}{$\tilde{x}_0$}
\psfrag{nu}{$\tilde{\nu}_n(\tilde\varsigma(t))$}
\psfrag{phi}{$\phi_n(\tilde\varsigma(t))$}
\psfrag{<}{$\langle \tilde{x},\tilde{\nu}_n(\tilde{\varsigma}(t))\rangle<0$}
\psfrag{>}{$\langle \tilde{x},\tilde{\nu}_n(\tilde{\varsigma}(t))\rangle>0$}
\psfrag{Fn}{$\tilde{F}_n$}
\psfrag{si}{$\phi_n(\tilde\varsigma)$}
\psfrag{AdS}{$\AdS^3$}
\includegraphics[width=0.8\textwidth]{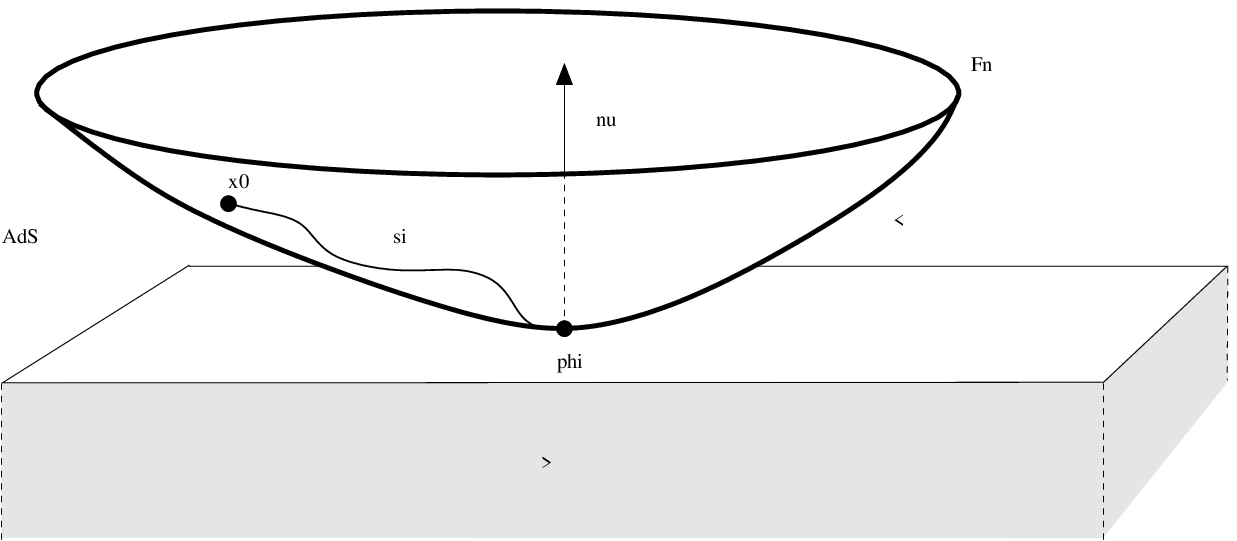}
\caption{{\small The product $\langle \tilde{x}_0,
    \tilde\nu_n(\tilde\varsigma(t))\rangle$ is
    negative.}}\label{supp:fig}
\end{figurehere}
\end{center}

We can decompose $\tilde{x}_0$ as
\[
   \tilde{x}_0=a(t)\tilde\phi_n(\tilde\varsigma(t))+
a_\perp(t)\tilde{\nu}_n(\tilde{\varsigma}(t))
+\tilde v(t)
\]
with $\tilde{v}(t)\in T_{\phi_n(\tilde\varsigma(t))}\tilde F_n$.
Imposing $\langle \tilde{x}_0, \tilde{x}_0\rangle=-1$, we deduce that
$\|\tilde{v}(t)\|\leq a+a_\perp$.

On the other hand,
$\|\tilde{v}(t)\|\leq \|\tilde{v}(t)\|_{\tilde{h}_n}$
and so
\[
\begin{array}{ccl}
\dot a & =&
\langle \tilde{x}_0,\dot \phi_n\circ\tilde\varsigma\rangle=
\langle \tilde{v}, d\phi_n(\dot{\tilde{\varsigma}}) \rangle\leq
(a+a_\perp)\|\dot{\tilde{\varsigma}}\|_{\tilde h_n}\\
\dot a_\perp & =&\langle \tilde{x}_0, B_{\tilde{F}_n}(d\phi_n(\dot{\tilde{\varsigma}}))
\rangle\leq (a+a_\perp)\| B_{\tilde{F}_n}(\dot{\tilde{\varsigma}})\|_{\tilde h_n}
\leq(a+a_\perp)\frac{\theta_n}{2}\|\dot{\tilde{\varsigma}}\|_{\tilde{h}_n^\dual}\,.
\end{array}
\]
Since $I_{S_n}^\#$ dominates both $\tilde h_n$ and
$\frac{\theta_n^2}{4} \tilde h_n^\dual$ we see that
\[
\begin{cases}
(a+a_\perp)(0)=2 \\
\dot{a}+\dot{a}_\perp \leq 2(a+a_\perp)\|\dot{\tilde{\varsigma}} \|_{I_{\tilde{S}_n}^\#}
\end{cases}
\]
and by a simple integration we have $a+a_\perp\leq 2e^{2D}$.

A similar argument can be applied to $a^\dual$ and $a_\perp^\dual$ using the fact
that the path $\tilde{\nu}_n(\tilde\varsigma(t))$
is contained in the dual surface $\tilde F_n^\dual$,
that is the surface made of normal vectors of $\tilde F_n$. 
Indeed, there is a  natural map
$\phi_n^\dual:\tilde S\rightarrow \tilde F_n^\dual$
that sends a point $\tilde{p}$ to the point dual to
the plane tangent to $\tilde F_n$ at $\phi_n(\tilde{p})$.
The corresponding embedding data are
\begin{equation}\label{dlemb:eq}
I_{\tilde{F}_n^\dual}=
\III_{\tilde{F}_n}=
\sin^2(\theta_n/2)h^\dual_n~, \qquad
B_{\tilde{F}_n^\dual}=-1/\tan(\theta_n/2)b_n^{-1}\,.
\end{equation}
In particular,
$F_n^\dual$ is a past-convex spacelike surface, and the previous  argument
shows that $a^\dual+a_\perp^\dual\leq 2e^{2D}$.
\end{proof}

\begin{proof}[Proof of Proposition \ref{splklimit:prop}]
We will consider the product model of $\AdS^3= \Hyp^2\times S^1$,
where the metric  at some point $(\xi, e^{i\vartheta})$ is 
\[
   g_{\AdS^3}=g_{ \Hyp^2}-\chi(\xi)d\vartheta
\]
where $\chi(\xi)=\ch d_{ \Hyp^2}(\xi, \xi_0)$, where $\xi_0$ is
some fixed point (see \cite{maximal}).

By a lemma of Mess \cite{mess},
the image  of $\phi_n$ is the graph of some function
$ \Hyp^2\ni\xi\mapsto e^{is_n(\xi)}\in S^1$ 
that satisfies the following spacelike
condition
\begin{equation}\label{splkcond:eq}
    \| \grad (s_n)\| <1/\chi\,.
\end{equation}

We can also suppose that $\phi_n(\tilde{p}_0)$
is the point $\tilde{x}_0=(\xi_0,0)$ and the
normal vector of $\tilde F_n$ at $\tilde{x}_0$ is the unit vertical vector.

By (\ref{splkcond:eq}), the functions
$s_n$ are uniformly Lipschitz on compact sets
of $ \Hyp^2$.  So, up to subsequence, $\tilde F_n$ converges to a
surface $\tilde F_\infty$
which is the graph of some limit function $s_\infty$, that verifies
$\| \grad(s_\infty)\| \leq 1/\chi$ almost every-where.

In order to show that $\tilde{F}_\infty$ is spacelike, we need to prove
that $s_\infty$
verifies the strict inequality (\ref{splkcond:eq}) almost
everywhere.  Notice that the projection $\pi_n:\tilde{F}_n\rightarrow
\Hyp^2$ increases the length, so the disk $D$ in $ \Hyp^2$ with
center $(\xi_0,0)$ and radius $r$ is contained in
$\pi_n\circ\phi_n(B_{\tilde{F}_n}(\tilde{x}_0,r))$.
By Lemma \ref{unifsplk:lm}, the normal
vectors of $\tilde{F}_n$ on the cylinder based on $D$ are contained in some
compact subset $K$ (independent of $n$).

Since the normal vector at $(\xi,s_n(\xi))$ is the vector
\[
\frac{1}{\sqrt{1-\chi^2\| \grad(s_n)\| ^2}}\left(\grad(s_n)+
\frac{\partial\,}{\partial\vartheta}\right)
\]
under the natural identification $T( \Hyp^2\times S^1)=T
\Hyp^2\oplus TS^1$, we deduce that there exists $\epsilon$ depending on
$K$, such that
\[
  \| \grad(s_n)\| \leq (1-\epsilon)/\chi
\]
for every $\xi\in D$ and every $n$. This shows that
$\tilde{F}_\infty$ is spacelike.

Moreover, the restriction of the projection $\pi_n\circ\phi_n:(\tilde S,
h)\rightarrow  \Hyp^2$ on $B_{\tilde{h}}(\tilde{p}_0,r)$ is
$C$-Lipschitz, for some constant $C$ depending only on
$r$.  Indeed, given a vector $\tilde{v}\in T_{\tilde{p}}\tilde S$, let
$\tilde{v}_n=d\phi_n(v)$ and $\tilde{u}_n=d\pi_n(\tilde{v}_n)$.
We have that $\tilde{v}_n=\tilde{u}_n+\langle
\grad(s_n), \tilde{u}_n\rangle \frac{\partial\,}{\partial\vartheta}$, so
$\cos^2(\theta_n/2)\tilde{h}_n(\tilde{v},\tilde{v})=
\langle \tilde{v}_n, \tilde{v}_n\rangle\geq
\| \tilde{u}_n\| ^2-\chi^2\| \grad(s_n)\| ^2\| \tilde{u}_n\| ^2
\geq \epsilon\| \tilde{u}_n\| ^2$.

Since $h_n\rightarrow h$, there is $C'$ such that the identity map between
$(S,h)$ and $(S,h_n)$ is $C'$-Lipschitz for every $n$. 
It follows, after taking a subsequence, $(\pi_n\circ\phi_n)$ converges to a map 
$\pi'_\infty:\tilde S\rightarrow \Hyp^2$,
so $(\phi_n)$ converges to the map
$\phi_\infty(\tilde{p})=(\pi'_\infty(\tilde{p}),
s_\infty(\pi'_\infty(\tilde{p}))$.
\end{proof}

We can prove now that the holonomy
$\rho_n:\pi_1(S)\rightarrow Isom_0(\AdS^3)$
of $N_n$ converges to a limit representation
$\rho_\infty$ for which $\phi_\infty$ equivariant.

\begin{lemma}\label{limitsptm:lm}
If $\phi_n$ converges to a space-like embedding $\phi_\infty$, then  
the representation
$\rho_n$ converges to a representation 
$\rho_\infty:\pi_1(S)\rightarrow Isom_0(\AdS^3)$
such that $\tilde{F}_\infty$ is $\rho_\infty$-equivariant.

Moreover, the left and right components of $\rho_\infty$
are discrete and faithful
representations of $\pi_1(S)$ into $\PSL_2(\mathbb R)$.
\end{lemma}

\begin{proof}
First we prove that, for every $\gamma\in\pi_1(S)$,
the sequence $\rho_n(\gamma)$ is bounded in $Isom_0(\AdS^3)$.

Recall that we are assuming that $\phi_n(\tilde{p}_0)=\tilde{x}_0$
for all $n$ and that
the normal vectors $\tilde{\nu}_n(\tilde{p}_0)$
are equal to $\tilde{\nu}_0$.
Now the $\rho_n(\gamma)(\tilde{x}_0)=\phi_n(\gamma \tilde{p}_0)$
form a sequence converging to
$\bar x_0=\phi_\infty(\gamma \tilde{p}_0)$
and the $\rho_n(\gamma)(\tilde{\nu}_0)$ converge to a
unit timelike vector $\bar\nu_0$ at $\bar x_0$ orthogonal to some
support plane of $\tilde{F}_\infty$.

This implies there is a bounded sequence of isometry of $\AdS^3$, says
$\tau_n$, such that
\[
    \tau_n \rho_n(\gamma)(\tilde{x}_0)=\tilde{x}_0~,
\qquad \tau_n\rho_n(\gamma)(\tilde{\nu}_0)=\tilde{\nu}_0
\]
Now the set of isometries that fix $\tilde{x}_0$ and $\tilde{\nu}_0$
is compact; so,
after taking a subsequence, $\tau_n\rho_n(\gamma)\rightarrow
\bar\tau$. Since up to passing to a subsequence we also have
$\tau_n\rightarrow \tau_\infty$, we can deduce that
$\rho_n(\gamma)\rightarrow \tau_\infty^{-1}\circ \bar\tau$.

To prove that $\rho_n$ is converging, it is sufficient to check that
two converging subsequences of $\rho_n$ share the same limit.  On
the other hand, suppose that $\rho_\infty$ is a limit of a subsequence of
$\rho_n$, then $\rho_\infty$ makes $\phi_\infty$ equivariant:
\[
   \phi_\infty(\gamma \tilde{p})=\rho_\infty(\gamma)\phi_\infty(\tilde{p})~.
\]
This relation uniquely determines the action of $\rho_\infty(\gamma)$ on
$\tilde{F}_\infty$.  Since two isometries of $\AdS^3$ that
coincide on a spacelike surface are equal the result follows.

The fact that the left and right representations $\pi_1(S)\rightarrow
\PSL_2(\mathbb R)$ corresponding to $\rho_\infty$ are faithful and
discrete is a consequence of the fact that they are limit of faithful
and discrete representations.
\end{proof}

If $(\phi_{n_i})$ is a convergent subsequence of $(\phi_n)$, then Lemma
\ref{limitsptm:lm} implies that $N_{n_i}$ is a convergent
sequence of spacetime. Let $N_\infty$ be the limit of such
spacetimes. Its holonomy is by definition the limit of the holonomies
of the $N_n$.  In particular, we can concretely realize $N_n$
as an AdS metric $g_{N_n}$ on $S\times\mathbb R$, in such a way that
$g_{N_{n_i}}$ converges to an AdS metric $g_{N_\infty}$
as tensors on $S\times\mathbb R$ and
$(S\times\mathbb R, g_{N_\infty})\cong N_\infty$.

\begin{prop}\label{isom:prop}
$F_{n_i}$ converges to the lower boundary
$\partial_- C(N_\infty)$ of the convex
  core of $N_\infty$.  Moreover, the induced map
\[
\bar\phi_\infty:(S,h)\rightarrow \partial_- C(N_\infty)
\]
is an isometry.
\end{prop}

The proof of this proposition is based on the following lemma.

\begin{lemma}\label{unif:lm}
If $N$ is a MGHC anti de Sitter spacetime
and $N(k)$ is a Cauchy surface of constant
curvature $k\leq -1$, then the Lorentzian distance of any point
of $N(k)$ from the convex core of $N$ is smaller than
$\arctan\sqrt{|1+k|}$.
\end{lemma}

\begin{proof}
We consider the point $x_0$ on $N(k)$ with the biggest distance from the
convex core.  If $d$ is the distance between $x_0$ and the convex
core, then it is well known that $d<\pi/2$ and there exists a timelike
geodesic segment $\varsigma$ joining the point $x_0$ to a point $y_0$ on the
boundary of the convex core with length equal to $d$ \cite{BeBo}.

\begin{center}
\begin{figurehere}
\psfrag{Xi}{$\Xi$}
\psfrag{Nk}{$\tilde{N}(k)$}
\psfrag{N}{$\tilde{N}$}
\psfrag{Xid}{$\Xi^d$}
\psfrag{CN}{$C(\tilde{N})$}
\psfrag{pCN}{$\partial C(\tilde{N})$}
\psfrag{x0}{$\tilde{x}_0$}
\psfrag{y0}{$\tilde{y}_0$}
\psfrag{si}{$\tilde{\varsigma}$}
\includegraphics[width=0.5\textwidth]{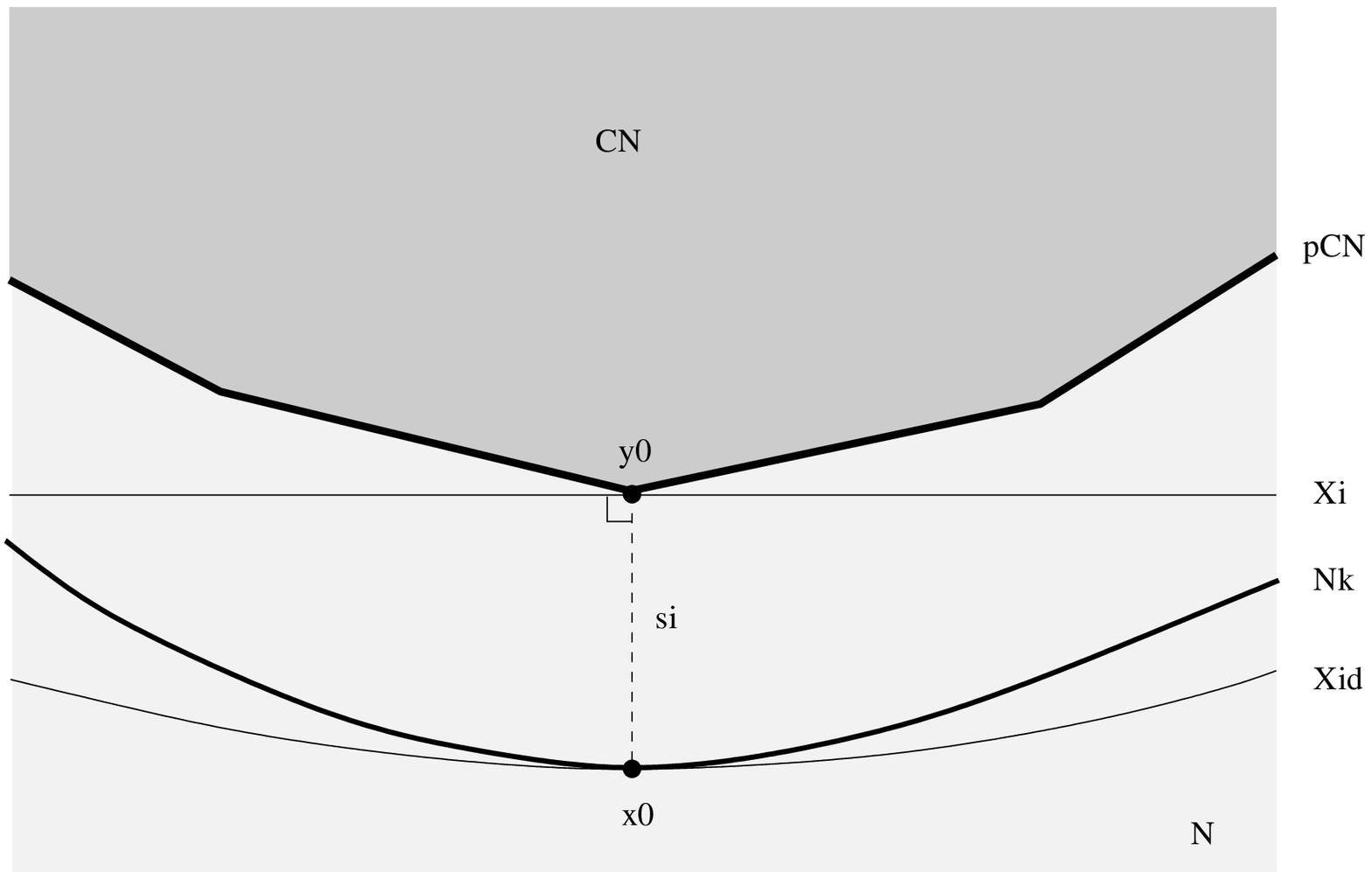}
\caption{{\small Estimating the distance between $\tilde{N}(k)$
and $C(\tilde{N})$.}}\label{close:fig}
\end{figurehere}
\end{center}

We consider now a lift $\tilde\varsigma$ of $\varsigma$
to the universal cover
$\tilde N\subset \AdS^3$ (see Figure \ref{close:fig}).
The plane $\Xi$ through $\tilde y_0$
orthogonal to $\tilde\varsigma$ is a support plane for the lifting of the
convex core.  Let $\Xi^d$ be the surface of points in $\AdS^3$ whose
distance from $\Xi$ is $d$.  It is a convex surface of constant
curvature $-1/\cos^2(d)$.  Clearly, $\tilde x_0$ lies on $\Xi^d$ and
$\tilde{N}(k)$ is contained in the convex side bounded by $\Xi^d$. In
particular, $\tilde{N}(k)$ and $\Xi^d$
are tangent at $\tilde x_0$ and, by
the maximum principle, the curvature of $\tilde{N}(k)$ is less than the
curvature of $\Xi^d$ at $\tilde x_0$.

We deduce that $k\leq -1/\cos^2(d)$, which implies that $\tan(d)\leq
\sqrt{|1+k|}$.
\end{proof}

\begin{proof}[Proof of Proposition \ref{isom:prop}]
Since the metrics $g_{N_{n_i}}$ converge to $g_{N_\infty}$,
the convex cores of the $N_{n_i}$
converge to the convex core of $N_\infty$. (Since the metrics
converge, the holonomy representations converge and so that their limit 
set in $\partial_\infty \AdS^3$ converge; therefore also their convex hulls, 
so the convex cores converge.)

In particular, the lower boundary
$\partial_- C(N_{n_i})$ of the convex core of $N_{n_i}$
converges to $\partial_- C(N_\infty)$.
By Lemma \ref{unif:lm}, the distance of any
point of $F_{n_i}$ from $\partial_- C(N_{n_i})$
is smaller that $\theta_{n_i}/2$. This
implies that $F_{n_i}$ converges to $\partial_- C(N_\infty)$.

In order to prove that the map $\bar\phi_\infty: (S,h)\rightarrow
\partial_- C(N_\infty)$
is an isometry, it is sufficient to show that $\bar\phi_\infty$
increases the distances.  Indeed, both $(S,h)$ and
$\partial_- C(N_\infty)$ are hyperbolic
surfaces and $\bar\phi_\infty$ is an homotopy equivalence.

We will prove that the lifting $\phi_\infty:(\tilde
S,h)\rightarrow \partial_- C(\tilde{N}_\infty)$
increases the lengths.  Given $\tilde{p},\tilde{q}\in\tilde{S}$,
we consider any path $\tilde\varsigma:[0,1]\rightarrow\tilde{S}$ connecting
$\tilde{p}$ and $\tilde{q}$ such that
\begin{itemize}
\item $\partial_- C(\tilde{N}_\infty)$ is smooth at almost all points
of $\tilde\varsigma_\infty:=\phi_\infty\circ \tilde\varsigma$;
\item $\ell(\tilde\varsigma_\infty)\leq d_\infty(\tilde{x}_\infty,
\tilde{y}_\infty)+\epsilon$.
\end{itemize}
where $\tilde{x}_\infty=\phi_\infty(\tilde{p})$ and $\tilde{y}_\infty=
\phi_\infty(\tilde{q})$ and
$d_\infty$ is the distance in $\partial_- C(\tilde{N}_\infty)$.

In the model $ \Hyp^2\times S^1$ of $\AdS^3$, the surfaces
$\tilde F_n$ are graphs of functions $e^{is_n}: \Hyp^2\rightarrow
S^1$ converging to $e^{is_\infty}: \Hyp^2\rightarrow S^1$ such that
$\partial_- C(\tilde{N}_\infty)$ is the graph of $e^{is_\infty}$.

We have $\tilde\varsigma_\infty(t)=(\xi(t), e^{is_\infty(\xi(t))} )$
with $\xi:[0,1]\rightarrow \Hyp^2$
Lipschitz function.  Take $\tilde\varsigma_n(t)=(\xi(t), e^{is_n(\xi(t))})$. 
For any smooth
point $\tilde{x}=(\xi, e^{is_\infty(\xi)})$ of
$\partial_- C(\tilde{N}_\infty)$ we have
$\grad(s_n)(\xi)\rightarrow\grad(s_\infty)(\xi)$.
Indeed, by convexity, tangent
planes of $\tilde F_n$ converge to support planes
of $\partial_- C(\tilde{N}_\infty)$.
 
By Lebesgue Theorem we have
\[
  \ell(\tilde\varsigma_n)=
\int_0^1\sqrt{\| \dot \xi\| ^2-\chi(\xi)\langle \dot \xi,
    \grad(s_n)\rangle^2}dt\rightarrow
\int_0^1\sqrt{\| \dot \xi\| ^2-\chi(\xi)\langle
    \dot \xi,\grad(s_\infty)\rangle^2}dt=\ell(\tilde\varsigma_\infty)
\]
since $\sqrt{\| \dot \xi\| ^2-\chi(\xi)\langle \dot \xi,
\grad(s_n)\rangle^2 }$
are all dominated by $\| \dot \xi\| $, which is an integrable function.

Since the endpoints of $\tilde\varsigma_n$
correspond to points $\tilde{x}_n=\varphi_n(\tilde{p}_n)$
and $\tilde{y}_n=\varphi_n(\tilde{q}_n)$
with $\tilde{p}_n\rightarrow\tilde{p}$ and $\tilde{q}_n\rightarrow\tilde{q}$,
we deduce that $d_{\tilde{h}}(\tilde{p},\tilde{q})\leq
d_\infty(\tilde{x}_\infty, \tilde{y}_\infty)+\epsilon$.  Since $\epsilon$ can be
chosen arbitrarly small,
\begin{equation}\label{inc:eq}
d_{\tilde{h}}(\tilde{p},\tilde{q})\leq d_\infty(\tilde{x}_\infty, \tilde{y}_\infty)\,.
\end{equation}
\end{proof}

So far, we have shown that $N_n$ is contained in a compact
subset of the space of MGH AdS structures, and any convergent subsequence
of $\bar\phi_n:S\rightarrow N_n$ converges to an isometric
embedding $\bar\phi_\infty:(S,h)\rightarrow N_\infty$, whose
image is the lower boundary $\partial_- C(N_\infty)$
of the convex core of $N_\infty$.

Let $\lambda_\infty$ be the bending lamination of
$\partial_- C(N_\infty)$. We will prove that
$\lambda_\infty=\lambda/2$.  Since the length spectrum of the third fundamental
form $\III_{F_n}$ converges to the intersection spectrum of
$\lambda/2$, it is sufficient to prove that it converges also to the
intersection spectrum of the bending lamination of $\partial_- C(N_\infty)$.

Now, let $\tilde F_n^\dual$ be the surface dual to $\tilde{F}_n$. Points
of $\tilde F_n^\dual$ are dual to tangent planes of $\tilde F_n$ and
$\tilde F_n^\dual$ is a past-convex surface
of constant curvature $-1/\sin^2(\theta_n/2)$
as (\ref{dlemb:eq}) shows.

Clearly, $\tilde F_n^\dual$ is invariant under the
holonomy action of $\pi_1(S)$, so it is contained in
$\tilde{N}_n^\dual$ and
its quotient is a Cauchy surface $F_n^\dual$ of $\tilde N_n^\dual$.
By (\ref{dlemb:eq}), the length spectrum of $\tilde F_n^\dual$ is equal to
the length spectrum of the third fundamental form
$\III_{F_n}$.

The boundary of the domain $\tilde{N}_\infty$ in $\AdS^3$ is the union of
two disjoint achronal surfaces: the past and the future singularities
of $\tilde{N}_\infty$, that are clearly invariant under the action
of $\pi_1(S)$.
 
By Proposition \ref{mehdi:prop}, the length spectrum of $F_n^\dual$ converges to
the length spectrum of the action of $\pi_1(S)$ on the future singularity
of $\tilde{N}_\infty^\dual$
(notice indeed that since $F_n^\dual$ is past-convex, in order to apply 
Proposition \ref{mehdi:prop} we need to exchange the time orientation).
On the other hand, by \cite{BeBo} the  length spectrum
of the future singularity of $\tilde {N}_\infty^\dual$
coincides with the intersection spectrum of the bending lamination
of the lower boundary of the convex core of $N_\infty$.

Combining these facts, we deduce that
\[
   \ell_{\III_{F_n}}(\gamma)\rightarrow\iota(\lambda_\infty,\gamma)
\]
so $\lambda_\infty=\lambda/2$.

\subsection{Asymptotic behavior of the measures $\tr(b_n)\omega_{h_n}$}

Let $(h_n,\,h_n^\dual)_{n\in\N}$ be a sequence of normalized hyperbolic
metrics on $S$ such that $h_n$ converges to $h$ and $h_n^\dual$ converges
to $[\lambda]$ in Thurston boundary of Teichm\"uller space, and
denote by $b_n$ the operator associated to $(h_n,\, h_n^\dual)$
provided by Corollary \ref{cr:b}.

Moreover, let $(\theta_n)_{n\in\N}$ be a sequence such that
$\theta_n \ell_{h_n^\dual}$ converges to $\iota(\lambda,\cbull)$
in the sense of spectra of closed curves.

In this section we study the asymptotic behavior of $\tr(b_n)$:
roughly speaking, it concentrates around the $h$-geodesic
representative of $\lambda$. Hence, we will always refer to $\lambda$ 
as to such an $h$-geodesic representative.

These results will turn useful in the proof of Theorem \ref{conv:thm}
and in Section \ref{ssc:limit_cn}.

\begin{prop}\label{tr:prop}
Let $V\subset S$ a closed subsurface such that $\partial V$ is
smooth and does not intersect $\lambda$.
Call $\lambda_V$ the $h$-geodesic sublamination $\lambda\cap V$.
Then
\[
\theta_n \int_V \tr(b_n)\omega_{h_n}\rightarrow \ell_h(\lambda_V)\, .
\]
where $\omega_g$ is the area form associated to $g$.
\end{prop}

In order to prove Proposition \ref{tr:prop}, we need the following lemma
that is analogous to Lemma \ref{unif:lm}.

\begin{lemma}\label{dist:lm}
Let $M$ be a hyperbolic end associated to some projective structure on $S$,
and let $M(k)$ be the surface of constant curvature $k$ with $k\in[-1,0)$.
Then the distance of any point of $M(k)$ from the boundary of $M$ is 
at most $\mathrm{arctanh}\sqrt{1+k}$.
\end{lemma}

The proof of Lemma \ref{dist:lm} is essentialy the same as in Lemma
\ref{unif:lm}. We leave the details to the reader.

\begin{cor}\label{bend:lm}
The family of isometric immersions $\sigma_n:(\tilde S,
\ch^2(\theta_n/2)\tilde h)\rightarrow \Hyp^3$ converges to a bending
map $\sigma_\infty:(\tilde S,h)\rightarrow \Hyp^3$, with bending
lamination $\lambda/2$.
\end{cor}
\begin{proof}
Since $\sigma_n$ are uniformly Lipschitz as maps $(\tilde{S},\tilde h)\rightarrow \Hyp^3$, 
they converge
up to subsequences to a locally convex surface. Combining Proposition \ref{im:prop}
and Lemma \ref{dist:lm}, we deduce that
this surface is the bent surface corresponding to $Gr_{\lambda/2}(S)$.
\end{proof}

\begin{proof}[Proof of Proposition \ref{tr:prop}]
We consider the embedding $\bar\sigma_n:S\rightarrow S_n\in M_n$
inside a hyperbolic end $M_n$ 
such that
\[
   I_{S_n}=\ch^2(\theta_n/2)h_n\qquad
\qquad B_{S_n}=\tanh(\theta_n/2)b_n\,.
\]
The projective structure $G_n$ at the ideal boundary of $M_n$
converges to $G_\infty=Gr_{\lambda/2}(h)$ (Proposition \ref{im:prop}).
Call $M_\infty$ the hyperbolic end determined by the $\mathbb{CP}^1$-surface
$G_\infty$.

A simple computation shows that the area element of $S_n$ with respect to
$I^\#_{S_n}$ is $\omega_{I_{S_n}^\#}=(\cosh^2(\theta_n/2)+\sinh^2(\theta_n/2)+
\sinh(\theta_n/2)\cosh(\theta_n/2)\tr(b_n))\omega_{h_n}$ and so the area of
$\bar\sigma_n(V)$ is
\begin{equation}\label{area:eq}
Area_{I^\#_{S_n}}(V)=
\left(
Area_{h_n}(V)+\frac{\theta_n}{2}\int_V \tr(b_n)
\omega_{h_n}
\right)
(1+o(\theta_n))
\end{equation}

As before,
we can identify $M_n\cup G_n\cong (S\times[0,\infty),g_{M_n})\cup S\times\{\infty\}$
in such a way that:
\begin{itemize}
\item
the developing map $dev_n:\tilde{S}\times[0,\infty]\rightarrow
\overline{\Hyp}^3$ converges to $dev_\infty$ (and so
$g_{M_n}\rightarrow g_{M_\infty}$);
\item
$\bar\sigma_n$
converges to the pleated surface
$\bar\sigma_\infty:S\rightarrow (S\times\{0\},g_{M_\infty})$
\end{itemize}

Call $\partial^d M_n$ the surface in $M_n$ at distance $d$
from the boundary and let
$\Pi_{\partial^d M_n}:G_n\rightarrow \partial^d M_n$ the projection
introduced in Lemma \ref{boundgr:lm}.

There exists two numbers $\epsilon_n<\delta_n$ such that
$S_n$ is contained between
$\partial^{\epsilon_n}M_n$ and
$\partial^{\delta_n}M_n$ and, by Lemma \ref{dist:lm}, 
$\delta_n\rightarrow 0$ as $n\rightarrow+\infty$.

By the monotonicity result proved in \cite{horo},
\[
\Pi_{\partial^{\epsilon_n}M_n}^*(\, I^\#_{\partial^{\epsilon_n}M_n}\,)
\leq
\Pi_{S_n}^* (\, I^\#_{S_n}\,)
\leq
\Pi_{\partial^{\delta_n} M_n}^* (\,I^\#_{\partial^{\delta_n}M_n}\,)\, .
\]
If $\lambda_n$ is the bending lamination of $M_n$,
the grafted metric on
$\partial^d M_n$ makes it isometric to
$e^{2d}g_{G_n}$, where $g_{G_n}$ is Thurston metric on
the projective surface
$G_n=Gr_{\lambda_n}(\partial M_n,g_{M_n})$.

So we deduce that
\[
     e^{2\epsilon_n}\omega_{G_n}\leq
\Pi_{S_n}^* (\,\omega_{I^\#_{S_n}}\,) \leq
e^{2\delta_n}\omega_{G_n}
\]
and so
\begin{equation}\label{comp:eq}
e^{2\epsilon_n}Area_{G_n}(\Pi_{S_n}^{-1}(\bar\sigma_n(V)))\leq
Area_{I^\#_{S_n}}(V) \leq
e^{2\delta_n}Area_{G_n}(\Pi_{S_n}^{-1}(\bar\sigma_n(V)))\, .
\end{equation}

Since $G_n\rightarrow G_\infty$, their Thurston metrics converge
to $g_{G_\infty}$. We claim that $\Pi_{S_n}^{-1}(\bar\sigma_n(V))$ converge
to $\Pi_{S_\infty}^{-1}(\bar\sigma_\infty(V))$ in the Hausdorff sense,
and so
\begin{equation}\label{area2:eq}
Area_{I^\#_{S_n}}(V)\rightarrow
Area_{G_\infty}(\Pi_{S_\infty}^{-1}(\bar\sigma_\infty(V)))=Area_h(V)+\frac{1}{2}\ell_h(\lambda_V)
\end{equation}
by Equation (\ref{comp:eq}).

The result will follow by
comparing Equations (\ref{area:eq}) and (\ref{area2:eq}).

In order to prove the claim, it is enough to prove that
$\partial \Pi_{S_n}^{-1}(\bar\sigma_n(V))\rightarrow \partial
\Pi_{S_\infty}^{-1}(\bar\sigma_\infty(V))$, which would follow
from the fact that $\Pi_{S_n}^{-1}\circ\bar\sigma_n|_{\partial V}$
converges to $\Pi_{S_\infty}^{-1}\circ\bar\sigma_\infty|_{\partial V}$.

Notice that $\Pi_{S_n}$ is a diffeomorphism, and
$\Pi_{S_\infty}^{-1}\circ\bar\sigma_\infty|_{\partial V}$ is well-defined
and continuous since $\partial V$ does not intersect $\lambda$.

Let $(p_n)_{n_\in\N}$ be a sequence of points in $\partial V$ such that
$p_n\rightarrow p$. The point $\Pi_{S_n}^{-1}(\bar\sigma_n(p_n))$
is the ideal point of the horocycle $U_n$
tangent to $S_n$ at $\bar\sigma_n(p_n)$. By convexity of the surfaces
$S_n$, one can easily see that $U_n$ converges to a horocycle
$U_\infty$ tangent to $S_\infty$ at $\bar\sigma_\infty(p)$.
Since $p\notin \lambda$, such a $U_\infty$ is unique and so
its ideal point is necessarily $\Pi^{-1}_{S_\infty}(p)$.
\end{proof}

\begin{cor}\label{asym:cor}
Let $h_n$ be a sequence of hyperbolic metrics converging to $h$, and let
$h_n^\dual$ be a diverging sequence of metrics.
Then, $(\theta_n\ell_{h_n^\dual}(\gamma))_{n\in\N}$
is bounded for every $\gamma\in\pi_1(S)$
if and only if
$\theta_n\int_S \tr(b_n)\omega_{h_n}$ is bounded.
\end{cor}

\subsection{Proof of Theorems \ref{tm:limit} and \ref{conv:thm}}

\begin{proof}[Proof of Theorem \ref{tm:limit}]
Let us fix a hyperbolic metric $h$ on $S$ and a sequence
of hyperbolic metrics $h_n^\dual$ converging to a point
$[\lambda]$ in Thurston boundary of $\mathcal T$.
Fixing a sequence $\theta_n$ such that $\theta_n\ell_{h_n^\dual}$ converges
to $\iota(\lambda,\cbull)$, let us set
\[
L_{e^{i\theta_n}}(h, h_n^\dual)=(h^1_n, h^2_n)\,.
\]

In subsection \ref{re:ssec} we have shown that $h^1_n\rightarrow E_{\lambda/2}(h)$.
To conclude the proof  we need to prove that
\begin{equation}\label{antip:eq}
\theta_n\ell_{h^2_n}\rightarrow\iota(\lambda,\cbull)\,,
\end{equation}
The main issue to prove (\ref{antip:eq}) is to show that
for every $\gamma\in\pi_1(S)$, there is $C=C(\gamma)$ such that
\begin{equation}\label{bant:eq}
\theta_n\ell_{h^2_n}(\gamma)<C\,.
\end{equation}

Corollary \ref{asym:cor} indicates that, 
in order to prove (\ref{bant:eq}), it is sufficient to bound
\[
   \theta_n\int_S  \tr (b'_n)\omega_{h^1_n}
\]
where $b'_n$ is the $h^1_n$-self-adjoint operator such that
$h^2_n=h^1_n(b_n'\cbull,b_n'\cbull)$.
Now we have that $b_n'=\beta'_n b_n\beta_n^{-1}$ where
$\beta_n=\cos(\theta_n/2)E+\sin(\theta_n/2)Jb_n$ and
$\omega_{h^1_n}=\det(\beta_n)\omega_h=\omega_h$.  So, $\theta_n\int_S\tr
(b_n')\omega_{h^1_n}=\theta_n\int_S\tr(b_n)\omega_h$ that in turn is bounded
since $\theta_n\ell_{h_n^\dual}(\gamma)$ is bounded for every $\gamma\in\pi_1(S)$
by hypothesis.

It follows that there exists a measured geodesic lamination $\mu$ such
that, up to passing to a subsequence, $\theta_n\ell_{h_n^2}\rightarrow\iota(\mu,\cbull)$.  To
show that $\mu=\lambda$, notice that by Proposition \ref{re:prop} we
have that
\begin{equation}\label{com:eq}
      L^1_{e^{i\theta_n}}(h^1_n,h^2_n)\rightarrow E_{\mu/2}(h_\infty)
\end{equation}
On the other hand, we have that 
$$L^1_{e^{i\theta_n}}(h_n^1,h^2_n)=L^1_{e^{i\theta_n}}(L_{e^{i\theta_n}}(h,h_n^\dual))=
L^1_{e^{2i\theta_n}}(h,h_n^\dual)~. $$ 
So, applying again Proposition \ref{re:prop}, we obtain that
\begin{equation}\label{com2:eq}
    L^1_{e^{i\theta_n}}(h_n^1,h_n^2)\rightarrow
E_{\lambda}(h)=E_{\lambda/2}(h_\infty)\,.
\end{equation}

Comparing (\ref{com:eq}) and (\ref{com2:eq}) we conclude that
$E_{\lambda/2}(h_\infty)=E_{\mu/2}(h_\infty)$ and so $\lambda=\mu$.
\end{proof}

\begin{proof}[Proof of Theorem \ref{conv:thm}]
Let $h$ and $h_n$ as above.
For all $z=t+is\in \overline{\HH}$, we have to prove
\begin{equation}\label{lll:eq}
P'_{\theta_n z}(h,h_n^\dual)\rightarrow Gr_{s\lambda/2}(E_{-t\lambda/2}(h))\,.
\end{equation}
 
Recall that
$P'_{\theta_n z}(h,h_n^\dual)=SGr_{\theta_n s}\circ L'_{-\theta_n t}(h,h_n^\dual)$.
Note that if we put
\[
       L'_{-\theta_n t}(h,h_n^\dual)=(h_n^1, h_n^2)
\]
Theorem \ref{tm:limit} shows that
$h_n^1\rightarrow h_\infty=E_{-t\lambda/2}(h)$ and $\theta_n\ell_{h_n^2}\rightarrow\lambda$.

Applying Proposition \ref{im:prop} we conclude that
\[
P'_{\theta_nz}(h, h_n^\dual)=SGr_{\theta_n s}(h_n^1,h_n^2)\rightarrow 
Gr_{s\lambda/2}(E_{-t\lambda/2}(h))\,.
\]
\end{proof}

\subsection{Convergence of the distances}

The aim of this section is to  study the asymptotic behavior of the sequence
of distances induced by the metrics $\theta_n^2 h_n^\dual$
introduced in the previous section.
By our assumption, the length spectrum of $h_n^\dual$
rescaled by $\theta_n$ converges to
the intersection with $\lambda$. Notice that
this assumption only concerns the isotopy
class of $h_n^\dual$.
On the other hand, once we concretely fix $h$,
the metric $h_n^\dual$ is uniquely
determined in its 
isotopy class by requiring the the identity map $(S,h)\rightarrow(S, h_n^\dual)$
is minimal Lagrangian.
The result we consider in this section deals with the asymptotic
behavior of $h_n^\dual$
considered as concrete metrics on $S$. Clearly, these results are valid
for this choice of gauge, and are no longer valid if
we change $h_n^\dual$ by some isotopy.

Notice that the
representative $\lambda$ of a point in Thurston boundary of $\mathcal T(S)$
can be chosen to be
a measured geodesic lamination for any hyperbolic metric on $S$.
In order to study the behavior of $h_n^\dual$, it is natural
to fix $\lambda$ as the concrete
measured geodesic realization with respect to the metric $h$.

%
Let us fix a universal cover $\tilde S\rightarrow S$ and
denote by $\tilde h$ and $\tilde h^\dual_n$ the pull-back of the metrics $h$
and $h^\dual_n$ on $\tilde S$.  Finally let $\tilde\lambda$  be the
pull-back of $\lambda$ on $\tilde S$.

\begin{prop}\label{limitest:pr}
For every $\tilde{p}, \tilde{q}\in \tilde
S\setminus\tilde\lambda$ we have
\begin{equation}\label{limitest:eq}
    \theta_n d_{\tilde{h}^\dual_n}(\tilde{p},\tilde{q})
\rightarrow \iota( \tilde{\alpha},\tilde\lambda)
\end{equation}
where $\tilde{\alpha}$ is any smooth path in $\tilde S$ joining $\tilde{p}$
to $\tilde{q}$ and
meeting each leaf of $\tilde\lambda$ at most once and
transversely.  Moreover, the
convergence is uniform on compact subsets of
$\tilde S\setminus\tilde\lambda$.
\end{prop}

Theorem \ref{tm:boundary} is a direct consequence of this statement, so that its proof
will be a consequence of the proof of this proposition.

First, notice that it is sufficient to
prove Proposition \ref{limitest:pr} after rescaling $\theta_n$ and
$\lambda$ by some arbitrary factor. In particular, we may assume that the 
projective surface $Gr_{\lambda/2}(S,h)$ is quasi-Fuchsian.
This technical assumption
will simplify some steps of the proof.

First we show that

\begin{lemma}\label{first:lm}
\begin{equation}\label{first:eq}
   \liminf \theta_n d_{\tilde{h}_n^\dual}(\tilde{p},\tilde{q})\geq
\iota( \tilde{\alpha},\tilde\lambda)\,.
\end{equation}
\end{lemma}
\begin{proof}
Let $M_n$ be the hyperbolic end introduced in Section \ref{ssc:imaginary}
with $h_n=h$
and let $\sigma_n:\tilde S\rightarrow\tilde M_n$ be the lifting of the
embedding $\bar\sigma_n:S\rightarrow M_n$.
By Proposition \ref{im:prop}, $M_n$ converges to the hyperbolic end
facing $Gr_{\lambda/2}(S,h)$. In particular, we may assume that $M_n$ are all
quasi-Fuchsian, so that $\tilde M_n$ is a concave region of $ \Hyp^3$.

By Corollary \ref{bend:lm}, the family of embeddings $(\sigma_n)_{n\in\N}$
converges
to the bending map $\sigma_\infty:(\tilde S,\tilde{h})\rightarrow\Hyp^3$,
with bending
lamination $\tilde\lambda/2$. 

Let $\tilde{r}_n^\dual: \tilde{S}^\dual_n\rightarrow \partial \tilde{M}^\dual_n$ be
the $1$-Lipschitz map defined in Lemma
\ref{cptlam:lm}. 
We denote by $\Xi_n(\tilde p)$ the plane in $\Hyp^3$ corresponding to 
$\tilde r_n^\dual(\sigma_n^\dual(\tilde p))$ and by $\Xi_n(\tilde q)$ the plane
corresponding to $\tilde r_n^\dual(\sigma_n^\dual(\tilde q))$.

We have that $\Xi_n(\tilde p)$ and $\Xi_n(\tilde q)$ are both support
planes of $\partial \tilde M_n$.  Let us put $
\iota_n=d_{\partial\tilde M_n^\dual}(\tilde
r_n^\dual(\sigma_n^\dual(\tilde p)), \tilde
r_n^\dual(\sigma_n^\dual(\tilde q)))\,.$
                                                        
 Since $\tilde r_n$ decreases the distances we deduce that 
\begin{equation}\label{bb:eq}
 \iota_n\leq d_{\III_n}(\sigma_n(\tilde p),\sigma_n(\tilde q))\sim \frac{\theta_n}{2}
 d_{\tilde h_n^\dual}(\tilde p,\tilde q)\,.
 \end{equation}
 We claim that
$\Xi_n(\tilde p)$ and $\Xi_n(\tilde q)$ converge to the support planes
 of $\partial \tilde{M}_\infty$ at
$\sigma_\infty(\tilde p)$ and $\sigma_\infty(\tilde q)$ respectively.
 
 Let us explain how the conclusion follows.
 The claim ensures that 
 it is possible to construct a sequence of arcs $\tilde\alpha_n:[0,1]\rightarrow\partial \tilde{M}_n$ such that
 \begin{itemize}
 \item each $\tilde \alpha_n$ intersects each leaf of $\tilde\lambda_n$ at
   most once and transversely; moreover end-points of $\tilde\alpha_n$
   are contained in $\partial \tilde{M}_n\setminus\tilde\lambda_n$;
 \item end-points of $\tilde\alpha_n$ are close to $\Xi_n(\tilde p)$
   and $\Xi_n(\tilde q)$ in the sense that $d_{\partial
     \tilde{M}_n^\dual}(\tilde\alpha^\dual_n(0),r_n^\dual(\sigma_n^\dual(\tilde
   p)))$ and $d_{\partial \tilde{M}_n^\dual}(\tilde\alpha^\dual_n(1),
   r_n^\dual(\sigma_n^\dual(\tilde q)))$ converge to $0$.
 \item $\tilde\alpha_n$ converges to a path $\tilde\alpha_\infty$
   connecting the stratum of $\partial \tilde{M}_\infty\setminus\tilde\lambda$
   containing $\sigma_\infty(\tilde p)$ with the stratum containing
   $\sigma_\infty(\tilde q)$
  \end{itemize}
 Thus we have
$ \iota(\tilde\lambda/2,\tilde\alpha)=\iota(\tilde\lambda/2,\tilde\alpha_\infty)=
 \lim_n\iota(\tilde\lambda_n,\tilde\alpha_n)$.
 By (\ref{initsing:eq})  $\iota(\tilde\lambda_n,\tilde\alpha_n)=
 d_{\partial \tilde{M}_n}
(\tilde \alpha_n^\dual(0), \tilde\alpha_n^\dual(1))$, so we conclude that
 \[
 \iota(\tilde\lambda/2,\tilde\alpha)=\lim_n\iota_n
 \]
 and the conclusion of the Lemma follows from (\ref{bb:eq}).

In order to prove the claim
recall that $\tilde r_n^\dual$ sends any point $\tilde x^\dual$ of $\tilde S_n^\dual$ to a point in 
$\partial \tilde{M}_n^\dual$ that lies in the past of $\tilde x^\dual$ (Lemma
\ref{cptlam:lm}).

\begin{center}
\begin{figurehere}
\psfrag{Xin}{$\Xi_n(\tilde{p})$}
\psfrag{TSn}{$T_{\sigma_n(\tilde{p})}\tilde{S}_n$}
\psfrag{sin}{$\sigma_n(\tilde{p})$}
\psfrag{Mn}{$\tilde{M}_n$}
\psfrag{Sn}{$\tilde{S}_n$}
\psfrag{pMn}{$\partial \tilde{M}_n$}
\includegraphics[width=0.5\textwidth]{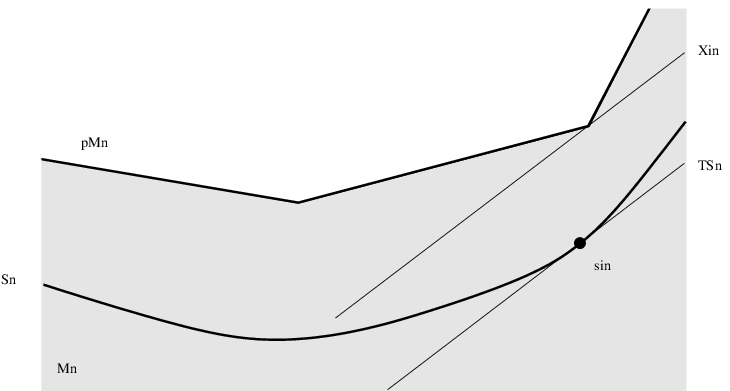}
\caption{{\small Convergence of the support planes.}}\label{xi:fig}
\end{figurehere}
\end{center}

This implies that $\Xi_n(\tilde p)$ and $T_{\sigma_n(\tilde p)}\tilde{S}_n$
are disjoint planes.
In particular we deduce that $\Xi_n(\tilde p)$ separates $\sigma_n(\tilde p)$ from
$\partial \tilde{M}_n$ (see Figure \ref{xi:fig}).
This easily implies that $\Xi_n(\tilde p)$ converges to the support plane
of $ \partial \tilde{M}_\infty$ at $\sigma_\infty(\tilde p)$
(that is unique by our assumption that $\tilde p$ does
not lie on $\tilde \lambda$).
 Analogously $\Xi_n(\tilde q)$ converges to the support plane of
$\partial \tilde{M}_\infty$
at $\sigma_\infty(\tilde q)$.
\end{proof}

To conclude the proof of Proposition \ref{limitest:pr} we need to show that
\begin{equation}\label{reverse:eq}
   \limsup \theta_n d_{\tilde{h}_n^\dual}(\tilde{p},\tilde{q})
\leq\iota( \tilde\alpha, \tilde\lambda)\,.
\end{equation}

In order to estimate $\theta_nd_{\tilde h_n^\dual}(\tilde{p},\tilde{q})$,
we need the following result.

\begin{lemma}\label{angles:lm}
Let $\tilde{U}$ be any convex smooth surface in $ \Hyp^3$ and let
$\tilde{x},\tilde{y}\in\tilde{U}$ such that
the support planes $\Xi_{\tilde{x}}$ and $\Xi_{\tilde{y}}$
at $\tilde{x}$ and $\tilde{y}$ intersect. Then the distance between
$\tilde{x}$ and $\tilde{y}$ with respect to the third fundamental
form of $\tilde{U}$ is less than the angle 
between $\Xi_{\tilde{x}}$ and $\Xi_{\tilde{y}}$.
\end{lemma}

We first prove a $2$-dimensional version of this lemma.
\begin{sublemma}\label{angles2d:lm}
Let $\tilde\varsigma$ be a convex curve in $ \Hyp^2$ joining
two points  $\tilde{x}, \tilde{y}$.
Suppose that the support lines $\tilde l_{\tilde{x}}$ and
$l_{\tilde{y}}$ at $\tilde{x}$ and $\tilde{y}$ intersect.
Then the angle they form is bigger than the integral of the 
curvature of $\tilde\varsigma$.
\end{sublemma}
\begin{proof}
By Gauss-Bonnet formula, the area bounded by $l_{\tilde{x}}$,
$l_{\tilde{y}}$ and $\tilde\varsigma$ is equal to
the difference between the wanted angle and the integral of the curvature.
\end{proof}

\begin{proof}[Proof of Lemma \ref{angles:lm}]
First we translate the condition that planes intersect in terms of a condition
of dual points $\tilde{x}^\dual$ and $\tilde{y}^\dual$.
Recalling that $\tilde{x}^\dual$ and $\tilde{y}^\dual$ are unit vector orthogonal to the planes (and pointing 
in the concave region bounded by $\tilde{U}$),
we easily derive that $\Xi_{\tilde{x}}$ and $\Xi_{\tilde{y}}$ intersect
if and only if the segment $\tilde{\kappa}^\dual$
joining $\tilde{x}^\dual$ to $\tilde{y}^\dual$ in $\dS^3$ is spacelike,
in which case the
angle between $\Xi_{\tilde{x}}$ and $\Xi_{\tilde{y}}$
coincides with the length of $\tilde{\kappa}^\dual$.

Let $\Pi$ be a timelike plane containing $\tilde{\kappa}^\dual$.
Since $\tilde{U}^\dual$ is an achronal surface,
$\tilde{U}^\dual\cap \Pi$ is a curve
containing $\tilde{x}^\dual$ and $\tilde{y}^\dual$.
Let $\tilde\varsigma^\dual$
be the segment on $\tilde{U}^\dual\cap \Pi$
connecting $\tilde{x}^\dual$ to $\tilde{y}^\dual$.
Clearly, the length of $\tilde\varsigma^\dual$
is greater than the distance on
$\tilde{U}^\dual$ between $\tilde{x}^\dual$ and $\tilde{y}^\dual$.
In order to conclude, it is sufficient to prove that the
length of $\tilde\varsigma^\dual$ is less than the length
of $\tilde{\kappa}^\dual$. 

Notice that this is a $2$-dimensional problem.
In fact, let $\Upsilon$ be the timelike linear $3$-space in $\mathbb R^{3,1}$
such that $\Pi=\Upsilon\cap \dS^3$.
The intersection $\Upsilon\cap  \Hyp^3$ is a hyperbolic plane denoted by
$\Lambda$.
Points on $\Pi$ correspond to planes of $ \Hyp^3$ that orthogonally meet
$\Lambda$. 
In particular, points on $\Pi$ bijectively correspond to lines on
$\Lambda$ and points of $\Lambda$ correspond to spacelike lines of $\Pi$.

Consider the curve $\tilde\varsigma$ on $\Lambda$
of points corresponding to support lines of $\tilde\varsigma^\dual$.
Notice that support lines at the endpoints of $\tilde\varsigma$
are $\Xi_{\tilde{x}}\cap \Lambda$ and $\Xi_{\tilde{y}}\cap\Lambda$, so the angles
these lines form is equal to the angle between
$\Xi_{\tilde{x}}$ and $\Xi_{\tilde{y}}$
and it is equal to the length of $\tilde{\kappa}^\dual$.
On the other hand the length of $\tilde\varsigma^\dual$
is equal to the integral of the curvature of $\tilde\varsigma$.
Thus the conclusion follows from Sublemma \ref{angles2d:lm} .
\end{proof}

Given any geodesic $\tilde\varsigma\in(\tilde S,\tilde h)$, 
we say that $\sigma_n(\tilde\varsigma)$ is a short path
if the support planes at $\sigma_n(\tilde\varsigma(0))$
and $\sigma_n(\tilde\varsigma(1))$ intersect. 
If $\sigma_n(\tilde\varsigma)$
is a short path, then we denote by $\eta_n(\tilde\varsigma)\in(0,\pi)$ 
the angle between the support planes at its endpoints.

Analogously, we say that $\sigma_\infty(\tilde\varsigma)$
is a short path if the endpoints of
$\sigma_\infty(\tilde\varsigma)$
are outside the bending lamination and the corresponding support planes
intersect. In this case, $\eta_\infty(\tilde\varsigma)$
is the angle between such planes.

Clearly, if $\sigma_\infty(\tilde\varsigma)$ is a short path,
then $\sigma_n(\tilde\varsigma)$ is definitively a short path and
$\eta_n(\tilde\varsigma)\rightarrow\eta_\infty(\tilde\varsigma)$.

\begin{lemma}\label{roof:lm}
There exists $\epsilon_0$ such that if $\tilde\alpha:[0,1]\rightarrow\tilde{S}$
is a geodesic path for $\tilde{h}$ of length
less than $\epsilon_0$, then
\begin{itemize}
\item $\sigma_\infty(\tilde\alpha)$ is a short path;
\item $\iota(\tilde\alpha,\tilde\lambda)\leq
\eta_\infty(\tilde\alpha)\leq 
(1+\ell_{\tilde{h}}(\tilde\alpha))\iota(\tilde\alpha,\tilde\lambda)$
\end{itemize}
\end{lemma}
\begin{proof}
The first point easily follows, since $\tilde{S}_\infty$
is invariant under the action of
a co-compact group of isometries of $ \Hyp^3$.

About the second point,
notice that the first inequality is given by Lemma \ref{angles:lm}.
The second inequality is more subtle.
Choosing $\epsilon_0$ sufficiently small,
we can suppose that either $\tilde\alpha$ intersects only one 
isolated leaf of $\tilde\lambda$ or it intersects no isolated leaf.
In the first case the second inequality is obvious.

Up to taking a smaller $\epsilon_0$, we may suppose that,
if $T$ is a hyperbolic triangle with an edge $e$ 
of length $l\leq\epsilon_0$ and the angles $\vartheta_1,\vartheta_2$
adjacent to $e$
less than $\pi/4$, then the area of $T$ is less than $l\vartheta_2$.
Now take a geodesic $\tilde\alpha$
on $\tilde S$ of length less than $\epsilon$ which
does not intersect the isolated leaves of the lamination. 
Taking any subdivision $\tilde{\alpha}(0)=\tilde{p}_1,\ldots,
\tilde{p}_{m+1}=\tilde\alpha(1)$
of $\tilde{\alpha}$, we consider the support planes
$\Xi_i$ of $\sigma_\infty(\tilde S)$ at $\sigma_\infty(\tilde{p}_i)$.
If the subdivision is sufficiently fine,
then the angles $\vartheta_i$ between $\Xi_i$ and $\Xi_{i+1}$ are
less than $\pi/4$.
Now consider the boundary $\partial K$ of the convex set $K$ obtained
by intersecting the half-spaces
bounded by $\Xi_1,\dots,\Xi_{m+1}$ and containing $\sigma_\infty(\tilde{S})$.
Notice that $\partial K$ is a finite bent surface:
indeed, its bending lines are $\Xi_i\cap \Xi_{i+1}$ for every $i$ such that
$\Xi_i$ and $\Xi_{i+1}$ are different.
We claim that $\eta_\infty(\tilde\alpha)\leq
(1+\ell_{\tilde{h}}(\tilde\alpha))\sum\vartheta_i$.
Taking a sequence of arbitrary fine subdivisions, we have that
$\sum\vartheta_i\rightarrow\iota(\tilde\alpha,\tilde\lambda)$, so
the second inequality follows from the claim.

In order to prove the claim we use an inductive argument.  Notice
that, if $\Xi_1,\Xi_2,\Xi_3$ are distinct, then
$\Xi_1\cap\Xi_2\cap\Xi_3=\emptyset$.  Thus, there is a plane $\Lambda$
orthogonal to all of them.  The triangle $T$ obtained by intersecting
$\Lambda\cap \Xi_1$, $\Lambda\cap \Xi_2$ and $\Lambda\cap \Xi_3$ has
angles $\vartheta_1,\vartheta_2$ and $\pi-\bar\vartheta_1$, where
$\bar\vartheta_1$ is the angle formed by $\Xi_1$ and $\Xi_3$.

\begin{center}
\begin{figurehere}
\psfrag{th1}{$\vartheta_1$}
\psfrag{th2}{$\vartheta_2$}
\psfrag{th3}{$\vartheta_3$}
\psfrag{bth2}{$\bar\vartheta_2$}
\psfrag{bth1}{$\bar\vartheta_1$}
\psfrag{p1}{$\tilde{p}_1$}
\psfrag{p2}{$\tilde{p}_2$}
\psfrag{p3}{$\tilde{p}_3$}
\psfrag{p4}{$\tilde{p}_4$}
\psfrag{p5}{$\tilde{p}_5$}
\psfrag{a}{$\tilde{\alpha}$}
\psfrag{pK}{$\partial K$}
\psfrag{Xi1}{$\Xi_1$}
\psfrag{Xi2}{$\Xi_2$}
\psfrag{Xi3}{$\Xi_3$}
\psfrag{Xi4}{$\Xi_4$}
\psfrag{Xi5}{$\Xi_5$}
\psfrag{T}{$T$}
\includegraphics[width=\textwidth]{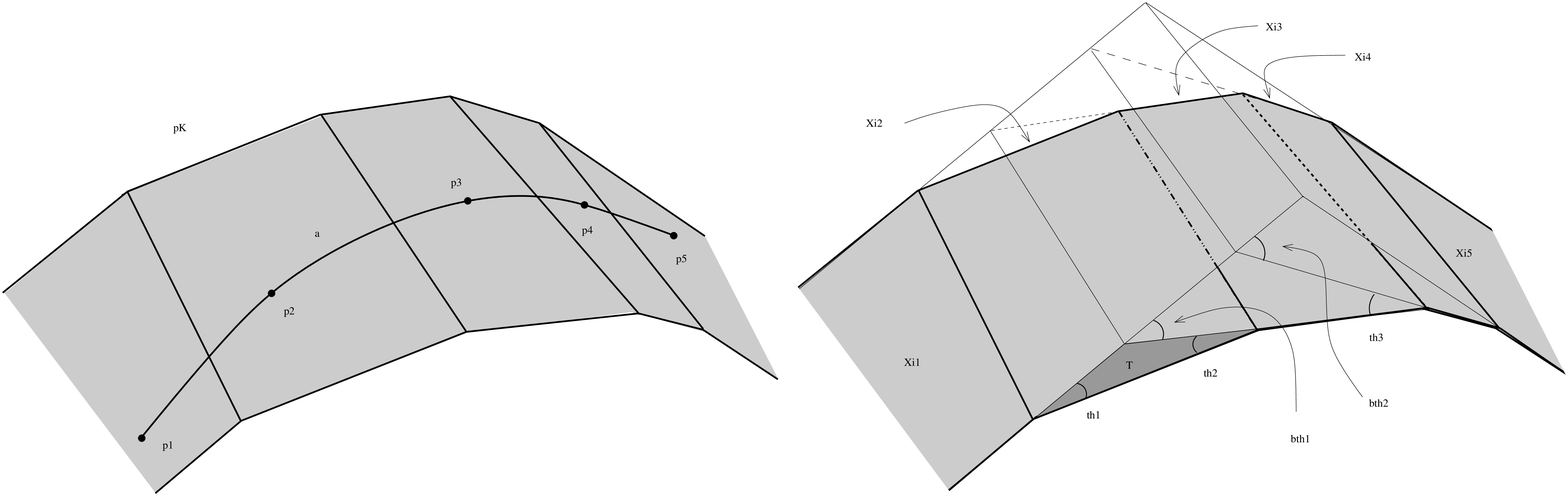}
\end{figurehere}
\end{center}

Notice
that the length of the edge between $\vartheta_1$ and $\vartheta_2$ is
less than the distance between $\tilde{p}_1$ and $\tilde{p}_2$, and so
it is smaller than $\epsilon_0$.  We conclude that the area of $T$ is
less than $d_{\tilde{h}}(\tilde{p}_1,\tilde{p}_2)\vartheta_2$.  In
particular, we deduce that $\bar\vartheta_1\leq
\vartheta_1+\vartheta_2 (1+d_{\tilde{h}}(\tilde{p}_1,\tilde{p}_2))$.

If at least two planes among $\Xi_1,\Xi_2,\Xi_3$ coincide, then the
area of $T$ is zero, and so $\bar{\vartheta}_1=\vartheta_1+\vartheta_2
\leq\vartheta_1+\vartheta_2(1+d_{\tilde{h}}(\tilde{p}_1,\tilde{p}_2))$.

Apply now the same argument to the planes $\Xi_1, \Xi_3, \Xi_4$.
If $\bar\vartheta_2$ is the angle formed by $\Xi_1$ and $\Xi_4$, then
$\bar\vartheta_2\leq \bar\vartheta_1+\vartheta_3
(1+d_{\tilde{h}}(\tilde{p}_1,\tilde{p}_3))\leq
\bar\vartheta_1+\vartheta_2(1+d_{\tilde{h}}(\tilde{p}_1,\tilde{p}_2))+
\vartheta_3(1+d_{\tilde{h}}(\tilde{p}_1,\tilde{p}_3))$.
Iterating this procedure, we deduce that the angle between
$\Xi_1$ and $\Xi_m$ is bounded by
$\vartheta_1+\vartheta_2(1+d_{\tilde{h}}(\tilde{p}_1,\tilde{p}_2)+
\vartheta_3(1+d_{\tilde{h}}(\tilde{p}_1,\tilde{p}_3))+\ldots+
\vartheta_m(1+d_{\tilde{h}}(\tilde{p}_1,\tilde{p}_m))$ and
this quantity is less than $(1+\ell_{\tilde{h}}(\tilde{\alpha}))
\sum\vartheta_i$.
\end{proof}

We can now prove (\ref{reverse:eq}). Fix $\epsilon<\epsilon_0$ and 
subdivide the geodesic $\tilde{\alpha}$ joining $\tilde{p}$ to
$\tilde{q}$ into segments $\tilde{\alpha}_i$ of length less than $\epsilon$.
Let $\tilde{p}=\tilde{p}_1,\tilde{p}_2,\ldots,\tilde{p}_{m+1}=\tilde{q}$
be endpoints of such subdivision.

For $n$ large, $\sigma_n(\tilde{\alpha}_i)$
are short paths and so by Lemma \ref{angles:lm}
\[
  \limsup d_{\III_{\tilde{S}_n}}(\tilde{p},\tilde{q})\leq 
\sum_{i=1}^{m}\eta(\tilde{\alpha}_i)
\]
On the other hand, applying Lemma \ref{roof:lm} we deduce that
\[
   \sum_{i=1}^m\eta(\tilde{\alpha}_i)\leq
\iota(\tilde\alpha,\tilde\lambda)(1+\epsilon)
\]
The uniform convergence follows from the fact that
the whole argument works as well,
if we consider sequences
of points $\tilde{p}_n\rightarrow \tilde{p}$
and $\tilde{q}_n\rightarrow \tilde{q}$
belonging to a compact subset of $\tilde{S}\setminus\tilde{\lambda}$.
This concludes the proof of Proposition \ref{limitest:pr}.

\section{Behavior of the centers}
\label{ssc:limit_cn}

In this section we want to discuss the behavior of the centers $c_n$ when
$h$ is fixed and $h_n^\dual\rightarrow[\lambda]$ in Thurston compactification.
Our aim is to prove that
the limit point(s) of $c_n$ does not only depend on $h$ and $[\lambda]$
but also on the sequence $h_n^\dual$. As a consequence, we will see that the
analog of Theorem \ref{tm:limit} does not hold if the sequence of centers $(c_n)$
converges to a projective measured lamination, see Corollary \ref{cor:center}.

Fix a hyperbolic metric $h$ on $S$ and let $c'$ be a point in the boundary
of the augmented Teichm\"uller space of $S$. Then $c'$ can be considered as
a complete hyperbolic metric of finite area on $S\setminus\Gamma$, where
$\Gamma$ is the disjoint union of simple closed curves $\gamma_1,\dots,
\gamma_l$. Up to isotopy, we can assume that $\gamma_i$ is a geodesic for $h$
and we denote by $\ell_i$ the $h$-length of $\gamma_i$.
We will also denote by $\overline{S}$
the surface obtained from $S$ by collapsing each
$\gamma_i$ to a node $\nu_i$.

We recall a construction of an infinite energy harmonic map
$f:(S\setminus\Gamma,c')\rightarrow(S\setminus\Gamma,h)$ by Wolf
\cite{wolf:infinite}.
For every $i=1,\dots,l$, choose a sequence $s_{i,n}\rightarrow+\infty$
with $s_{i,n}>1$.
Let $U_{i,+}(s)$ and $U_{i,-}(s)$ be the cusps of
$(S\setminus\Gamma,c')$ bounded by horocycles of length $1/s$ near
$\nu_i$ and let $U_i(s):=U_{i,+}(s)\cup U_{i,-}(s)$ and
$\overline{U}_i(s):=U_i(s)\cup\gamma_i$.

Fix an isometry $\xi_{i,\bullet}:U_{i,\bullet}(1)\rightarrow
U/(0,y)\sim(1,y)$ with $U=[0,1]\times[1,+\infty)\subset\mathbb{H}$ and
put on $U_{i,\bullet}(1)$ the
flat metric $|\Psi_i|$ induced by the quadratic differential
$\Psi_i=\xi_{i,\bullet}^*(dz^2)$.

Define a new metric $c$ on $S\setminus\Gamma$ which
\begin{itemize}
\item agrees with $c'$ outside $\bigcup_i U_i(1)$
and with $|\Psi_i|$ on each $U_{i,\bullet}(1)$,
\item is conformally equivalent to $c'$,
\item $c$ is smooth away from $\partial U_i(1)$.
\end{itemize}

Let $U_{i,\bullet}^n(s)$ be the annulus $U_{i,\bullet}(s)\setminus
U_{i,\bullet}(s_{i,n})$
and call $U_i^n(s):=U_{i,+}^n(s)\cup U_{i,-}^n(s)$.
We denote by $\xi_{i,\bullet}^n:U_{i,\bullet}^n(1)\rightarrow U/\sim$
the restriction of $\xi_{i,\bullet}$.

Then $S_n$ is obtained from $S$ by
removing $\overline{U}_i(s_{i,n})$
from the cusps adjacent to $\gamma_i$.
Gluing the seams together, we obtain a compact surface $\overline{S}_n$
with quadratic differentials $\Psi_i^n$ on $\overline{U}_i^n(1)$
obtained restricting $\Psi_i$,
with distinguished geodesics $\gamma_i^n$ corresponding to the seams and
collars $\overline{U}_i^n(s)=U_i^n(s)\cup\gamma_i^n$.
We will also define $\xi_i^n:\overline{U}_i^n(1)\rightarrow U/\sim$ as
\[
\xi_i^n(p)=
\begin{cases}
\xi_{i,+}^n(p) & \text{if $p\in\overline{U}_{i,+}^n$} \\
1+2i(s_{i,n}+1)-\xi_{i,-}^n(p) &  \text{if $p\in\overline{U}_{i,-}^n$}
\end{cases}
\]
Notice that the metric $c_n$ induced by $c$ on
$\overline{S}_n$ detemines a point in $\mathcal{T}(S)$:
we will denote by $c'_n$ the hyperbolic metric conformally equivalent to $c_n$.

Notice that, hidden in this construction, there is an arbitrary choice of
twists associated to the gluings or, equivalently, to the charts
$\xi_{i,\bullet}$.

Call $f_n:(\overline{S}_n,c_n)\rightarrow (S,h)$ the unique harmonic
map in the given homotopy class \cite{eells-sampson:harmonic}.

\begin{theorem}[Wolf \cite{wolf:infinite}]
Up to subsequences, $f_n$ converges $C^{2,\alpha}$
to a harmonic map $f:(S\setminus\Gamma,c)
\rightarrow (S\setminus\Gamma,h)$ on the compact subsets
of $S\setminus\Gamma$.

In each cusp $U_{i,\bullet}(1)$ of $(S\setminus\Gamma,c)$,
the Hopf differential $\Phi$ of $f$ 
looks like
$\Phi=(\ell_i^2/4+O(e^{-2\pi y\circ\xi_{i,\bullet}}))\Psi_i$. Moreover,
the energy
density $e(f;c,h)\rightarrow \ell_i^2/2$
and the holomorphic energy density $H(f;c,h)\rightarrow \ell_i^2/4$
as $y\circ\xi_{i,\bullet}\rightarrow+\infty$.
\end{theorem}

We will assume throughout this section that we have already extracted
a good subsequence (which we will still call $f_n$)
so that the above theorem holds.

Now let $h^\dual_n$ be the metric on $S$ antipodal to $h$ with respect to
$(f_n)_* c_n$. Because of the theorem, $h^\dual_n$ converges
to some $h^\dual$ smoothly away from $\Gamma$.

\begin{prop}\label{prop:centers}
Let $1=b_1=\dots=b_r>b_{r+1}\geq\dots\geq b_l>0$ and put
$s_{i,n}=(a_i/\ell_i)t_n^{b_i}$, where $a_i>0$ and $t_n\rightarrow+\infty$.
Then, up to subsequences,
$c_n\rightarrow [b_1\gamma_1+\dots+b_l\gamma_l]$ and
$h_n^\dual\rightarrow [a_1\gamma_1+\dots+a_r\gamma_r]$ in
Thurston compactification
of $\mathcal{T}(S)$. 
\end{prop}

\begin{cor} \label{cor:center}
If $h$ is fixed and $h^\dual_n\rightarrow[\lambda]$ in Thurston
compactification, then $c_n$ does not
necessarily converge to $[\lambda]$.
If $h$ is fixed and $c_n\rightarrow[\lambda]$ in Thurston
compactification, then $h_n^\dual$ does not
necessarily converge to $[\lambda]$ and so the cyclic flow centered at $c_n$
does not necessarily converge to an earthquake along $\lambda$
(with any normalization).
\end{cor}

In order to prove Proposition~\ref{prop:centers},
we need to estimate the {\it{transversal length}}
$\mathrm{trl}_{\gamma_i}(c'_n)$,
that is, the width of a standard $c'_n$-collar of $\gamma_i^n$ bounded
by hypercycles of length $1$.

\begin{lemma}\label{lemma:easy}
The extremal length of $\gamma_i$ at $c_n$ satisfies
\[
\frac{1}{C_1+2(a_i/\ell_i)t_n^{b_i}}\leq
\mathrm{Ext}_{\gamma_i}(c_n)\leq \frac{\ell_i}{2a_it_n^{b_i}}
\]
and so $\mathrm{trl}_{\gamma_i}(c'_n)\asymp 2b_i\log t_n$.
\end{lemma}
\begin{proof}
By construction, $(\overline{U}_i^n(1),c_n)$
contains a flat cylinder
of circumference
$1$ and height $2(a_i/\ell_i) t_n^{b_i}$ and so the extremal length satisfies
\[
\mathrm{Ext}_{\gamma_i}(c_n)\leq \frac{\ell_i}{2a_i t_n^{b_i}}~. 
\]
On the other hand,
%
%
consider a metric $\hat{c}_n$
on $\overline{S}_n$ which is conformally equivalent to $c_n$,
which agrees with $c_n$ on $\overline{S}_n\setminus
\bigcup_{j\neq i}\overline{U}_j^n$ and 
such that the $\hat{c}_n$-area of $\overline{U}_j^n$ is bounded by
a fixed constant for $j\neq i$ and the distance between the two
boundary components of $\overline{U}_j^n$ is at least $1$.
For instance, one can define $\hat{c}_n$ be rescaling
$c_n$ by a factor which is constantly
$1/s_{j,n}$ on the regions
$\overline{U}_j^n(2)$ for all $j\neq i$, which interpolates between
$1$ and $1/s_{j,n}$ on $\overline{U}_j^n(1)\setminus \overline{U}_j^n(2)$
for $j\neq i$, and which is constantly $1$ elsewhere.
 
Then $\ell_{\gamma_i}(\hat{c}_n)=1$
and the area $Area_{\hat{c}_n}(S)\leq C_1+2(a_i/\ell_i)t_n^{b_i}$
where $C_1$ is a constant that depends only on $\chi(S)$ and $k$.
Hence,
\[
\mathrm{Ext}_{\gamma_i}(c_n)\geq \frac{1}{C_1+2(a_i/\ell_i)t_n^{b_i}}
\]
As $\mathrm{Ext}_{\gamma_i}(c_n)\rightarrow 0$, Maskit's estimate
(see \cite{maskit})
gives
$\ell_{\gamma_i}(c_n)\asymp \pi\mathrm{Ext}_{\gamma_i}(c_n)$ and so
$\mathrm{trl}_{\gamma_i}(c'_n)\asymp
-2\log\ell_{\gamma_i}(c'_n) \asymp 2b_i \log t_n$.
\end{proof}

For each $i$, fix
an open neighbourhood $A_i\subset (\overline{S},c)$
of $\nu_i$ whose closure does
not contain any zero of $\Phi$, and such that
$(A_i\setminus\nu_i,|\Phi|)$ is the union
of two annuli.
Moreover, choose standard $h$-collars $R_i\subset (S,h)$ around $\gamma_i$
such that $R_i\subset f(A_i)$.
%
%
%
By Wolf's construction (see \cite{wolf:infinite}),
outside $\bigcup_i f^{-1}(R_i)$
the Hopf differential $\Phi_n$ of $f_n$ converges
$C^{1,\alpha}$ to $\Phi$.

Here we recall that, by definition,
\[
h=2(f_n)_*\mathrm{Re}(\Phi_n)+
e(f_n;c_n,h) \, c_n
\]
where $e(f_n;c_n,h)$ is the energy density; moreover,
by Equation (\ref{eq:hopf}),
\[
(f_n)_*\mathrm{Re}(\Phi_n)=\frac{1}{4}h((E-b_n^2)\cbull,\cbull)
\]
and so $b_n$ converges $C^{1,\alpha}$ to $b$ outside $\bigcup_i R_i$.

Notice that the horizontal (resp. vertical)
directions of $\Phi_n$ are exactly the
eigenspaces of $b_n$ corresponding to the smaller (resp. bigger) eigenvalue.

%
%

\begin{lemma}\label{lemma:stable}
Fix a small $\varepsilon>0$.
Up to shrinking $R_i$ and for $n$ large enough,
\[
\left|\frac{4}{\ell_i^2}\frac{\Phi_n}{\Psi_i^n}-1\right|
<\varepsilon^2
\]
and
\[
\left|\frac{2}{\ell_i^2}e(f_n;c_n,h)-1\right|<2\varepsilon^2
\]
in every $f_n^{-1}(R_i)$.
\end{lemma}

\begin{proof}
Up to shrinking $R_i$, we can assume that
$R_i\subset f_n(\overline{U}_i^n(3))$ and there
\[
\left|\frac{4}{\ell_i^2}\frac{\Phi}{\Psi_i}-1\right|
<\varepsilon^2/2
\quad\text{and}\quad
\left|\frac{4}{\ell_i^2}H(f;c,h)-1\right|<\varepsilon^2/2
\]
where $H(f;c,h)=\frac{1}{2}\|\partial f\|^2$
is the holomorphic energy density of $f$.

As $f^{-1}_n|_{\partial R_i}\rightarrow f^{-1}|_{\partial R_i}$,
for $n$ large enough $f^{-1}_n(\partial R_i)$ is contained inside
$\overline{U}_i^n(2)$.
Moreover,
$(\xi^n_{i,\bullet})_*\Phi_n\rightarrow
(\xi_{i,\bullet})_*\Phi$ in a compact neighbourhood $K_{i,\bullet}$ of
$\xi^n_{i,\bullet}(f_n^{-1}(\partial R_i)\cap U_{i,\bullet}^n)$ and so
\[
\frac{4}{\ell_i^2}
\left|(\xi_{i,\bullet}^n)_*\frac{\Phi_n}{\Psi_i^n}-
(\xi_{i,\bullet})_*\frac{\Phi}{\Psi_i} \right|<\varepsilon^2/2
\]
in $K_{i,\bullet}$. Thus,
\[
\frac{4}{\ell_i^2}
\left|(\xi_i^n)_*\frac{\Phi_n}{\Psi_i^n}-1 \right|<\varepsilon^2
\]
in $K_{i,\bullet}$. Because $\xi_{i,\bullet}^n$ is holomorphic and
$\Phi_n/\Psi_i^n$ is a holomorphic function on $f_n^{-1}(R_i)$,
the same estimate holds in $f_n^{-1}(R_i)$ for $n$ large enough.

In a similar way, if $H_n=H(f_n;c_n,h)$, then
the energy density and the Jacobian of $f_n$ are
\[
e(f_n;c_n,h)=H_n+\frac{|\Phi_n|^2}{|\Psi_i^n|^2 H_n} \quad\text{and}\quad
J(f_n;c_n,h)=H_n-\frac{|\Phi_n|^2}{|\Psi_i^n|^2 H_n}>0
\]
which implies $H_n> |\Phi_n|/|\Psi_i^n|>(\ell_i^2/4)(1-\varepsilon^2)$
on $f_n^{-1}(R_i)$.

On the other hand (see \cite{wolf:teichmuller}, for instance)
\[
\Delta_{c_n}\log H_n=2J_n> 0
\]
and so $\Delta_{c_n}\log(4H_n/\ell_i^2)\geq 0$ on $\overline{U}_i^n(2)$.
As $\log(4H_n/\ell_i^2)<\varepsilon^2/2$ on $\partial f_n^{-1}(R_i)$
for $n$ large, we obtain
$H_n\leq(\ell_i^2/4)(1+\varepsilon^2/2)$ on $f_n^{-1}(R_i)$
and so the wished estimate for $e(f_n;c_n,h)$.
\end{proof}

Thanks to the previous lemma, we can draw a few consequences.

\begin{cor}\label{eigen:cr}
\begin{itemize}
\item[(i)]
The metric $h^\dual$ has nodes at $\Gamma$ and $h_n^\dual\rightarrow h^\dual$ in
the augmented Teichm\"uller space.
\item[(ii)]
In the whole $f_n^{-1}(R_i)$ found in the above lemma,
the bigger eigenvalue $\kappa$ of $b_n$ is greater than $1/2\varepsilon$
for large $n$.
\end{itemize}
\end{cor}

\begin{proof}
As for (i), notice that $h^\dual_n=-2\mathrm{Re}(\Phi_n)+e(f_n;c_n,h)c_n$.
Because of Lemma~\ref{lemma:stable}, for every $\varepsilon>0$
there exists $R_i$ such that, for $n>n(\varepsilon)$,
the $h^\dual_n$-norm of
$(f_n\circ (\xi_{i}^n)^{-1})_*\partial_x$ is at most $2\varepsilon\ell_i$
the $h_n$-norm of
$(f_n\circ(\xi_i^n)^{-1})_*\partial_y$ is at most $2\varepsilon\ell_i$.
Thus, $\ell_{\gamma_i}(h^\dual_n)\leq 2\varepsilon\ell_i$ and so
$\ell_{\gamma_i}(h^\dual_n)\rightarrow 0$.
As $h^\dual_n\rightarrow h^\dual$ on the compact subsets of $S\setminus \Gamma$,
we conclude that $h^\dual$ has nodes at $\gamma_i$ and the sequence converges
in the augmented Teichm\"uller space.

As for (ii), the two identities $f_n^* h=2\mathrm{Re}(\Phi_n)+e(f_n;c_n,h)c_n$
and $f_n^*h((E-b_n^2)\cbull,\cbull)=4\mathrm{Re}(\Phi_n)$,
already seen in Section \ref{section:definition} Equation (\ref{eq:hopf}),
give
\[
\frac{h((E-b_n^2)(f_n)_* v,(f_n)_* v)}{h((f_n)_*v,(f_n)_*v)}=
\frac{4\mathrm{Re}(\Phi_n)(v,v)}{2\mathrm{Re}(\Phi_n)(v,v)+e(f_n;c_n,h)c_n(v,v)}
\]
for any vector tangent $v$ to $\overline{U}_i^n(1)$.
Choosing $v=(\xi_{i}^n)^*\partial_y$, we obtain
\[
\frac{h((E-b_n^2)(f_n)_* v,(f_n)_* v)}{h((f_n)_*v,(f_n)_*v)}\leq
\frac{\ell_i^2(-1+3\varepsilon^2)}
{4\ell_i^2\varepsilon^2}
\]
on $\overline{U}_i^n(1)\cap f_n^{-1}(R_i)$ 
and so $1-\kappa^2\leq (-1/4\varepsilon^2)+1$ there, that is
$\kappa\geq 1/2\varepsilon$.
%
%
\end{proof}

Up to subsequences, we can assume that $h_n^\dual$ converges to
a point $[\lambda]$ in Thurston boundary. Notice that $\lambda$
must be supported on $\gamma_1\cup\dots\cup\gamma_l$ and so
$\lambda=w_1\gamma_1+\dots+w_l\gamma_l$ with $w_1,\dots,w_l\geq 0$.

We will show that $[\lambda]=[a_1\gamma_1+\dots+a_r\gamma_r]$,
and so the result
will not depend on the chosen subsequence.

Let $(\theta_n)_{n\in\N}$ be a sequence of positive numbers such that
$\theta_n\ell_{h_n^\dual}\rightarrow \iota(\lambda,\cbull)$.

\begin{lemma}\label{weight:lm}
For every $i$, we have $\displaystyle
\theta_n\ell_i s_{i,n}\rightarrow w_i
$ as $n\rightarrow\infty$.
\end{lemma}
\begin{proof}
Let $\varepsilon>0$. We can choose a collar
$\gamma_i\subset R_i\subset S$ such that
$\tr(b_n)>2/\varepsilon$ on $R_i$ and
$|\frac{4}{\ell_i^2}\Phi_n-\Psi_i^n|<\varepsilon^2|\Psi_i^n|$ on
$f_n^{-1}(R_i)$ for $n$ large.

As $f^{-1}_n\rightarrow f^{-1}$ on $S\setminus\Gamma$,
we can assume that there exists $\bar{y}>1$ such that
$\overline{U}^i_n(1)\supset f_n^{-1}(R_i)
\supset \overline{U}^i_n(\bar{y})$ for $n$ large.

By abuse of notation, denote just
by $(x,y)$ the Euclidean coordinates
on $R_i$ given by the parametrization $\xi_i^n\circ f_n$.

By writing the relation $4Re(\Phi_n)=h((E-b_n^2)\cbull,\cbull)$
in coordinates on $R_i$ and taking the determinant, we obtain
\begin{equation}\label{xy:eq}
(1+\eta(\varepsilon))\ell_i^2
\ dx\wedge dy=\sqrt{\tr(b_n)^2-4}\ \omega_h
\end{equation}
where $|\eta(\varepsilon)|<\varepsilon^2$ for $n$ large, by
Lemma \ref{lemma:stable}.

By Corollary \ref{eigen:cr}(ii), $(1-\varepsilon^2/2)\tr(b_n)
\leq \sqrt{\tr(b_n)^2-4}\leq \tr(b_n)$ for $n$ large.

Multiplying by $\theta_n$ both hand sides of (\ref{xy:eq})
and integrating over $R_i$, we get
\begin{equation}\label{pippo:eq}
\theta_n(1-\varepsilon^2/2)
\int_{R_i}\tr(b_n)\omega_h\leq
\theta_n(1+\eta(\varepsilon))\ell_i^2\int_{R_i}dx\wedge dy\leq
\theta_n\int_{R_i}\tr(b_n)\omega_h
\end{equation}
Now,
\begin{equation}\label{pluto:eq}
\theta_n\ell_i^2(s_{i,n}-\bar{y})
\leq \theta_n\ell_i^2\int_{R_i}dx\wedge dy
\leq \theta_n \ell_i^2 s_{i,n}
\end{equation}
and $\theta_n\int_{R_i}\tr(b_n)\omega_h
\rightarrow w_i\ell_i$ by Proposition \ref{tr:prop}.

The result follows by comparing Equations (\ref{pippo:eq}) and (\ref{pluto:eq}).
\end{proof}

\begin{proof}[Proof of Proposition~\ref{prop:centers}]
By Lemma~\ref{lemma:easy} and Maskit's estimate, we immediately
obtain that $c_n\rightarrow[b_1 \gamma_1+\dots+b_l\gamma_l]$.

By Lemma \ref{weight:lm},
\[
\frac{w_i}{w_j}=\lim_n \frac{\ell_i s_{i,n}}{\ell_j s_{j,n}}=
\lim_n \frac{a_i t_n^{b_i}}{a_j t_n^{b_j}}
\]
which shows that $[\lambda]=[a_1\gamma_1+\dots+a_r\gamma_r]$.
\end{proof}


\section{The landslide flow on the universal Teichm\"uller space}
\label{sc:universal}

In this section we show how the construction of the landslide flow $L$ extends to
the universal Teichm\"uller space. We believe that this $S^1$ action on the product
of two copies of the universal Teichm\"uller space can be of independent
interest, but limit our investigations here to its definition and to checking 
that it is non-trivial.

\subsection{Minimal Lagrangian maps and the universal Teichm\"uller space}

The universal Teichm\"uller space $\cT_U$ can be defined as the quotient of the group
$\cQS$ of quasi-symmetric homeomorphisms of the circle by composition on the left 
by projective transformations, see e.g. \cite{gardiner-harvey}.
We will show here that the map $L$ defined above extends to a circle action
$\mathcal{L}$ on $\cT_U\times \cT_U$.
This is based on the following statement.

\begin{theorem}[\cite{maximal}]
Let $\bar\psi\in \cQS$. There exists a unique quasi-conformal minimal Lagrangian
diffeomorphism $m:\Hyp^2\rightarrow \Hyp^2$ such that $\partial m=\bar\psi$.
\end{theorem}

We call $g$ the hyperbolic metric on $\Hyp^2$, and $\nabla$ its Levi-Civita
connection. 
It follows from the basic facts on minimal Lagrangian diffeomorphisms, as
recalled in Section \ref{ssc:minilag}, that there exists a unique bundle morphism
$b:T\Hyp^2\rightarrow T\Hyp^2$ such that:
\begin{itemize}
\item $m^*g=g(b\bullet, b\bullet)$,
\item $\det(b)=1$,
\item $b$ is self-adjoint for $g$,
\item $b$ satisfies the Codazzi equation: $d^\nabla b=0$.
\end{itemize}
Since $m$ is quasi-conformal, $b$ has eigenvalues in $[\epsilon,
1/\epsilon]$ for some $\epsilon>0$, see \cite{maximal}.

\subsection{An extension of $L$ to the universal Teichm\"uller space}

The construction of the $S^1$-action is based on the following definition and lemma.

\begin{defi}
Let $\bar\psi\in \cQS$ and let
$e^{i\theta}\in S^1$. Let $m:\Hyp^2\rightarrow \Hyp^2$ be the
unique quasi-conformal minimal Lagrangian diffeomorphism such that
$\dr m=\bar\psi$,
and let $b:T\Hyp^2\rightarrow T\Hyp^2$ be as above.
We call 
$$ \beta_\theta = \cos(\theta/2)E+\sin(\theta/2)Jb $$
and set $g_\theta:=g(\beta_\theta \bullet,\beta_\theta \bullet)$ and $g^\dual_\theta=
g(\beta_{\theta+\pi} \bullet,\beta_{\theta+\pi} \bullet)$.
\end{defi}

\begin{lemma} \label{lm:quasi}
With the notations above, $g_\theta$ and $g^\dual_\theta$ are complete hyperbolic metrics on $\Hyp^2$. 
The identity map from $(\Hyp^2,g_\theta)$ to $(\Hyp^2,g^\dual_\theta)$ is minimal Lagrangian and quasi-conformal.  
\end{lemma}

\begin{proof}[Sketch of the proof]
The fact that $g_\theta$ and $g^\dual_\theta$ have curvature $-1$
follows from the same argument as in the
proof of Proposition \ref{pr:basics}; we do not repeat it here.
Moreover, the argument given in the proof of Theorem \ref{thm:c}
also shows that the identity map from
$(\Hyp^2,g_\theta)$ to $(\Hyp^2,g^\dual_\theta)$ is
minimal Lagrangian.

To check that the identity between $(\Hyp^2,g)$ and $(\Hyp^2,g_\theta)$ is quasiconformal, it is sufficient to prove that 
the eigenvalues of the bundle morphism ${}^t\beta_{\theta}\cdot\beta_{\theta}$ are between $\epsilon'$
and $1/\epsilon'$, for some $\epsilon'>0$ depending on $\epsilon$. However $\det(\beta_{\theta})=1$
by definition, so that ${}^t\beta_{\theta}\cdot\beta_{\theta}$ has determinant $1$. To compute its 
trace, notice that 
$$ {}^t\beta_{\theta}\cdot\beta_{\theta} = \cos^2(\theta/2)E+\sin^2(\theta/2)b^2 + \cos(\theta/2)\sin(\theta/2)(Jb-bJ) $$
and that $\tr(Jb)=\tr(bJ)=0$. It follows that 
$$ \tr({}^t\beta_{\theta}\cdot\beta_{\theta}) = 2\cos^2(\theta/2)+\sin^2(\theta/2)\tr(b^2)~. $$
Since the eigenvalues of $b$ are in $[\epsilon, 1/\epsilon]$, it follows that the identity between
$(\Hyp^2,g)$ and $(\Hyp^2,g_\theta)$ is quasi-conformal, and that $g_\theta$ and $g_\theta^\dual$ are complete.
\end{proof}

As a consequence we can give the definition of the action considered here. 
Let $\psi,\psi^\dual\in\cQS$ represent points $[\psi],[\psi^\dual]\in\cT_U$.
Then $\bar\psi:=\psi^\dual\circ \psi^{-1}:S^1\rightarrow S^1$
is a quasisymmetric homeomorphism:
let $m:\Hyp^2\rightarrow \Hyp^2$ be the unique quasiconformal
minimal Lagrangian diffeomorphism with
$\partial m=\bar\psi$. We can then define a
bundle morphism $b:T\Hyp^2\rightarrow T\Hyp^2$ as above,
as well as two hyperbolic metrics $g_\theta$
and $g^\dual_\theta$ associated to $\theta$ and $b$,
as in Section 3. Lemma \ref{lm:quasi} shows that
the identity between $(\Hyp^2, g_\theta)$ and $(\Hyp^2,g^\dual_\theta)$
is minimal Lagrangian.

Since $g_\theta$ is hyperbolic, $(\Hyp^2,g_\theta)$ is isometric
to the hyperbolic plane. The 
identity between $(\Hyp^2, g)$ and $(\Hyp^2,g_\theta)$ therefore
determines a quasiconformal diffeomorphism
$\Psi_\theta$ between $(\Hyp^2,g)$ and $\Hyp^2$, well-defined up to composition on the left, with boundary value $\psi_\theta\in\cQS$.
Similary, the identity map from $(\Hyp^2,g)$ to $(\Hyp^2,g^\dual_\theta)$
determines a quasiconformal map $\Psi_\theta^\dual$ 
between the hyperbolic plane and itself,
with boundary value $\psi_\theta^\dual\in\cQS$.
Then $\Psi_\theta^\dual\circ \Psi_\theta^{-1}$
is a quasiconformal minimal Lagrangian
diffeomorphism by Lemma \ref{lm:quasi}. 
We define $\cL$ as:
$$ \cL_{e^{i\theta}}([\psi],[\psi^\dual])=([\psi_\theta],[\psi_\theta^\dual])~. $$

To establish a relation between $\cT_U$ and $\cT_S$,
fix a hyperbolic metric $h_0$ on $S$ and a universal covering map
$\Hyp^2\cong(\tilde{S},\tilde{h}_0)\rightarrow (S,h_0)$,
so that $\rho_0:\pi_1(S)
\rightarrow\PSL_2(\R)$ is the associated holonomy representation.
Given $[h]\in\cT_S$, we can consider the lift $\tilde{h}$ of $h$
to $\Hyp^2$ and let $\rho$ be its holonomy representation.
The identity map from $(\Hyp^2,\tilde{h}_0)$ to 
$(\Hyp^2,\tilde{h})$ determines a quasi-conformal diffeomorphism
$\Psi$ of $\Hyp^2$ to itself, with boundary value $\psi\in\cQS$,
that conjugates the action of $\rho_0$ on $\Hyp^2$ to the 
action of $\rho$.

Let $i_S:\cT_S\hookrightarrow\cT_U$ be the canonical embedding
of Teichm\"uller space of $S$ in the universal Teichm\"uller space
defined as $i_S([h])=[\psi]$.

\begin{prop} \label{pr:finite}
The restriction via $i_S$ of $\cL$ to $\cT_S\times \cT_S\subset
\cT_U\times \cT_U$ is the landslide action $L$.
\end{prop}

\begin{proof}
Let $[\psi]=i_S([h])$ and $[\psi^\dual]=i_S([h^\dual])$ be points of
$i_S(\cT_S)\subset \cT_U$. Let $\rho$ and $\rho^\dual$ be the holonomy
representations of $h$ and $h^\dual$ respectively.
Then $\psi$ and $\psi^\dual$ are the boundary values of
quasiconformal maps $\Psi,\Psi^\dual:\Hyp^2\rightarrow \Hyp^2$
which are conjugating $\rho_0$
to actions $\rho$ and $\rho^\dual$ on $\Hyp^2$.
By construction, $g_\theta=\tilde{h}_\theta$ is $\rho$-invariant and
$g_\theta^\dual=\tilde{h}^\dual_\theta$ is $\rho^\dual$-invariant;
moreover, $\Psi_\theta$ (resp. $\Psi_\theta^\dual$) conjugates the
action of $\rho_0$ to the action of the holonomy
representation $\rho_\theta$ of $h_\theta$
(resp. the holonomy representation $\rho_\theta^\dual$ of $h_\theta^\dual$).
Hence,
$([\psi_\theta],[\psi_\theta^\dual])\in
i_S(\cT_S)\times i_S(\cT_S)\subset \cT_U\times \cT_U$,
and the restriction of $\cL$ to $\cT_S\times \cT_S$ coincides with $L$,
as claimed.
\end{proof}

\begin{theorem} \label{tm:univ}
The map $\cL$ defines a non-trivial action of $S^1$ on $\cT_U$. 
\end{theorem}

\begin{proof}[Sketch of the proof]
To prove that $\cL$ determines an action of $S^1$, it is sufficient to check that,
for all $\theta,\theta'\in \R$ and all $[\psi],[\psi^\dual] \in \cT_U$, 
$\cL_{e^{i\theta}}(\cL_{e^{i\theta'}}([\psi],[\psi^\dual]))=\cL_{e^{i(\theta+\theta')}}([\psi],[\psi^\dual])$.

However, this follows from the fact that $\Psi^\dual_\theta\circ \Psi_\theta^{-1}$ is
minimal Lagrangian, followed by the same argument used
in the proof of Theorem \ref{tm:cyclic}; so we do not repeat them here.

The nontriviality of $\cL$ is clear, since $L$ is non-trivial in all the copies
of Teichm\"uller spaces of surfaces of finite genus, see Proposition \ref{pr:finite}.
\end{proof}

\section{Applications, extensions, and questions}

This section contains a brief outline of some possible applications of the 
landslide flow developed here, and of some open questions.

\subsection{Constant curvature surfaces in globally hyperbolic AdS manifolds}

As already mentioned, a recent result of Barbot, B\'eguin and Zeghib \cite{BBZ2}
states that
the complement of the convex core in a globally hyperbolic AdS manifold has a unique
foliation by constant Gauss curvature surfaces. This was used above in the proof
of Theorem \ref{tm:earthquake}. 

However we believe that this result might not be necessary to prove Theorem \ref{tm:earthquake},
and that a direct, albeit longer, proof could be given, based on a deformation argument. Using
the same argument as in the proof of Theorem \ref{tm:earthquake}, but backwards, it should then
be possible to recover parts of the main result of \cite{BBZ2}: the proof of the existence, in
a given globally hyperbolic AdS 3-dimensional manifold, of a unique
surface of prescribed constant curvature.

\subsection{Holomorphic disks in Teichm\"uller space}

One obvious consequence of Theorems \ref{tm:earthquake} and \ref{tm:complex} is the existence
of many holomorphic disks in Teichm\"uller space of $S$:
given $h,h'\in \cT$ with $h\neq h'$ and given $\zeta\in S^1\setminus \{ 1\}$
there exists $h^\dual\in\cT$
and a holomorphic map
$C_\cbull(h,h^\dual):\overline{\Delta}\rightarrow \cT$
from the unit disk in $\C$ to
Teichm\"uller space of $S$, such that $C_1(h,h^\dual)=h$
and $C_\zeta(h,h^\dual)=h'$. 

It is of course conceivable that the disks obtained in this manner for two different
values of $h'$ have the same image. However, there are reasons to believe that it is
not often the case. If this is correct, it would mean that the landslide disks provide
a $(12g-11)$-dimensional family of holomorphic disks in $\cT$. 


\subsection{Other questions}

There are many remaining questions concerning the landslide flow or its complex extension,
mostly motivated by the analogy with the earthquake flow. Some of those statements can be
translated in terms of 3-dimensional hyperbolic or AdS geometry. We give here a short list
of example of possible questions. 

\subsubsection*{Smooth grafting as homeomorphism}

Recall that Scannell and Wolf \cite{scannell-wolf} proved that, for $\lambda\in \cML$ fixed, 
the map $h\mapsto gr_\lambda (h)$ is a homeomorphism of $\cT$. When $h\in \cT$ is fixed, the
map $\lambda\mapsto gr_\lambda(h)$ is also a homeomorphism from $\cML$ to $\cT$, see 
\cite{dumas-wolf}. 

It is tempting to ask whether those statements can be extended to the smooth grafting map $sgr$. Note that 
in this setting the two statements above concerning the grafting map -- with the 
measured lamination fixed, and with the hyperbolic metric fixed -- are now merged into one,
since the two hyperbolic metrics that occur in the map $sgr$ play symmetric roles.

\begin{question} \label{q:sgr}
Let $s\in (0,1)$, and let $h\in \cT$. Is the map
$h^\dual\mapsto sgr_s(h,h^\dual)$ a 
homeomorphism from $\cT$ to $\cT$?
\end{question}

This statement can be translated in terms of the geometry of hyperbolic ends, in the
following, essentially equivalent question. 

\begin{question} \label{q:sgr2}
Let $h,c\in \cT$ and let $K\in (-1,0)$. 
\begin{itemize}
\item Is there a unique hyperbolic end with conformal
structure at infinity $c$, and containing an embedded surface of constant curvature $K$
with induced metric proportional to $h$? 
\item Is there a unique hyperbolic end with conformal structure at infinity $c$, containing
an embedded surface of constant curvature $K$ with third fundamental form proportional to
$h^\dual$?
\end{itemize}
\end{question}

\subsubsection*{The action of the landslide flow at infinity}

It is quite natural to wonder to what extend the landslide flow can be extended
to Thurston boundary of Teichm\"uller space. One side of this question is
already answered above in Section 6, concerning the limit of $L$ to the earthquake
flow when one of the parameter converges to Thurston boundary and the other
is fixed. However other questions can be asked, in particular in light of the
results of Wolf \cite{wolf:infinite} on the behavior of harmonic maps at
the boudary of $\cT$.

\subsubsection*{The landslide flow as a Hamiltonian flow}

Consider a fixed measured lamination $\lambda\in \cML$. 
The flow of earthquakes along $\lambda$ is the Hamiltonian flow of the length of $\lambda$,
considered as a function on $\cT$,
with respect to the Weil-Petersson symplectic structure.
In a similar way, is the landslide flow the Hamiltonian
flow of some functional on $\cT\times \cT$?

\subsubsection*{The data at infinity of hyperbolic ends}

For all $K\in (-1,0)$, there is a parameterization of $\cCP$ by
$\cT\times \cT$, with a complex projective structure $P$ corresponding to
$(h,h^\dual)$ if the hyperbolic end $E$ with complex projective structure
$P$ at
infinity contains a surface of constant curvature $K$ with induced metric
proportional to $h$ and third fundamental form proportional to $h^\dual$.

There is also another parameterization of $\cCP$ by the space of couples
$(h,b)$, where $h\in \cT$ and where $b$
is a bundle morphism which is self-adjoint
for $h$ and satisfies the Codazzi equation and $\det(b)=1$. 

Given a $\mathbb{CP}^1$-structure $P$,
we can also consider the data at infinity
$I^*$ and $\III^*$ of the
corresponding hyperbolic end, as defined in \cite{volume}, 
and take the limit as $K\rightarrow 0$. Is it true that
$h$ and $h^\dual$ limit to $I^*$ and $\III^*$ respectively?
And that the traceless part of 
$B$, suitably renormalized, converges to $B^*$?

\subsubsection*{Landslides on the universal Teichm\"uller space}

Section 8 on the universal Teichm\"uller space leaves a number of elementary
questions unanswered. One natural question is whether for fixed
$[\psi],[\psi^\dual]\in \cT_U$
the map $e^{i\theta}\rightarrow \cL^1_{e^i\theta}([\psi],[\psi^\dual])$ extends to a 
holomorphic disk in $\cT_U$,
as for the landslide action on the Teichm\"uller
space of a closed surface.

Another natural question is whether all fixed points of the
landslide action on $\cT_U\times \cT_U$ are on the diagonal.

\subsubsection*{Cone singularities}

It appears possible that all the results obtained here extend from closed 
hyperbolic surfaces to finite volume hyperbolic surfaces, and more generally to
hyperbolic surfaces with cone singularities (perhaps of angle less than $\pi$).
The 3-dimensional AdS or hyperbolic part of the picture would then be filled
with 3-dimensional AdS or hyperbolic manifolds with ``particles'', as considered
e.g. in \cite{minsurf,qfmp,conebend}.



\bibliographystyle{amsplain}
\bibliography{/home/schlenker/papiers/outils/biblio}

\def\cprime{$'$}
\providecommand{\bysame}{\leavevmode\hbox to3em{\hrulefill}\thinspace}
\providecommand{\MR}{\relax\ifhmode\unskip\space\fi MR }
\providecommand{\MRhref}[2]{%
  \href{http://www.ams.org/mathscinet-getitem?mr=#1}{#2}
}
\providecommand{\href}[2]{#2}
\begin{thebibliography}{10}


\bibitem{mess-notes}
Lars Andersson, Thierry Barbot, Riccardo Benedetti, Francesco Bonsante,
  William~M. Goldman, Fran{\c{c}}ois Labourie, Kevin~P. Scannell, and Jean-Marc
  Schlenker, \emph{Notes on: ``{L}orentz spacetimes of constant curvature''
  [{G}eom. {D}edicata {\bf 126} (2007), 3--45; mr2328921] by {G}. {M}ess},
  Geom. Dedicata \textbf{126} (2007), 47--70. \MR{MR2328922}

\bibitem{BBZ}
Thierry Barbot, Fran{\c{c}}ois B{\'e}guin, and Abdelghani Zeghib,
  \emph{Constant mean curvature foliations of globally hyperbolic spacetimes
  locally modelled on {${\rm AdS}\sb 3$}}, Geom. Dedicata \textbf{126} (2007),
  71--129. \MR{MR2328923 (2008j:53041)}

\bibitem{BBZ2}
Thierry Barbot, Francois Beguin, and Abdelghani Zeghib, \emph{Prescribing
  {Gauss} curvature of surfaces in 3-dimensional spacetimes, application to the
  {Minkowski} problem in the {Minkowski} space}, 2008, arXiv.org:0804.1053.
  To appear, {\it Ann. Institut Fourier}.

\bibitem{mehdi}
Mehdi Belraouti, PhD thesis, in preparation.

\bibitem{BeBo}
Riccardo Benedetti and Francesco Bonsante, \emph{Canonical {Wick} rotations in
  3-dimensional gravity}, Memoirs of the American Mathematical Society
  \textbf{198} (2009), 164pp, math.DG/0508485.

\bibitem{beguad}
Riccardo Benedetti and Enore Guadagnini, \emph{Cosmological time in
  {$(2+1)$}-gravity}, Nuclear Phys. B \textbf{613} (2001), no.~1-2, 330--352.
  \MR{MR1857817 (2002i:83047)}

\bibitem{maximal}
Francesco Bonsante and Jean-Marc Schlenker, \emph{Maximal surfaces and the
  universal {T}eichm\"uller space}, Invent. Math. \textbf{182} (2010), no.~2,
  279--333. \MR{2729269}

\bibitem{dumas-wolf}
David Dumas and Michael Wolf, \emph{Projective structures, grafting and
  measured laminations}, Geom. Topol. \textbf{12} (2008), no.~1, 351--386.
  \MR{MR2390348 (2009c:30114)}

\bibitem{eells-sampson:harmonic}
James Eells, Jr. and J.~H. Sampson, \emph{Harmonic mappings of {R}iemannian
  manifolds}, Amer. J. Math. \textbf{86} (1964), 109--160. \MR{0164306 (29
  \#1603)}

\bibitem{gardiner-harvey}
Frederick~P. Gardiner and William~J. Harvey, \emph{Universal {T}eichm\"uller
  space}, Handbook of complex analysis: geometric function theory, {V}ol.\ 1,
  North-Holland, Amsterdam, 2002, pp.~457--492. \MR{1966201 (2004a:30041)}

\bibitem{kerckhoff}
Steven~P. Kerckhoff, \emph{The {N}ielsen realization problem}, Ann. of Math.
  (2) \textbf{117} (1983), no.~2, 235--265. \MR{MR690845 (85e:32029)}

\bibitem{minsurf}
Kirill Krasnov and Jean-Marc Schlenker, \emph{Minimal surfaces and particles in
  3-manifolds}, Geom. Dedicata \textbf{126} (2007), 187--254. \MR{MR2328927}

\bibitem{volume}
\bysame, \emph{On the renormalized volume of hyperbolic 3-manifolds}, Comm.
  Math. Phys. \textbf{279} (2008), no.~3, 637--668. \MR{MR2386723}

\bibitem{kulkarni-pinkall}
Ravi~S. Kulkarni and Ulrich Pinkall, \emph{A canonical metric for {M}\"obius
  structures and its applications}, Math. Z. \textbf{216} (1994), no.~1,
  89--129.

\bibitem{LL}
F.~Labourie, \emph{Probl\`eme de Minkowski et surfaces \`a
courboure constante dans les vari\'et\'es hyperboliques}, Bull. Soc.
Math. Fr. \textbf{119}(1991), 307--325

\bibitem{L5}
F.~Labourie, \emph{Surfaces convexes dans l'espace hyperbolique et
  {CP1}-structures}, J. London Math. Soc., II. Ser. \textbf{45} (1992),
  549--565.

\bibitem{conebend}
Cyril Lecuire and Jean-Marc Schlenker, \emph{The convex core of quasifuchsian
  manifolds with particles}, arXiv:0909.4182, 2009.

\bibitem{maskit}
Bernard Maskit, \emph{Comparison of hyperbolic and extremal lengths},
Ann. Acad. Sei. Fenn. Ser. A I Math. {\textbf{10}} (1985), 381--386.
\MR{87c:30062}

\bibitem{mcmullen:complex}
Curtis~T. McMullen, \emph{Complex earthquakes and {T}eichm\"uller theory}, J.
  Amer. Math. Soc. \textbf{11} (1998), no.~2, 283--320. \MR{1478844
  (98i:32030)}

\bibitem{mess}
Geoffrey Mess, \emph{Lorentz spacetimes of constant curvature}, Geom. Dedicata
  \textbf{126} (2007), 3--45. \MR{MR2328921}

\bibitem{mondello}
Gabriele Mondello, \emph{Flows of $\mathrm{SL}_2(\mathbb{R})$-type on
the cotangent space to Teichm\"uller space}, in preparation.

\bibitem{qfmp}
Sergiu Moroianu and Jean-Marc Schlenker, \emph{Quasi-{F}uchsian manifolds with
  particles}, J. Differential Geom. \textbf{83} (2009), no.~1, 75--129.
  \MR{MR2545031}

\bibitem{RH}
Igor Rivin and Craig~D. Hodgson, \emph{A characterization of compact convex
  polyhedra in hyperbolic 3-space}, Invent. Math. \textbf{111} (1993), 77--111.

\bibitem{sampson}
J.~H. Sampson, \emph{Some properties and applications of harmonic mappings},
  Ann. Sci. \'Ecole Norm. Sup. (4) \textbf{11} (1978), no.~2, 211--228.
  \MR{510549 (80b:58031)}

\bibitem{scannell}
Kevin~P. Scannell, \emph{Flat conformal structures and the classification of de
  {S}itter manifolds}, Comm. Anal. Geom. \textbf{7} (1999), no.~2, 325--345.
  \MR{1685590 (2000a:53122)}

\bibitem{scannell-wolf}
Kevin~P. Scannell and Michael Wolf, \emph{The grafting map of {T}eichm\"uller
  space}, J. Amer. Math. Soc. \textbf{15} (2002), no.~4, 893--927 (electronic).
  \MR{MR1915822 (2003d:32011)}

\bibitem{shu}
Jean-Marc Schlenker, \emph{M\'etriques sur les poly\`edres hyperboliques
  convexes}, J. Differential Geom. \textbf{48} (1998), no.~2, 323--405.
  \MR{MR1630178 (2000a:52018)}

\bibitem{horo}
\bysame, \emph{Hypersurfaces in {$H\sp n$} and the space of its horospheres},
  Geom. Funct. Anal. \textbf{12} (2002), no.~2, 395--435. \MR{MR1911666
  (2003d:53108)}

\bibitem{schoen-Yau:78}
Richard Schoen and Shing~Tung Yau, \emph{On univalent harmonic maps between
  surfaces}, Invent. Math. \textbf{44} (1978), no.~3, 265--278. \MR{0478219 (57
  \#17706)}

\bibitem{schoen:role}
Richard~M. Schoen, \emph{The role of harmonic mappings in rigidity and
  deformation problems}, Complex geometry ({O}saka, 1990), Lecture Notes in
  Pure and Appl. Math., vol. 143, Dekker, New York, 1993, pp.~179--200.
  \MR{MR1201611 (94g:58055)}

\bibitem{spivak}
M.~Spivak, \emph{A comprehensive introduction to geometry, vol.i-v}, Publish or
  perish, 1970-1975.

\bibitem{thurston-notes}
William~P. Thurston, \emph{Three-dimensional geometry and topology.},
  Originally notes of lectures at Princeton University, 1979. Recent version
  available on http://www.msri.org/publications/books/gt3m/, 1980.

\bibitem{wolf:teichmuller}
Michael Wolf, \emph{The {T}eichm\"uller theory of harmonic maps}, J.
  Differential Geom. \textbf{29} (1989), no.~2, 449--479. \MR{982185
  (90h:58023)}

\bibitem{wolf:infinite}
\bysame, \emph{Infinite energy harmonic maps and degeneration of hyperbolic
  surfaces in moduli space}, J. Differential Geom. \textbf{33} (1991), no.~2,
  487--539. \MR{1094467 (92b:58055)}

\end{thebibliography}

\end{document}